\documentclass{amsart}

\usepackage{amssymb}

\usepackage{stmaryrd}
\usepackage{multirow}
\usepackage{mathtools}

\newtheorem{theorem}{Theorem}[section]
\newtheorem{lem}[theorem]{Lemma}
\newtheorem{prop}[theorem]{Proposition}
\newtheorem{coro}[theorem]{Corollary}

\newcounter{NoMain}

\newtheorem{mainthm}[NoMain]{Theorem}

\theoremstyle{definition}
\newtheorem{definition}[theorem]{Definition}

\newtheorem{Conj}[theorem]{Conjecture}

\theoremstyle{remark}
\newtheorem{remark}[theorem]{Remark}

\usepackage[bookmarks=true]{hyperref}

\numberwithin{equation}{section}

\newcommand{\R}{\mathbb{R}}
\newcommand{\N}{\mathbb{N}}
\newcommand{\Z}{\mathbb{Z}}

\newcommand{\C}{\mathbb{C}}

\newcommand{\fonction}[5]{#1: \begin{array}{ccc}
#2 & \rightarrow & #3 \\
 #4 & \longmapsto & #5 \end{array}}

\DeclareMathOperator{\ad}{ad}
\DeclareMathOperator{\Ad}{Ad}
\DeclareMathOperator{\vspan}{span}
\DeclareMathOperator{\Id}{Id}
\DeclareMathOperator{\rank}{rank}
\DeclareMathOperator{\Aut}{Aut}
\DeclareMathOperator{\Stab}{Stab}
\DeclareMathOperator{\Hom}{Hom}
\DeclareMathOperator{\Tr}{Tr}
\DeclareMathOperator{\pr}{pr}

\begin{document}

\title[Regularity of matrix coefficients]{Regularity of matrix coefficients of a compact symmetric pair of Lie groups}

\author{Guillaume Dumas}
\address{Université Claude Bernard Lyon 1, ICJ UMR5208, CNRS, Ecole Centrale de Lyon, INSA Lyon, Université Jean Monnet,
69622 Villeurbanne, France.}

\email{gdumas@math.univ-lyon1.fr}
\thanks{}

\subjclass[2010]{Primary 22E46; Secondary 43A85, 43A90}

\date{July 19, 2023}

\dedicatory{}

\begin{abstract}
We consider Gelfand pairs $(G,K)$ where $G$ is a compact Lie group and $K$ a subgroup of fixed points of an involutive automorphism. We study the regularity of $K$-bi-invariant matrix coefficients of unitary representations of $G$. The results rely on the analysis of the spherical functions of the Gelfand pair $(G,K)$. When the symmetric space $G/K$ is of rank $1$ or isomorphic to a Lie group, we find the optimal regularity of $K$-bi-invariant matrix coefficients of unitary representations. Furthermore, in rank $1$ we also find the optimal regularity of $K$-bi-invariant Herz-Schur multipliers of $S_p(L^2(G))$. We also give a lower bound for the optimal regularity in some families of higher rank symmetric spaces. From these results, we make a conjecture in the general case involving the root system of the symmetric space. Finally, we prove that if all $K$-bi-invariant matrix coefficients of unitary representations of $G$ are $\alpha$-Hölder continuous for some $\alpha>0$, then all $K$-finite matrix coefficients of unitary representations are also $\alpha$-Hölder continuous.
\end{abstract}

\maketitle

\section{Introduction}
In this article, we investigate the regularity and local behaviour of $K$-bi-invariant and $K$-finite matrix coefficients of unitary representations of some Lie groups $G$ with compact subgroup $K$. More precisely, we want to find some $\alpha>0$ such that any $K$-finite matrix coefficients of a unitary representation of $G$ is locally $\alpha$-Hölder continuous. Since continuous group morphisms of Lie groups are smooth, if $\pi$ is a finite dimensional unitary representation of $G$, then its coefficients are smooth. In particular, by the Peter-Weyl theorem, every matrix coefficient of an irreducible unitary representation of a compact Lie group $G$ is smooth. In the non-compact setting, if $G$ is semisimple and $K$ is a maximal compact subgroup of $G$, it is known by the work of Harish-Chandra (\cite{harishchandra}) that matrix coefficients associated to $K$-finite vectors of irreducible unitary representations of $G$ (and more generally admissible representations) are $C^\infty$. Unitary representations of $G$ decompose as direct integrals of irreducible representations, but since estimates depend on the representations, it does not provide any estimates for arbitrary representations.

The previous example is part of a larger class of pairs where the subgroup $K$ is such that $(G,K)$ is a Gelfand pair. In that case, there is a 1-1 correspondence between positive-definite spherical functions of $(G,K)$, which arise as characters of the (abelian) convolution algebra of $K$-bi-invariant functions, and irreducible unitary representations with non-zero $K$-invariant vectors. In this setting, any $K$-bi-invariant matrix coefficient of a unitary representation of $G$ will decompose as a direct integral of positive-definite spherical functions (Section \ref{sec:gelfand}). Then, proving estimates of positive-definite spherical functions that are uniform on the whole family will produce estimates on any $K$-bi-invariant matrix coefficient of a unitary representation of $G$.

This idea was used by Lafforgue to show that $SO(2)$-bi-invariant coefficients of $SO(3)$ are $\frac{1}{2}$-Hölder outside of singular points, which was a key ingredient in the proof of his strengthening of property (T) for $SL(3,\R)$ (\cite{lafforgue}). With this result, he showed asymptotic estimates of coefficients of representations on Banach spaces, far more general than unitary representations on Hilbert spaces (see also \cite{de2022analysis} for a survey). This was further improved in \cite{delaatdls} for all higher rank simple Lie groups on a large class of Banach spaces. Various problems in operator algebra were solved using this idea, applied to different pairs, but also in other fields. The same ideas extend to the study of Fourier multipliers (\cite{dlSPR}). In \cite{laffdls}, \cite{Haagerup_2013}) and \cite{haagerupdelaat2}, it is used to show that higher rank simple Lie groups do not have the approximation property. It also implies results on non-coarse embeddability of families of expanders (\cite{deLaatdelaSalle+2018+49+69}). In dynamical systems, strong property $(T)$ was also an important part in the Zimmer's program (\cite{zimmer},\cite{zimmer2},\cite{zimmer3},\cite{brown2022lattice},\cite{fisher2022rigidity}).

In this article, we want to study more systematically this local regularity for symmetric Gelfand pairs, that is when $G/K$ is a symmetric space. In this framework, a lot is known on the spherical functions of the pair (\cite{helgason1979differential},\cite{helgason2000groups}). These functions are parameterized by a subset $\Lambda$ of $\mathfrak{g}_\C^*$, where $\mathfrak{g}$ is the Lie algebra of $G$. We call $\lambda\in \Lambda$ the spectral parameter (see Section \ref{sec:root}). The asymptotics of spherical functions in the group variable have been studied a lot, in particular in the work of Harish-Chandra, but much less is known for local behaviour. To get information on local behaviour, the main tool is to study the asymptotic behaviour as the spectral parameter goes to infinity, while the group variable remains in a compact subset. An important remark is that the estimates we want to show fail at the identity of $G$, and in fact at all other singular points (see Definition \ref{def:regularpoints} for the definition of singular and regular points). Thus, we will only obtain regularity of matrix coefficients on the dense open subset of regular points.

In \cite{clerc}, Clerc showed estimates on spherical functions of compact pairs, but that are uniform only in a cone with compact basis in the Weyl chamber, and do not apply to the differentials of the spherical functions. Cowling and Nevo showed precise estimates, uniform in the spectral parameter for the directional derivatives of the spherical functions of any pair $(G_\C,G)$ where $G_\C$ is a complex semisimple Lie group and $G$ a maximal compact subgroup (\cite{cowling2001uniform}). These estimates are very precise, but they are not enough to obtain regularity results as they are only directionals. Furthermore, we will restrict ourselves to the case where $G$ is compact, thus we will not study these pairs. However, in Section \ref{caracteres}, we will show similar rates of decay for spherical functions of the pair $(G\times G,G)$, which is such that $(G\times G)/G$ is the compact dual of the symmetric space $G_\C/G$. This indicates a possible link between regularities in the compact and non-compact settings.

\subsection*{Main results}
\begin{definition}\label{holdermultivar}Let $(X,d)$ be a metric space and $U$ open subset of $X$, $(E,\Vert.\Vert)$ a normed vector space, $\alpha\in ]0,1]$. A function $f:U\to E$ is $\alpha$-Hölder if for any compact subset $K$ of $U$, there is $C_K>0$ such that $\forall x,y\in K$, $\Vert f(x)-f(y)\Vert \leq C_Kd(x,y)^\alpha$.

If $(X,d)$ is furthermore a Riemannian manifold and $r\in \N$, we say that the map $f$ belongs to $C^{(r,\alpha)}(U,E)$ if $f\in C^r(U,E)$ and the $r$-th differential $D^rf$ is $\alpha$-Hölder. We extend to $\alpha=0$ by $C^{(r,0)}(U,E)=C^r(U,E)$.

For $K$ a compact subset of $U$ and $f\in C^{(r,\alpha)}(U,E)$, define $$\Vert f\Vert_{C^{(r,\alpha)}(K,E)}=\max \left\{\underset{k\leq r}{\max}\, \underset{x\in K}{\sup}\Vert D^kf(x)\Vert,\underset{x,y\in K,x\neq y}{\sup}\frac{\Vert D^rf(x)-D^rf(y)\Vert}{d(x,y)^\alpha}\right\}.$$
The family of semi-norms $\Vert.\Vert_{C^{(r,\alpha)}(K,E)}$ for $K$ a compact subset of $U$ makes the space $C^{(r,\alpha)}(U,E)$ into a Fréchet space.
\end{definition}
\begin{remark}
    If $U$ is locally compact, a function $f:U\to E$ is $\alpha$-Hölder in the sense of the previous definition if and only if $f$ is locally $\alpha$-Hölder - that is to say that for any $x\in U$, there exists a neighborhood $U_x$ of $x$ and a constant $C_x>0$ such that for any $y,z\in U_x$, $\Vert f(y)-f(z)\Vert\leq C_xd(y,z)^{\alpha}$.
\end{remark}
We will denote $C^{(r,\alpha)}(U,\C)$ by $C^{(r,\alpha)}(U)$.

Let $G$ be a Lie group and $K$ a compact subgroup. We want to find $(r,\alpha)$ such that for any $K$-finite unitary matrix coefficient $\varphi$ of $G$, $\varphi\in C^{(r,\alpha)}(G_1)$, where $G_1$ will be the dense open subset of regular points (see Definition \ref{def:regularpoints}). If we assume that $(G,K)$ is a Gelfand pair, the main tool for this will be to study the boundedness of the family of spherical functions in Hölder spaces. If furthermore $(G,K)$ is a symmetric pair, then $G/K$ is a symmetric space and there are results on the spherical functions in the literature.

Our first results involve symmetric pairs of rank $1$ (see Setion \ref{sec:root} for the notion of rank).
\begin{mainthm}
\label{mainrank1}Let $(G,K)$ be a compact symmetric pair of rank $1$ and $G_1$ the dense open subset of regular points. Let $$\alpha=\frac{\dim G/K-1}{2}.$$Then any $K$-finite unitary matrix coefficient $\varphi$ of $G$ is in $C^{(\lfloor\alpha\rfloor,\alpha-\lfloor\alpha\rfloor)}(G_1)$. Furthermore, this regularity is optimal in the sense that for any $(r,\delta)>(\lfloor\alpha\rfloor,\alpha-\lfloor\alpha\rfloor)$, there exists a $K$-finite (and even $K$-bi-invariant) unitary matrix coefficient of $G$ not in $C^{(r,\alpha)}(G_1)$.
\end{mainthm}
We also show in Corollary \ref{spmult} a regularity result for Herz-Schur multipliers of $S_p(L^2(G))$ for any $p$.

We then turn to higher rank symmetric pairs. Among higher rank pairs, a simple class is given by the pairs $(G\times G,G)$ with $G$ a compact semisimple Lie group. The symmetric space associated to these pairs is isomorphic to the Lie group $G$. Given $\Phi$ a root system for $G$, $\Phi^+$ a choice of positive roots and $\Delta=\{\alpha_1,\cdots,\alpha_\ell\}$ a basis, we can write $\alpha=\sum_{i=1}^\ell n_i(\alpha)\alpha_i$ for $\alpha\in \Phi^+$ (see Section \ref{sec:root} for more details about root systems). Then our second result gives the regularity for these pairs.
\begin{mainthm}\label{maingroupeLie}Let $\gamma=\underset{1\leq i\leq \ell}{\min} \vert \{\alpha \in \Phi^+ \vert n_i(\alpha)\geq 1\}\vert$ and $G_1$ the subset of regular points of $G\times G$. Then any $G$-finite matrix coefficient $\varphi$ of a unitary representation of $G\times G$ is in $C^{(\gamma,0)}(G_1)$ and this regularity is optimal.
\end{mainthm}

We also give regularity results for some other families of higher rank pairs in Theorems \ref{roptGrassC2} and \ref{highrankbis}. However, we cannot prove optimality for these pairs. Given all these results, we make a conjecture on the optimal regularity in the general case. Given $(G,K)$ a compact symmetric pair, there is a decomposition $\mathfrak{g}=\mathfrak{k}\oplus \mathfrak{m}$ of the Lie algebra. Let $\mathfrak{a}$ be a maximal abelian subspace of $\mathfrak{m}$ and $\Sigma$ the root system pair. Let $\Lambda=\{\mu\in i\mathfrak{a}^* \vert \forall \alpha\in \Sigma^+,\frac{\langle\mu,\alpha\rangle}{\langle\alpha,\alpha\rangle}\in \N\}$. For $\alpha\in \Sigma$, let $m_\alpha$ be the multiplicity of the root. 
\begin{Conj}Let $G_1$ be the subset of regular points of $G$. Let $$\alpha=\underset{\mu\in \Lambda\setminus \{0\}}{\inf} \underset{\alpha\in \Sigma^+,<\alpha,\mu>\neq 0}{\sum} \frac{m_\alpha}{2}.$$Then any $K$-finite matrix coefficient $\varphi$ of a unitary representation of $G$ is in $C^{(\lfloor\alpha\rfloor,\alpha-\lfloor\alpha\rfloor)}(G_1)$ and this regularity is optimal.
\end{Conj}

\subsection*{Organisation of the paper}
In Section \ref{sectionpreli}, we recall some results that will be used throughout the paper. Until Section \ref{kfinitesection}, we restrict ourselves to $K$-bi-invariant matrix coefficients of unitary representations.

In Section \ref{sectionrank1}, we give a complete answer in the case of a symmetric pair of rank $1$. There are well-known descriptions of the spherical functions in this case involving Jacobi polynomials (Section \ref{sectionspheriques}). The main result is Theorem \ref{mainthm}, where we use analysis of these polynomials to show that the family of spherical functions is bounded in some Hölder space. Similar estimates were used before to study different families of Jacobi polynomials (\cite{haagschli},\cite{Haagerup_2013}). From this, we deduce Theorem \ref{mainrank1} in the case of $K$-bi-invariant matrix coefficients of unitary representations. This generalizes results used in \cite{lafforgue},\cite{Haagerup_2013},\cite{delaatdls} and \cite{deLaatdelaSalle+2018+49+69} for some specific pairs. In Section \ref{schattensection}, using the aforementioned results, we prove in Corollary \ref{spmult} regularity results for any $K$-bi-invariant Herz-Schur multipliers of $S_p(L^2(G))$, which improve on \cite{dlSPR}. We then show in Section \ref{optimalitysection} that our results are optimal (Theorem \ref{optimality}).

In Section \ref{higherranksection}, we try to extend the results to higher rank. For the case $(G\times G,G)$, we can give the optimal regularity (Theorem \ref{regcar} and Theorem \ref{opticar}), which proves Theorem \ref{maingroupeLie} for $K$-bi-invariant coefficients. For this, we rely on the description of spherical functions of these pairs with characters of $G$ and use the Weyl character formula. The techniques involved in the optimality are very similar to the rank $1$ case, but require the study of the root system associated to $G$. We then show, for some of the remaining pairs, results of regularity (that may not be optimal) in Section \ref{sectionpartialresult}. We study these specific pairs because once again their spherical functions can be described with Jacobi polynomials. These different cases allow us to formulate the above conjecture.

Finally, in Section \ref{kfinitesection}, we show that knowing the regularity of any \emph{$K$-bi-invariant} unitary matrix coefficient is sufficient to get a regularity result on any \emph{$K$-finite} unitary matrix coefficient. The crucial result is Lemma \ref{dlmdls}, which was already used in \cite{AIF_2016__66_5_1859_0} but only to get Hölder continuity results. Given the previous results, it completes the proof of Theorem \ref{mainrank1} and \ref{maingroupeLie} for $K$-finite matrix coefficients.

\subsection*{Acknowledgement}
I would like to thank my Ph.D supervisor Mikael de la Salle for his involvement and support. I am also grateful to the referees for their constructive comments.

\section{Preliminaries}\label{sectionpreli}
\subsection{Hölder spaces}
We recall here standard results on Hölder spaces. The following lemma is a generalisation of \cite[Proposition 4.1]{dlSPR}.
\begin{lem}
\label{precomposition}Let $\alpha>0$ and $(X,d),(Y,d')$ be two Riemannian manifolds and $U,V$ open subsets of $X,Y$ respectively. Let $\varphi:U\to V$ be a function of class $C^{\infty}$. Then $\varphi_*:f\mapsto f\circ \varphi$ maps $C^{(r,\alpha)}(V)$ to $C^{(r,\alpha)}(U)$. Furthermore, if $(f_n)$ is bounded in $C^{(r,\alpha)}(V)$, then $(f_n\circ \varphi)$ is bounded in $C^{(r,\alpha)}(U)$.
\end{lem}

\begin{proof}
    Let $f\in C^{(r,\alpha)}(V)$. Then $f\circ \varphi\in C^r(U)$. By \cite[Thm. 6.17]{lee2019introduction}, for any $x\in X$, there exists $\varepsilon_0$ such that for any $\varepsilon\leq \varepsilon_0$, the geodesic ball $B(x,\varepsilon)$ is geodesically convex. Then there is $U_x$ open neighborhood of $x$, geodesically convex and such that $L=\overline{U_x}$ is compact and a subset of $U$.\\
The map $D^r(f\circ \varphi)(x)$ is a sum of terms of the form $$D^if(\varphi(x))\circ (D ^{j_1}\varphi(x),\cdots,D^{j_i}\varphi(x))$$ with $i\leq r$ and $\sum {j_i}=r$. Take $y,z\in U_x$, then for the term $i=r$, we have \begin{multline*}
    \Vert D^rf(\varphi(y))\circ (D\varphi(y),\cdots,D\varphi(y))-D^rf(\varphi(z))\circ (D\varphi(z),\cdots,D\varphi(z))\Vert \\ \begin{aligned}
        & \leq \Vert [D^rf(\varphi(y))-D^rf(\varphi(z))]\circ \left(D\varphi(y),\cdots,D\varphi(y)\right)\Vert\\
        & \phantom{\leq}+ \sum_{k=1}^r\Vert D^rf(\varphi(y)) \circ (D\varphi(y),\cdots,D\varphi(y)-D\varphi(z),\cdots,D\varphi(z))\Vert\\
        & \leq \Vert f\Vert_{C^{(r,\alpha)}(\varphi(L))} \Vert \varphi\Vert_{C^1(L)}^r d(\varphi(y),\varphi(z))^\alpha\\
        & \phantom{\leq} + \sum_{k=1}^r \Vert f\Vert_{C^{(r,\alpha)}(\varphi(L))}\Vert \Vert\varphi\Vert_{C^1(L)}^{r-1}\Vert D\varphi(y)-D\varphi(z)\Vert\\
        & \leq \Vert f\Vert_{C^{(r,\alpha)}(\varphi(L))} \Vert \varphi\Vert_{C^1(L)}^{r+\alpha}d(y,z)^\alpha\\
        & \phantom{\leq}+ \sum_{k=1}^r \Vert f\Vert_{C^{(r,\alpha)}(\varphi(L))}\Vert \Vert\varphi\Vert_{C^1(L)}^{r-1}\Vert\varphi\Vert_{C^2(L)}d(y,z)
    \end{aligned}
\end{multline*}
We can use the mean value theorem to get the last inequality because $U_x$ is geodesically convex. There is a constant $C(U_x,f,\varphi)>0$ such that any term with $i<r$ is bounded in a similar way by $C(U_x,f,\varphi)d(y,z)$. The diameter of $U_x$ is bounded so there is $C>0$ such that $d(y,z)\leq C d(y,z)^\alpha$. Thus, $f\circ \varphi \in C^{(r,\alpha)}(U)$.

Finally, $\varphi_*$ is a linear map between Fréchet spaces. If $(f_n,f_n\circ \varphi)$ converges to $(f,g)$ in $C^{(r,\alpha)}(V)\times C^{(r,\alpha)}(U)$, then $f_n$ converge uniformly to $f$ on compact subset and $f_n\circ \varphi$ converge uniformly to $g$ on compact subset. Thus, $g=f\circ \varphi$ and the graph of $\varphi_*$ is closed. By the closed graph theorem, $\varphi_*$ is continuous. Thus, the image of a bounded subset is bounded.
\end{proof}

Using Leibniz rule, we can prove the following lemma.
\begin{lem}
 \label{multsmooth}Let $U$ be an open subset of a normed vector space $E$ and $g:U\to \C$ a smooth function. If $(f_i)_{i\in I}$ is bounded in $C^{(r,\alpha)}(U)$, then $(gf_i)$ is bounded in $C^{(r,\alpha)}(U)$. 
\end{lem}
\begin{proof}
    It is clear that if $f\in C^{(r,\alpha)}(U)$, $gf$ is $r$ times differentiable and for $x\in U$, $H=(H_1,\cdots,H_r)\in E^n$, $$D^r(gf)(x)(H)=\sum_{j\in J} D^{k_j}g(x)(H_{p_1(j)},\cdots,H_{p_{k_j}(j)}) D^{r-k_j}f(x)(H_{p_{k_j+1}(j)},\cdots,H_{p_r(i)})$$where $k_j\leq r$ for each $j\in J$ and $J$ finite. Thus, we see that $gf\in C^{(r,\alpha)}(U)$.\\
    Furthermore, $m:f\mapsto gf$ is a well-defined linear map between Fréchet spaces. If $(f_n,gf_n)$ converges uniformly to $(f,h)$ in $C^{(r,\alpha)}(U)\times C^{(r,\alpha)}(U)$, then $f_n$ converges uniformly to $f$ on compact subsets of $U$ and $gf_n$ converges uniformly on compact subsets to $h$, and to $gf$. Thus $gf=h$, so the graph of $m$ is closed. By the closed graph theorem, $m$ is continuous so the image of a bounded subset if bounded.
\end{proof}

Given a family of functions which are eigenvalues of a map $T$ into $B(\mathcal{H})$, the regularity of $T$ can be linked to the boundedness of the eigenvalues in Hölder spaces.
\begin{lem}
\label{existencederiv}Let $U$ be an open subset of $\R^d$ and $\mathcal{H}$ a Hilbert space. Consider a map $T:U\mapsto B(\mathcal{H})$ such that there is an orthonormal basis $(e_n)$ of $\mathcal{H}$ and a family of maps $f_n:U\mapsto \R$, such that for any $X\in U$, $T(X)$ is diagonal in the basis $(e_n)$ with eigenvalues $(f_n(X))$. If $(f_n)$ is bounded in $C^{(r,\alpha)}(U)$ for some $\alpha>0$, then $T$ lies in $C^{(r,\alpha)}(U,B(\mathcal{H}))$.
\end{lem}
\begin{proof}If $T$ is $C^k$, then for any $X\in U$, $H_1,\cdots,H_k\in \R^d$, we must have $f_n\in C^k$ for any $n\in \N$ and $$D^kT(X)(H_1,\cdots,H_k)e_n=D^kf_n(X)(H_1,\cdots,H_k)e_n.$$
We prove the result by induction on $r$. If $S\in B(\mathcal{H})$, let $\Vert S\Vert_\infty=\underset{\Vert H\Vert=1}{\sup}\Vert S(H)\Vert$ be its operator norm. If $r=0$, we have for $X,Y\in L$ compact subset of $U$, $$
    \Vert T(X)-T(Y)\Vert_\infty  =  \sup_{n} \vert f_n(X)-f_n(Y)\vert \leq  C_L \Vert X-Y\Vert^\alpha$$so $T$ is $\alpha$-Hölder.

Assume the result is true for $r-1$. Since $C^{(r,\alpha)}(U,B(\mathcal{H}))\subset C^{(r-1,1)}(U,B(\mathcal{H}))$, by the induction hypothesis, we have $T\in C^{r-1}(U,B(\mathcal{H}))$. Let $A$ be the $r$-linear map such that $$A(H_1,\cdots,H_r)e_n=D^rf_n(X)(H_1,\cdots,H_r)e_n.$$ Note that since $(D^rf_n)$ is bounded in $n$ on compacts, $A(H_1,\cdots,H_r)\in B(\mathcal{H})$.

Denote $H=(H_1,\cdots,H_{r-1})$ and $H^+=(H_1,\cdots,H_{r})$. We must show that $$\underset{H_r\to 0}{\lim} \frac{\Vert D^{r-1}T(X+H_r)(H)-D^{r-1}T(X)(H)-A(H^+)\Vert_\infty}{\Vert H_r\Vert_1}=0$$uniformly for $H_1,\cdots,H_{r-1}$ in bounded sets. We have \begin{multline*}
        \Vert D^{r-1}T(X+H_r)(H)-D^{r-1}T(X)(H)-A(H^+)\Vert_\infty\\=\underset{n}{\sup} \vert D^{r-1}f_n(X+H_r)(H)-D^{r-1}f_n(X)(H)-D^rf_n(X)(H^+)\vert.
\end{multline*}

Assume that $H_r$ is small enough, so that $B(X,\Vert H\Vert_\infty)\subset U$. Let $$\fonction{g_n}{[0,1]}{\R}{t}{D^{r-1}f_n(X+tH_r)(H_1,\cdots,H_{r-1})}.$$Then $$g_n(1)-g_n(0)=g_n'(t)=D^rf_n(X+tH_r)(H_1,\cdots,H_r)$$for some $t\in ]0,1[$.

Then \begin{multline*}
    \vert D^{r-1}f_n(X+H_r)(H)-D^{r-1}f_n(X)(H)-D^rf_n(X)(H^+)\vert \\
    \begin{aligned}
       & = \vert D^rf_n(X+tH_r)(H^+)-D^rf_n(X)(H^+)\\
       & \leq \Vert D^rf_n(X+tH_r)-D^rf_n(X)\Vert \prod_{i=1}^r \Vert H_i\Vert\\
       & \leq \Vert X+tH_r-X\Vert^\alpha \prod_{i=1}^r \Vert H_i\Vert\\
       & \leq \Vert H_r\Vert^\alpha \prod_{i=1}^r \Vert H_i\Vert
    \end{aligned}
\end{multline*}
so since $\alpha>0$, we get what we want.

So we have that $D^rT(X)$ exists for any $X\in U$, and thus for $X,Y\in L$ compact subset of $U$, $$\Vert D^r T(X)-D^rT(Y)\Vert = \underset{n}{\sup} \Vert D^rf_n(X)-D^rf_n(Y)\Vert \leq C_L\Vert X-Y\Vert^\alpha$$because $D^rf_n$ are uniformly $\alpha$-Hölder on $L$.
\end{proof}

\subsection{Gelfand pairs}\label{sec:gelfand}
\begin{definition}Let $G$ be a locally compact topological group with a left Haar measure $dg$ and $K$ a compact subgroup with normalized Haar measure $dk$. The pair $(G,K)$ is a Gelfand pair if the algebra of continuous $K$-bi-invariant functions on $G$ with compact support is commutative for the convolution.

A spherical function of $(G,K)$ is a continuous $K$-bi-invariant non-zero function on $G$ such that for all $x,y\in G$, $$\int_K \varphi(xky)\,dk=\varphi(x)\varphi(y).$$
\end{definition}

A standard result (see \cite[Coro. 6.3.3]{Dijk+2009}) gives a link between spherical functions of $(G,K)$ and unitary representations of $G$.
\begin{prop}
If $(G,K)$ is a Gelfand pair, then for any irreducible unitary representation $\pi$ of $G$ on a Hilbert space $\mathcal{H}$, the subspace $\mathcal{H}^K$ of $K$-invariant vectors is of dimension at most $1$.

The positive-definite spherical functions of $G$ are exactly the matrix coefficients $g\mapsto \langle\pi(g)v,v\rangle$ with $\pi$ an irreducible unitary representation of $G$ and $v$ a $K$-invariant unit vector.

If $G$ is compact, any spherical function is positive-definite.
\end{prop}
\begin{remark}If we assume that $G,K$ are Lie groups, positive-definite spherical functions have a geometric interpretation. Let $D(G/K)$ be the algebra of differential operators on $G/K$ invariant by the action of $G$ by translation on $G/K$. Then $\varphi:G/K\to \C$ is a positive-definite spherical function if and only if $\varphi(K)=1$, $\varphi$ is invariant by the action of $K$ and $\varphi$ is an eigenvalue of all operators of $D(G/K)$.

When $G/K$ is a compact symmetric space of rank $1$ (a sphere or a projective space, see Section \ref{sectionrank1}), then $D(G/K)$ is generated by the Laplacian. Hence in that case, positive-definite spherical functions are normalized $K$-invariant eigenvalues of the Laplacian.
\end{remark}
More details on Gelfand pairs can be found in \cite[Ch. 5,6,7]{Dijk+2009}.

Given a Gelfand pair $(G,K)$, it is natural to study spherical functions in order to get results on $K$-bi-invariant matrix coefficients of unitary representations. Indeed, any matrix coefficient of a unitary representation decomposes into an integral of spherical functions - an infinite sum if $G$ is compact.
\begin{lem}
Let $(G,K)$ be a Gelfand pair with $G$ second countable. Let $\varphi$ be a $K$-bi-invariant matrix coefficient of a unitary representation $\pi$ on an Hilbert space $\mathcal{H}$. Then, there exists a standard Borel space $X$ and a $\sigma$-finite measure $\mu$ on $X$ such that $$\varphi=\int_X c_x\varphi_x d\mu(x)$$where $\varphi_x$ is a positive-definite spherical function of $(G,K)$ for any $x\in X$ and $c\in L^1(X,\mu)$.
\end{lem}
\begin{proof}\label{deccoef}If $\varphi(g)=\langle \pi(g)u,v\rangle$, we can replace $\mathcal{H}$ by $\overline{\vspan(\pi(G)u,\pi(G)v)}$ which is a $G$-invariant separable subspace, since $G$ is second-countable hence separable. Thus, we can assume that $\mathcal{H}$ is separable.

Then, by \cite[Section 8.4]{kirillov1976elements}, there exists $(X,\mu)$ and an isometry $U:\mathcal{H}_X\to \mathcal{H}$ where $\mathcal{H}_X$ is the direct integral of the collection of Hilbert spaces $(\mathcal{H}_x)_{x\in X}$, such that $\forall g\in G$, $\pi(g)=U\circ \Tilde{\pi}(g)\circ U^{-1}$, where $(\Tilde{\pi}(g)\xi)_x=\pi_x(g)\xi_x$ and $(\pi_x,\mathcal{H}_x)$ is an irreducible unitary representation of $G$.

Let $P$ denote the projection on the space of $K$-invariant vectors in $\mathcal{H}$. Since $\varphi$ is $K$-bi-invariant, we have $$\varphi(g)=\int_K\int_K \varphi(kgk')\,dk\,dk'=\langle\pi(g)Pu,Pv\rangle.$$Thus, we can assume that $u,v$ are $K$-invariant, so $\xi=U^{-1}u,\eta=U^{-1}v$ are $K$-invariant.

So for $\mu$ almost every $x$, $\xi_x,\eta_x$ are $K$-invariant. Now, if $\pi_x$ is such that $0$ is the only $K$-invariant vector, $\langle\pi_x(g)\xi_x,\eta_x\rangle_{\mathcal{H}_x}=0$ for every $g\in G$. On the other hand, assume $\pi_x$ has non-zero $K$-invariant vectors. Then we know that the space of $K$-invariant vector is one-dimensional and that there is $c_x\in \C$ such that $\langle\pi_x(g)\xi_x,\eta_x\rangle_{\mathcal{H}_x}=c_x\varphi_{x}$, where $\varphi_{x}$ is the spherical function associated to $\pi_x$ (and so is positive-definite).

Setting $c_x=0$ if $\pi_x$ has no non-zero $K$-invariant vectors and $\varphi_x$ the constant spherical function, we have $$\varphi(g)=\langle\Tilde{\pi}(g)\xi,\eta\rangle_{\mathcal{H}_X}=\int_X \langle\pi_x(g)\xi_x,\eta_x\rangle_{\mathcal{H}_x}\,d\mu(x)=\int_{X}c_x\varphi_{x}(g)\,d\mu(x).$$Since $\vert c_x\vert \leq \Vert \xi_x\Vert_{\mathcal{H}_x} \Vert \eta_x\Vert_{\mathcal{H}_x}$, the function $x\mapsto c_x$ is in $L^1(X,\mu)$.
\end{proof}

\begin{lem}
\label{roptid}Let $(G,K)$ be a Gelfand pair with $G$ a Lie group endowed with a Riemannian metric $d$ and $U$ any open subset of $G$. Let $(\varphi_\lambda)_{\lambda\in \Lambda}$ be the the family of positive-definite spherical functions of $(G,K)$. Then $(\varphi_\lambda)_{\lambda\in \Lambda}$ is bounded in $C^{(r,\delta)}(U)$ if and only if any $K$-bi-invariant matrix coefficient of a unitary representation of $G$ is in $C^{(r,\delta)}(U)$.
\end{lem}
\begin{proof}
Assume that $(\varphi_\lambda)_{\lambda\in \Lambda}$ is bounded in $C^{(r,\delta)}(U)$. Let $\varphi$ be a $K$-bi-invariant matrix coefficient of a unitary representation of $G$. Since $G$ is a Lie group, $G$ is second countable, thus by the above lemma, there exists a $\sigma$-finite measured space $(X,\mu)$ such that $$\varphi=\int_X c_x\varphi_{\lambda_x}d\mu(x)$$with $c\in L^1(X,\mu)$ and $\lambda_x\in \Lambda$ for any $x\in X$.

Let $g_0\in U$ and $L$ a compact neighborhood of $g_0$. Then there exists $C_{L}>0$ such that $\Vert D^k \varphi_{\lambda_x}(g)\Vert\leq C_L$ for any $x\in X,g\in L,k\leq r$. Thus, $(x,g)\mapsto D^k (c_x\varphi_{\lambda_x})(g)$ is bounded by the integrable function $x\mapsto C_Lc_x$. Hence $\varphi$ is $r$ times differentiable in a neighborhood of $g_0$, and $$D^k\varphi(g)=\int_X c_x D^k\varphi_{\lambda_x}(g) d\mu(x)$$for any $k\leq r$. Since this hold for any $g_0\in U$, $\varphi\in C^r(U)$. Finally, let $L$ be any compact subset of $U$. There exists $D_L>0$ such that for any $g,h\in L$ and $x\in X$, $$\Vert D^r\varphi_{\lambda_x}(g)-D^r\varphi_{\lambda_x}(h)\Vert \leq D_L d(g,h)^\delta.$$Thus, $$\Vert D^r\varphi(g)-D^r\varphi(h)\Vert \leq \int_x \vert c_x\vert \Vert D^r\varphi_{\lambda_x}(g)-D^r\varphi_{\lambda_x}(h)\Vert d\mu(x) \leq D_L\Vert c\Vert_1 d(g,h)^\delta.$$Hence, we showed that $\varphi\in C^{(r,\delta)}(U)$.

For the other direction, assume that any $K$-bi-invariant matrix coefficient of a unitary representation of $G$ is in $C^{(r,\delta)}(U)$. Let $E$ be the space of $K$-bi-invariant matrix coefficient of unitary representations, endowed with the norm $$\Vert \varphi\Vert = \inf \left\lbrace \Vert \xi\Vert \Vert \eta\Vert \, \vert \, \exists \pi \textrm{ such that } \forall g\in G,\varphi(g)=\langle\pi(g)\xi,\eta\rangle\right\rbrace.$$Then $E$ is a Banach space, and $\Vert \varphi\Vert \geq \Vert \varphi\Vert_\infty$. Consider $f:E\mapsto C^{(r,\delta)}(U)$ the linear map sending $\varphi$ to its restriction to $U$. By the assumption on regularity, $f$ is well-defined. Let $\mathcal{G}_f$ be the graph of $f$ in $E\times C^{(r,\delta)}(U)$. We claim that it is closed. Indeed, if $(\varphi_n,\varphi_n)\to (\varphi,\psi)$, then in particular, $\varphi_n$ converges to $\varphi$ uniformly on compact subsets of $G$, since $\Vert \varphi_n-\varphi\Vert_\infty \leq \Vert \varphi_n-\varphi\Vert$. On the other hand, by definition of the seminorms on $C^{(r,\delta)}(U)$, $\varphi_n$ converges uniformly on compact subset of $U$ to $\psi$. Thus $\psi=\varphi\vert_U$. Now, since $E$ is a Banach space and $C^{(r,\delta)}(U)$ a Fréchet space, by the closed graph theorem, $f$ is continuous. Finally, since the family of positive-definite spherical functions is in the unit ball of $E$, its image in $C^{(r,\delta)}(U)$ is bounded.
\end{proof}
This result shows that studying boundedness of (positive-definite) spherical functions is enough to obtain regularity for all $K$-bi-invariant matrix coefficients of unitary representations, and even that the optimal regularity of such coefficients is exactly the optimal uniform regularity of spherical functions.

\subsection{Symmetric spaces}
Let $G$ be a connected Lie group $G$ and $\sigma$ an involutive automorphism of $G$. Let $G^\sigma$ denote the subgroup of fixed points of $\sigma$ and $(G^\sigma)_0$ its identity component. For a subgroup $K$ of $G$ such that $(G^\sigma)_0\subset K\subset G^\sigma$, the quotient space $X=G/K$ can be given a structure of symmetric space. Using a characterization of Gelfand pairs (\cite[Prop. 6.1.3]{Dijk+2009}), we see that $(G,K)$ is a Gelfand pair if and only if $K$ is compact. We call such pairs symmetric Gelfand pairs.

It turns out that all symmetric spaces arise in this way (see \cite[Ch. II, Thm. 3.1]{loos1969symmetric}). Given $M$ a connected symmetric space and $o\in M$, there is a canonical connected Lie group $G(M)$ associated to it, called the group of displacements. The group $G(M)$ is a subgroup of $\mathrm{Aut}(M)$, the group of automorphisms of symmetric spaces of $M$, thus acts of $M$. If $K(M)$ is the isotropy subgroup of $o$, we have $M\simeq G(M)/K(M)$. Furthermore, $K(M)$ is compact if and only if the symmetric space is Riemannian (\cite[Ch. IV, Prop. 1.7]{loos1969symmetric}).

Since there is a classification of Riemannian symmetric spaces (\cite[Ch. VII]{loos1969symmetric2}), it seems natural to study symmetric spaces. We say that a Riemannian symmetric space $M$ is\begin{itemize}
    \item euclidean if its sectional curvature is identically zero,
    \item of compact type if its sectional curvature is non-positive and not identically zero,
    \item of non-compact type if its sectional curvature is non-negative and not identically zero.
\end{itemize}If $M$ is a simply connected symmetric space, then there are $M_0$ euclidean, $M_+$ of non-compact type and $M_-$ of compact type such that $M=M_0\times M_+\times M_-$ (\cite[Ch. IV, Coro. 1]{loos1969symmetric}). In this paper, we will study symmetric Gelfand pairs, and more precisely those where $M$ is of compact type. 
However, a question arises: if two pairs represent the same symmetric space, are their spherical functions the same ?

We say that a symmetric space $M$ is semisimple if $G(M)$ is a semisimple Lie group.  
\begin{lem}
\label{indepgelfandpair}Let $(G,K)$ be a symmetric pair and $M=G/K$ the associated symmetric space. If $M$ is semisimple, then there is a bijection between spherical functions of $(G,K)$ and spherical functions of $(G(M),K(M))$, such that the image of $\varphi$ induces the same function as $\varphi$ on $M$.
\end{lem}
\begin{proof}
Let $\tau:G\to \mathrm{\Aut(M)}$ be the morphism defined by $\tau(g):xK\mapsto gxK$. Then $\ker \tau= \bigcap_{g\in G} gKg^{-1}$. Let $G_\sigma=\{x\sigma^{-1}(x)\vert x\in G\}$. By \cite[Ch. II, Thm. 1.3]{loos1969symmetric}, $G(M)=\langle\tau(G_\sigma)\rangle$ is a subgroup of $\tau(G)\simeq G/\ker \tau$.

By \cite[Ch. IV, Prop. 1.4]{loos1969symmetric}, since $M$ is a semisimple symmetric space, we have that $G(M)=(\mathrm{Aut}(M))_0$. So we have $G(M)=(\mathrm{Aut}(M))_0<\tau(G) < \mathrm{Aut}(M)$ and $\tau(G)$ is connected because $G$ is, so $G(M)\simeq G/\ker \tau$, and $K(M)\simeq K/\ker \tau$.

Let $\pi$ be an irreducible unitary representation of $G(M)$ with a $K(M)$-invariant vector $\xi$. By composition with the isomorphism and projection, it induces an irreducible unitary representation of $G$, with $\xi$ which is a $K$-invariant vector.

Conversely, let $\pi$ be an irreducible unitary representation of $G$ on $V$ with a $K$-invariant vector $\xi$. By irreducibility, $\overline{\vspan(\pi(G)\xi)}=V$. Let $x\in \ker \tau$. If $g\in G$, there is $k\in K$ such that $x=gkg^{-1}$. Hence, $$\pi(x)\pi(g)\xi=\pi(g)\pi(k)\xi=\pi(g)\xi.$$So for any $g\in G$, $\pi(g)\xi$ is $\pi(x)$-invariant. By density of the vector space generated by these vectors, $\pi(x)=\Id_V$. So $\ker \tau\subset \ker \pi$, thus $\pi$ induces an irreducible representation of the quotient $G(M)$, with a $K(M)$-invariant vector $\xi$.
\end{proof}

This results says that in the case of a semisimple space, the spherical functions depend essentially only on the symmetric space. But the assumption that $M$ is semisimple is crucial, as the following example illustrates. Let $M$ be the euclidean space $\R^n$. The canonical pair associated is $(\R^n,\{0\})$. But $\R^n\neq (\mathrm{Aut}(M))_0$. For example, $(\R^n\rtimes SO(n),SO(n))$ is a symmetric Gelfand pair and $\R^n$ is the associated symmetric space. In that case, spherical functions of the two pairs are not related.

If $M$ is of compact type or non-compact type, then $M$ is semisimple. In fact (\cite[Ch. IV, Thm. 3.5]{loos1969symmetric}), $M$ is of compact type if and only if $M$ is compact and semisimple, if and only if its universal cover is compact, if and only if $G(M)$ is compact and semisimple.

Thus, as we restrict ourselves to symmetric pairs associated to symmetric spaces of compact type, by the previous discussion, we can study only one pair for each symmetric space (for example, the canonical pair).

\subsection{Spherical functions of compact symmetric pairs}\label{sec:root}
In this section, let $M$ be a compact connected simply connected symmetric space and $(G,K)$ the associated canonical compact symmetric Gelfand pair. Since $G$ is compact, any spherical functions is positive-definite and thus corresponds to an irreducible unitary representation of $G$ with a non-zero $K$-invariant vector.

We know that the finite dimensional irreducible representations of $G$ are classified by the highest weights $\mu$. Let $\mathfrak{g}$ be the Lie algebra of $G$ and $\mathfrak{k}$ the Lie algebra of $K$. Then we have $\mathfrak{g}=\mathfrak{k} \oplus \mathfrak{m}$ where $\mathfrak{k}$ (resp. $\mathfrak{m}$) is the eigenspace of $+1$ (resp. $-1$) of $\sigma$. The space $\mathfrak{m}$ is also the Lie triple system of $M$ (see \cite{loos1969symmetric}, Ch. II, Prop 2.3). Let $\mathfrak{a}$ be a maximal abelian subspace of $\mathfrak{m}$, and $\Sigma_\mathfrak{a}^+$ a choice of positive root system of $\mathfrak{a}_\C$ in $\mathfrak{g}_\C$. Let $\mathfrak{k}^{\mathfrak{a}}=\{x\in \mathfrak{k} \vert [x,\mathfrak{a}]=0\}$ and $\mathfrak{t}$ a Cartan subalgebra of $\mathfrak{k}^\mathfrak{a}$. Then $\mathfrak{h}=\mathfrak{t}+\mathfrak{a}$ is a Cartan subalgebra of $\mathfrak{g}$. Consider $\Sigma^+$ a choice of positive root system of $\mathfrak{h}_\C$, such that $\Sigma_\mathfrak{a}^+=\{\phi\vert_\mathfrak{a} \vert \phi\in \Sigma^+, \phi\vert_\mathfrak{a}\neq 0\}$. Also, recall that by definition, $\rank G=\dim \mathfrak{h}$ and $\rank M=\dim \mathfrak{a}$.

We say that $\mu\in\mathfrak{h}_\C^*$ is a dominant element if for any $\alpha\in \Sigma^+$, $\langle \mu,\alpha\rangle\geq 0$. We say that $\mu$ is analytically integral if for any $H\in \mathfrak{h}$ such that $\exp(H)=1$, $\mu(H)\in 2\pi i\Z$. By the theorem of the highest weight (\cite[Thm 5.110]{knapp2002lie}), there is a one-one correspondence between finite-dimensional irreducible representations of $G$ and dominant analytically integral elements. Let $\mu$ be a dominant analytically integral element in $\mathfrak{h}_\C^*$ and $(\pi_\mu,V_\mu)$ a representative of the associated class of finite dimensional irreducible representations of $G$. Let $$\rho=\frac{1}{2}\sum_{\alpha\in \Sigma^+}\alpha.$$ We know by the Weyl formula (\cite[Thm. 10.18]{hall2003lie}) that $$d_\mu=\dim V_\mu= \frac{\prod_{\alpha\in \Sigma^+}\langle\alpha,\mu+\rho\rangle}{\prod_{\alpha\in \Sigma^+}\langle\alpha,\rho\rangle}.$$

The first question we want to ask is, given $\mu$ a dominant analytically integral element, what are the conditions for $\pi_\mu$ to have a non-zero $K$-invariant vector. In this context, the answer is given by the Cartan-Helgason theorem (\cite[Ch. V, Thm. 4.1]{helgason2000groups}). Let $\hat{G}_K$ denote the set of classes of irreducible finite dimensional representations with a non-zero $K$-invariant vector. 
\begin{theorem}[Cartan-Helgason]
\label{cartanhelgason}Let $\Lambda=\{\mu\in i\mathfrak{a}^* \vert \forall \alpha\in \Sigma_\mathfrak{a}^+,\frac{\langle\mu,\alpha\rangle}{\langle\alpha,\alpha\rangle}\in \N\}$. Then the map which sends a representation to its highest weight is a bijection from $\widehat{G}_K$ to $\Lambda$.
\end{theorem}

\begin{remark}
\label{cartanhelgremark}In \cite{helgason2000groups}, the result is stated for $G$ simply connected. However, let $\Tilde{G}$ be the universal cover of $G$ and $p:\Tilde{G}\twoheadrightarrow G$ the covering map. Then $\ker(p)\subset Z(\Tilde{G})$. Let $\Tilde{K}=p^{-1}(K)$. We have that $\Tilde{G}/\Tilde{K}\simeq M$. Since $\Tilde{G}$ is connected and $M$ is simply connected, $\Tilde{K}$ is connected by the long exact sequence of homotopy groups. $\Tilde{K}$ contains $\Tilde{G}^\sigma$ which is connected (by \cite[Ch.V, Thm 3.3]{borel1998semisimple}) and has the same Lie algebra, thus $\Tilde{K}=\Tilde{G}^\sigma$.

This means that $(\Tilde{G},\Tilde{K})$ is another symmetric pair for the symmetric space $G/K$, so by Lemma \ref{indepgelfandpair}, there is a bijection between $\widehat{G}_K$ and $\widehat{\Tilde{G}}_{\Tilde{K}}$. Finally, $\Lambda$ depends only on the Lie algebra, and so by the theorem and the previous bijection, $\widehat{G}_K\to \Lambda$ is a bijection too.

Also, note that $\ker p\subset \Tilde{G}^\sigma$ since its elements act trivially on $M$. Hence, $$\ker p\subset Z(\Tilde{G})^\sigma=\{g\in Z(\Tilde{G})\vert \sigma(g)=g\}=\Tilde{K}\cap Z(\Tilde{G}).$$
\end{remark}

For example, if $\rank M=\rank G$, $\mathfrak{a}=\mathfrak{h}$ so $\Lambda$ is twice the set of dominant analytically integral elements. In general, let $\ell$ be the rank of $M$. Then the choice of $\Sigma_\mathfrak{a}^+$ gives a choice of a basis of the root system $\{\alpha_1,\cdots,\alpha_\ell\}$. By \cite[Thm. 2.1]{vretare_orth}, there are fundamental weights $\mu_i,1\leq i\leq \ell$ such that $\Lambda=\left\{\sum m_i\mu_i, m_i\in \N\right\}\simeq \N^\ell$. These fundamental weights verify $\langle\mu_i,\alpha_j\rangle=0$ if $i\neq j$.

Let $\mathfrak{a}_r=\{H\in \mathfrak{a} \vert \forall \lambda\in \Sigma_\mathfrak{a}, \lambda(H)\not\in i\pi\Z \}$ and $Q$ the connected component of $\mathfrak{a}_r$ contained in the positive Weyl chamber $\mathcal{C}=\{H\in \mathfrak{a}\vert \forall \lambda\in \Sigma_\mathfrak{a}^+, -i\lambda(H)>0\}$ and such that $0\in \overline{Q}$. Then in \cite[Prop. 3.2]{clerc}, we have the following result:
\begin{prop}
\label{kakcpt}For any $g\in G$, there exists $k_1,k_2\in K$ and a unique $H\in \overline{Q}$ such that $g=k_1\exp(H)k_2^{-1}$.
\end{prop}
\begin{remark}Again, this result is given for $G$ simply connected. If $G$ is not simply connected, note first that $Q$ depends only on the Lie algebra. Let $\Tilde{G}$ be the universal cover of $G$ and $p:\Tilde{G}\twoheadrightarrow G$ the covering map. With the notation of Remark \ref{cartanhelgremark}, for $g\in G$, there is $\Tilde{g}\in \Tilde{G}$ such that $g=p(\Tilde{g})$. The decomposition gives $k_1,k_2\in \Tilde{K}$ and $H\in Q$ such that $$\Tilde{g}=k_1\exp_{\Tilde{G}}(H)k_2^{-1},$$so $$g=p(k_1)p(\exp_{\Tilde{G}}(H))p(k_2)^{-1}=p(k_1)\exp_G(H)p(k_2)^{-1}$$gives the decomposition for $G$.
Furthermore, if there are $H,H'$ such that $\exp_G(H)=\exp_G(H')$, then $\exp_{\Tilde{G}}(H)=\exp_{\Tilde{G}}(H')k$, $k\in \ker p\subset \Tilde{K}$. By uniqueness for simply connected groups, $H=H'$.
\end{remark}

\begin{definition}
    \label{def:regularpoints}Let $G_1=K\exp(Q)K$, then $G_1$ is a dense open subset of $G$. We say that a point $g\in G_1$ is regular and $g\in G\setminus G_1$ is singular. 
\end{definition}

Now, let $(\varphi_\mu)_{\mu\in \Lambda}$ be the family of spherical functions of $(G,K)$. Consider the projection map $$\fonction{\pr}{G}{K\backslash G/K}{g}{KgK}.$$ The functions $\varphi_\mu$ are $K$-bi-invariant, thus the value of $\varphi_\mu(g)$ depends only on the double coset $\pr(g)=KgK$. Let $\psi_\mu=\varphi_\mu \circ \exp\vert_Q$. The family $(\psi_\mu)$ is a family of functions defined on a open subset of $\R^{\ell}$. We are interested in the regularity of these functions. Let $r_{opt}(M)= \sup \left\lbrace (r,\alpha) \vert (\varphi_\mu) \textrm{ bounded in } C^{(r,\alpha)}(Q)\right\rbrace$. By Lemma \ref{roptid}, this is also the supremum of $(r,\alpha)$ such that all $K$-bi-invariant matrix coefficients of unitary representations of $G$ are in $C^{(r,\alpha)}(Q)$.
 
\begin{remark}Since the spherical functions can be defined on $K\backslash G/K$, we can see them as functions on $\overline{Q}$ by the previous proposition without losing information on the function. However, we restrict to $Q$ because on the singular points, the behaviour cannot be controlled.
Furthermore, Proposition \ref{kakcpt} will be refined in Proposition \ref{regkak}, which will allow to recover the regularity obtained on the Lie algebra at the level of the group itself.
\end{remark}

Let $M_1,M_2$ be two simply connected symmetric spaces of compact type and $(G_1,K_1),(G_2,K_2)$ the associated canonical Gelfand pairs. Then $M=M_1\times M_2$ is a simply connected symmetric space of compact type, whose canonical Gelfand pair is $(G,K)$ with $G=G_1\times G_2$ and $K=K_1\times K_2$. Let $\pi$ be an irreducible representation of $G$. The irreducible representations of $G$ are $\pi=\pi_1\otimes \pi_2$ where $\pi_i$ is an irreducible representation of $G_i$, and $\pi$ has a non-zero $K$-invariant vector if and only if $\pi_i$ has a non-zero $K_i$-invariant vector, $i=1,2$. In that case, we know that the space of $K_i$-invariant vectors is one-dimensional, let $e_{K_i}$ be a unitary generator. Then, $e_K=e_{K_1}\otimes e_{K_2}$ is non-zero, unitary, $K$-invariant and generates the one-dimensional space of $K$-invariant vectors of $\pi$. Thus, if $\varphi_i$ is the spherical function associated to $\pi_i$, and $\varphi$ associated to $\pi$, we get $$\varphi(g_1,g_2)=\varphi_1(g_1)\varphi_2(g_2).$$Similarly, at the level of the Lie algebra, we have $Q=Q_1\times Q_2$ so $$\psi(H_1,H_2)=\psi_1(H_1)\psi_2(H_2).$$Denote $E_i$ vector space such that $Q_i\subset E_i$, and on $E_1\times E_2$, we consider the norm $\Vert (x,y)\Vert= \max(\Vert x\Vert,\Vert y\Vert)$. We write $\Lambda,\Lambda_1,\Lambda_2$ the set of highest weights of representations with invariant vectors for $M,M_1,M_2$. Since the constant function $1$ is a spherical function of any pair, we can take $\psi_2=1$ and we get that for any $\mu\in \Lambda_1$, the function $(g_1,g_2)\mapsto \psi_\mu(g_1)$ is a spherical function of $(G,K)$.

Let $(r,\alpha)<r_{opt}(M)$. Then $(\psi_\mu)_{\mu\in \Lambda_1}\subset (\psi_\mu)_{\mu\in\Lambda}$ is bounded in $C^{(r,\alpha)}(Q_1)$, thus we get $r_{opt}(M_1)\geq r_{opt}(M)$. Symmetrically, $r_{opt}(M_2)\geq r_{opt}(M)$.

Conversely, let $(r,\alpha)<\min(r_{opt}(M_1),r_{opt}(M_2))$. Let $L$ be a compact subset of $Q$. There are $L_1,L_2$ compact subsets of $Q_1,Q_2$ such that $L\subset L_1\times L_2$. For $k\leq r$, let $$C_{L_i,k}=\underset{\mu\in \Lambda_i}{\sup}\underset{x\in L_i}{\sup}\Vert D^k\psi_\mu(x)\Vert$$and let $$C_{L_i}=\underset{\mu\in \Lambda_i}{\sup}\underset{x\in L_i}{\sup}\frac{\Vert D^rf(x)-D^rf(y)\Vert}{\Vert x-y\Vert^\alpha}.$$These are finite numbers since by assumptions, $(\psi_\mu)_{\mu\in \Lambda_i}$ is bounded in $C^{(r,\alpha)}(Q_i)$, $i=1,2$. We have $\Lambda\simeq \Lambda_1\times \Lambda_2$. Consider $(\mu_1,\mu_2)\in \Lambda$, and $$\psi_{(\mu_1,\mu_2)}:(x_1,x_2)\mapsto \psi_{\mu_1}(x_1)\psi_{\mu_2}(x_2).$$Then clearly, $\psi_{(\mu_1,\mu_2)}$ is $r$ times differentiable and we have 
    \begin{multline*}
        D^r\psi_{(\mu_1,\mu_2)}(x_1,x_2)((H_1,K_1),\cdots,(H_r,K_r))=\\ \sum_{i\in I} D^{k_i}\psi_{\mu_1}(x_1)(H_{j_1(i)},\cdots,H_{j_{k_i}(i)}) D^{r-k_i}\psi_{\mu_2}(x_2)(K_{j_{k_i+1}(i)},\cdots,K_{j_r(i)})
    \end{multline*}
where $k_i\leq r$ for each $i\in I$ and $I$ finite. 

Thus, for all $(x_1,x_2),(y_1,y_2)\in L$, we have \begin{multline*}
    \Vert D^r\psi_{(\mu_1,\mu_2)}(x_1,x_2)-D^r\psi_{(\mu_1,\mu_2)}(y_1,y_2)\Vert\leq \\ \sum_{i\in I} \Vert D^{k_i}\psi_{\mu_1}(x_1) D^{r-{k_i}}\psi_{\mu_2}(x_2)-D^{k_i}\psi_{\mu_1}(y_1) D^{r-{k_i}}\psi_{\mu_2}(y_2)\Vert.
\end{multline*}

So we have \begin{align*}
    \MoveEqLeft[8]\Vert D^k\psi_{\mu_1}(x_1) D^{r-k}\psi_{\mu_2}(x_2)-D^k\psi_{\mu_1}(y_1) D^{r-k}\psi_{\mu_2}(y_2)\Vert & \\
     & \leq \Vert D^k\psi_{\mu_1}(x_1) D^{r-k}\psi_{\mu_2}(x_2)-D^k\psi_{\mu_1}(y_1) D^{r-k}\psi_{\mu_2}(x_2) \\
     & \phantom{\leq} +D^k\psi_{\mu_1}(y_1) D^{r-k}\psi_{\mu_2}(x_2)-D^k\psi_{\mu_1}(y_1) D^{r-k}\psi_{\mu_2}(y_2)\Vert\\
     & \leq \Vert D^k\psi_{\mu_1}(x_1)-D^k\psi_{\mu_1}(y_1)\Vert \Vert D^{r-k}\psi_{\mu_2}(x_2)\Vert \\
     & \phantom{\leq} + \Vert D^k\psi_{\mu_1}(y_1)\Vert \Vert D^{r-k}\psi_{\mu_1}(x_2)-D^{r-k}\psi_{\mu_1}(y_2)\Vert
\end{align*}

If $k\neq 0,r$, this gives \begin{multline*}\Vert D^k\psi_{\mu_1}(x_1) D^{r-k}\psi_{\mu_2}(x_2)-D^k\psi_{\mu_1}(y_1) D^{r-k}\psi_{\mu_2}(y_2)\Vert\\ \leq C_{L_1,k+1}C_{L_2,r-k}\Vert x_1-y_1\Vert +C_{L_1,k}C_{L_2,r-k+1}\Vert x_2-y_2\Vert,\end{multline*}if $k=0$, \begin{multline*}
    \Vert D^k\psi_{\mu_1}(x_1) D^{r-k}\psi_{\mu_2}(x_2)-D^k\psi_{\mu_1}(y_1) D^{r-k}\psi_{\mu_2}(y_2)\Vert \\ \leq C_{L_1,1}C_{L_2,r}\Vert x_1-y_1\Vert +C_{L_1,0}C_{L_2}\Vert x_2-y_2\Vert^\alpha
\end{multline*}and if $k=r$, \begin{multline*}\Vert D^k\psi_{\mu_1}(x_1) D^{r-k}\psi_{\mu_2}(x_2)-D^k\psi_{\mu_1}(y_1) D^{r-k}\psi_{\mu_2}(y_2)\Vert \\ \leq C_{L_1}C_{L_2,0}\Vert x_1-y_1\Vert^\alpha +C_{L_1,r}C_{L_2,1}\Vert x_2-y_2\Vert.\end{multline*}

But since $\alpha\in [0,1]$, there is $C_i>0$ such that for all $x,y\in L_i$, $$\Vert x-y\Vert \leq C_i \Vert x-y\Vert^\alpha$$and so there is a constant $C_L>0$ which does not depend on $(\mu_1,\mu_2)\in \Lambda$, such that 
\begin{alignat*}{1}
\Vert D^r\psi_{(\mu_1,\mu_2)}(x_1,x_2)-D^r\psi_{(\mu_1,\mu_2)}(y_1,y_2)\Vert &\leq  C_L\max(\Vert x_1-y_1\Vert^\alpha,\Vert x_2-y_2\Vert^\alpha)\\&=C_L\Vert (x_1,x_2)-(y_1,y_2)\Vert^\alpha.
\end{alignat*}
And thus, we showed that $(\psi_\mu)_{\mu\in\Lambda}$ is bounded in $C^{(r,\alpha)}(Q)$, so that $$r_{opt}(M)\geq \min(r_{opt}(M_1),r_{opt}(M_2)).$$ By induction, we get the following :
\begin{prop}
Let $M_i$, $1\leq i\leq n$ be a simply connected symmetric space of compact type, and $M=\prod_{i=1}^n M_i$. Then $$r_{opt}(M)=\min(r_{opt}(M_i)).$$
\end{prop}

This result tells us that we can study only the irreducible simply connected symmetric spaces of compact type. By the classification in \cite[Ch. VII]{loos1969symmetric2}, such a space $M$ is either \begin{enumerate}
    \item one of the spaces in \cite[Ch. VII, Table 4 and Table 8]{loos1969symmetric2}, 
    \item $(G\times G)/\Delta(G)$, where $G$ is a simply connected simple compact Lie group and $\Delta(G)$ the diagonal subgroup of $G\times G$.
\end{enumerate}
The first case contains the symmetric spaces of rank $1$ for which we will solve the question in Section \ref{sectionrank1}. In Section \ref{caracteres}, we will study $(G\times G)/\Delta(G)$ (even for reducible and not simply connected). Finally, we give some partial results for some other higher rank symmetric spaces in Section \ref{sectionpartialresult}.

\section{Symmetric spaces of rank \texorpdfstring{$1$}{1}}\label{sectionrank1}
\subsection{Spherical functions of symmetric spaces of rank 1}\label{sectionspheriques}
A classification of symmetric spaces of compact type can be found in \cite[Ch. VII]{loos1969symmetric2} and from this classification, we extract the canonical compact Gelfand pairs associated to compact symmetric spaces of rank $1$. Table \ref{tab:rank1} lists these symmetric spaces with $(G,K)$ the canonical pair associated to $M$, the dimension of $M$ and two real parameters $\alpha,\beta$ used later.

\begin{table}[ht]
\renewcommand{\arraystretch}{1.3}
$$\begin{array}{|c|c|c|c|c|c|}
     \hline 
     M  & G & K & \dim M & \alpha & \beta \\
     \hline
     S^{k-1},k\geq 3 & SO(k) & SO(k-1) & k-1 & \frac{k-3}{2} &  \frac{k-3}{2} \\
     \hline
     \R P^{k-1},k\geq 3 & SO(k) & S(O(1)\times O(k-1)) & k-1 &  \frac{k-3}{2} &  -\frac{1}{2}\\
     \hline
     \C P^{k-1},k\geq 2 & SU(k) & S(U(1)\times U(k-1)) & 2(k-1) & k-2& 0\\
    \hline
    \mathbb{H} P^{k-1},k\geq 2 & Sp(k) & Sp(1)\times Sp(k-1) & 4(k-1) & 2k-3 & 1\\
    \hline
    \mathbb{O}P^2 & F_4 & \mathrm{Spin}(9) & 16 &  7& 3 \\
    \hline
     \end{array}$$
\caption{Compact symmetric pairs of rank $1$}
\label{tab:rank1}
\end{table}
\begin{remark}We can see that $\alpha=\frac{\dim G/K}{2}-1$.\end{remark}

The spherical functions of these pairs are well-known and can be found in \cite[Ch. V, Theorem 4.5]{helgason2000groups}, while the dimension of the associated representation are found in \cite[Theorem 2.4, 3.2, 4.2, 5.2 and 6.2]{wolfcahn}. They can be expressed in terms of Jacobi polynomials.
\begin{definition}[Jacobi polynomials] Let $\alpha > -1, \beta > -1$, the Jacobi polynomials of parameters $(\alpha,\beta)$ are defined as the unique polynomials $\left(P_n^{(\alpha,\beta)}\right)_{n\in \N}$ such that for all $n\in \N$, $P_n^{(\alpha,\beta)}$ is of degree $n$, for all $m\neq n$, $$\int_{-1}^1 P^{(\alpha,\beta)}_m(x)P_n^{(\alpha,\beta)}(x)(1-x)^\alpha (1+x)^\beta \,dx=0$$and for all $n\in \N$, $$P_n^{(\alpha,\beta)}(1)=\binom{n+\alpha}{n}.$$
\end{definition}

Since $G/K$ is of rank $1$, we have either $\Sigma_{\mathfrak{a}}^+=\{\alpha\}$ or $\Sigma_{\mathfrak{a}}^+=\{\alpha,2\alpha\}$. In both cases, $\mathfrak{a}\simeq \R$ by $H\mapsto -i\alpha(H)$. By this identification, we have $Q\simeq ]0,\pi[$. If $\varphi$ is a spherical function of $(G,K)$, denote $\psi=\varphi\circ \exp\vert_{\overline{Q}}$.
\begin{theorem}If $(G,K)$ is a compact symmetric pair of rank $1$, then its spherical functions are the functions $\varphi_n$ defined at the level of the Lie algebra by $$\psi_n:\theta\in \overline{Q}\mapsto \frac{\Gamma(\alpha+1)\Gamma(n+1)}{\Gamma(n+\alpha+1)}P_n^{(\alpha,\beta)}(\cos\theta).$$

Furthermore, the dimension $m_n$ of the representation associated to $\varphi_n$ is a polynomial in $n$ of degree $(\dim G/K)-1$.
\end{theorem}
\begin{remark}
To get the function $\varphi_n$ itself, we need to understand the projection $H:G\mapsto \overline{Q}$ since $\varphi_n=\psi_n\circ H$ by $K$-bi-invariance.

For $(SO(n),SO(n-1))$, we have $H(g)=\arccos(g_{1,1})$.
For the other non-exceptional rank $1$ pairs, we have $H(g)=\arccos({2\vert g_{1,1}\vert^2-1})$.

We delay this study at the level of the group until Section \ref{kfinitesection}. By convention, "the family of spherical functions" will refer to the functions $(\psi_n)_{n\in \N}$.
\end{remark}

\subsection{Regularity of matrix coefficients}\label{mainsection}
\begin{theorem}
\label{mainthm}Let $(G,K)$ be one of the symmetric pairs of rank $1$. Let $$\alpha_\infty=\frac{\dim G/K-1}{2}.$$Then the family of spherical functions of the Gelfand pair $(G,K)$ is bounded in $C^{(\lfloor \alpha_\infty\rfloor,\alpha_\infty-\lfloor \alpha_\infty\rfloor)}(Q)$.
\end{theorem}
Given that the spherical functions of the Gelfand pairs we are interested in are all Jacobi polynomials, of parameters $(\alpha,\beta)$ fixed by the pair, we can derive the theorem from the following result:
\begin{theorem}
 \label{jacobi} Let $\alpha\geq 0,\beta >-1$ be two reals. Then the family $$\left(\frac{\Gamma(\alpha+1)\Gamma(n+1)}{\Gamma(n+\alpha+1)}P_n^{(\alpha,\beta)}\right)_{n\in \N}$$ is bounded in $C^{\left(\lfloor \alpha+\frac{1}{2}\rfloor,\alpha+\frac{1}{2}-\lfloor \alpha+\frac{1}{2}\rfloor\right)}(]-1,1[)$.   
\end{theorem}

\begin{proof}[Proof of Theorem \ref{mainthm} using Theorem \ref{jacobi}]According to Section \ref{sectionspheriques} and up to a reparametrization by Lemma \ref{precomposition}, the spherical functions of $(G,K)$ are normalized Jacobi polynomials of parameters $(\alpha,\beta)$ with $$\alpha=\frac{\dim G/K}{2}-1.$$Since $$\alpha+\frac{1}{2}=\frac{\dim G/K-1}{2}=\alpha_\infty,$$we get the result by Theorem \ref{jacobi}.
\end{proof}

The proof of Theorem \ref{jacobi} relies on two ingredients that can be found in \cite[(4.21.7) and (8.21.10)]{szego1939orthogonal}.
\begin{prop}
\label{derivee}For all $\alpha,\beta\geq 0$, and for all $n\geq k$, we have $$\frac{d^k}{dx^k}P_n^{(\alpha,\beta)}(x)=\frac{\Gamma(\alpha+\beta+n+1+k)}{2^k\Gamma(\alpha+\beta+n+1)}P_{n-k}^{(\alpha+k,\beta+k)}(x).$$
\end{prop}

\begin{prop}[Darboux's formula]
\label{darboux}We have $$P_n^{(\alpha,\beta)}(\cos \theta)=n^{-\frac{1}{2}}k(\theta)\cos(N\theta+\gamma)+O(n^{-\frac{3}{2}})$$where $$\left\lbrace{\begin{aligned}k(\theta )&=\pi ^{{-{\frac  {1}{2}}}}\sin ^{{-\alpha -{\frac  {1}{2}}}}{\tfrac  {\theta }{2}}\cos ^{{-\beta -{\frac  {1}{2}}}}{\tfrac  {\theta }{2}},\\N&=n+{\tfrac  {1}{2}}(\alpha +\beta +1),\\\gamma &=-{\tfrac  {\pi }{2}}\left(\alpha +{\tfrac  {1}{2}}\right),\end{aligned}}\right.$$and $O(n^{-\frac{3}{2}})$ is uniform for $\theta$ in the compact $[\varepsilon,\pi-\varepsilon]$, for any $\varepsilon>0$.
\end{prop}

\begin{proof}
[Proof of Theorem \ref{jacobi}]Let $\varphi_n(x)=\frac{\Gamma(\alpha+1)\Gamma(n+1)}{\Gamma(n+\alpha+1)}P_n^{(\alpha,\beta)}(x)$. Set $L$ a compact subset of $]-1,1[$. Up to replacing $L$ by its convex hull we can assume that $L$ is an interval. There exists $\varepsilon>0$ such that if $\cos \theta\in L$ with $\theta\in [0,\pi]$, then $\theta \in [\varepsilon,\pi-\varepsilon]$. Thus by Proposition \ref{darboux}, because the function $\theta\mapsto k(\theta)$ of Darboux's formula is bounded on compacts in $]0,\pi[$, there exists $C_{L,\alpha,\beta}>0$ such that for all $x\in L$ and $n\in \N^*$, \begin{equation}\label{majjacobi}\vert P_n^{(\alpha,\beta)}(x)\vert \leq C_{L,\alpha,\beta}n^{-\frac{1}{2}}.\end{equation}Thus, by Proposition \ref{derivee}, we have for all $x\in L$ and $n> k$, \begin{equation}
\label{majderiv}\left\vert \frac{d^k}{dx^k}P_n^{(\alpha,\beta)}(x)\right\vert \leq \frac{\Gamma(\alpha+\beta+n+1+k)}{2^k\Gamma(\alpha+\beta+n+1)}C_{L,\alpha+k,\beta+k}(n-k)^{-\frac{1}{2}}.
\end{equation}For any $k\in \N$, and $n>k$, we have $n-k\geq \frac{n}{k+1}$. Using that $\Gamma(x+1)=x\Gamma(x)$ for $x>0$, we get that \begin{equation}
    \left\vert \frac{d^k}{dx^k}P_n^{(\alpha,\beta)}(x)\right\vert \leq \frac{\sqrt{k+1}\prod_{j=1}^k(\alpha+\beta+n+j)}{2^k}C_{L,\alpha+k,\beta+k}n^{-\frac{1}{2}},
\end{equation}so \begin{equation}\label{majphi}\begin{array}{rcl}
    \vert \varphi_n^{(k)}(x)\vert & = &  \frac{\Gamma(\alpha+1)\Gamma(n+1)}{\Gamma(n+\alpha+1)}\left\vert \frac{d^k}{dx^k}P_n^{(\alpha,\beta)}(x) \right\vert \\
     & \leq & \frac{C_{L,\alpha+k,\beta+k}\sqrt{k+1}\Gamma(\alpha+1)}{2^k}\frac{\Gamma(n+1)\prod_{j=1}^k(\alpha+\beta+n+j)}{\Gamma(n+\alpha+1)} n^{-\frac{1}{2}}.
\end{array}\end{equation}Now, by the Stirling approximation for the Gamma function (see \cite[Thm. 5]{gordon}), there are $C_1,C_2>0$ such that if $f(x)=\sqrt{2\pi}x^{x+\frac{1}{2}}e^{-x}$, then for $x> 0$, $$C_1f(x)\leq \Gamma(x+1)\leq C_2f(x).$$Thus for all $n> 0$, \begin{equation}\label{gamma}\frac{\Gamma(n+1)}{\Gamma(n+\alpha+1)}\leq \frac{C_2}{C_1}\frac{n^{n+\frac{1}{2}}e^{-n}}{(n+\alpha)^{n+\alpha+\frac{1}{2}}e^{-(n+\alpha)}}\leq \frac{C_2e^\alpha}{C_1}(n+\alpha)^{-\alpha}\leq \frac{C_2e^\alpha}{C_1}n^{-\alpha} \end{equation}since $\alpha>0$. Furthermore, for each $j\in \llbracket 1,k\rrbracket$ and $n>0$, we have $$(\alpha+\beta+n+j)\leq (\alpha+\beta+k+n)\leq (\alpha+\beta+k+1)n$$thus \begin{equation}
    \label{majcoef} \prod_{j=1}^k(\alpha+\beta+n+j)\leq (\alpha+\beta+k+1)^k n^k.
\end{equation}Setting $$\Tilde{C}_{L,\alpha,\beta,k}=\frac{C_{L,\alpha+k,\beta+k}\sqrt{k+1}\Gamma(\alpha+1)}{2^k}\frac{C_2e^\alpha}{C_1} (\alpha+\beta+k+1)^k$$and inserting (\ref{gamma}) and (\ref{majcoef}) into (\ref{majphi}), we finally get that for all $n>k$ and $x\in L$, \begin{equation}
    \label{majderivphi} \vert  \varphi_n^{(k)}(x)\vert \leq \Tilde{C}_{L,\alpha,\beta,k}n^{k-\alpha-\frac{1}{2}}.
\end{equation}

From this inequality, we see that the derivatives of the family of spherical functions are bounded in $n$ up to order $\lfloor \alpha+\frac{1}{2}\rfloor$. If $\alpha+\frac{1}{2}\in \Z$, this shows that $\left(\varphi_n\right)_{n\in \N}$ is bounded in $C^{(\alpha+\frac{1}{2},0)}$. If $\alpha+\frac{1}{2}\not\in \Z$, we must now verify the Hölder part.

There are now two cases to consider.

\textbf{Case 1:} assume that $\alpha-\lfloor \alpha\rfloor<\frac{1}{2}$, which is equivalent to $\lfloor \alpha+\frac{1}{2}\rfloor=\lfloor \alpha\rfloor$. Set $$\lambda=\alpha-\lfloor \alpha\rfloor+\frac{1}{2}=\alpha+\frac{1}{2}-\lfloor \alpha+\frac{1}{2}\rfloor.$$We have $\lambda \in [\frac{1}{2},1[$. Then by (\ref{majderivphi}) applied to $k=\lfloor \alpha\rfloor$, we get for $x\in L$ that $$\vert \varphi_n^{(\lfloor \alpha\rfloor)}(x)\vert \leq \Tilde{C}_{L,\alpha,\beta,\lfloor \alpha\rfloor} n^{-\lambda},$$so that for $x,y\in L$, \begin{equation}
    \label{case11}
\vert \varphi_n^{(\lfloor \alpha\rfloor)}(x)-\varphi_n^{(\lfloor \alpha\rfloor)}(y)\vert \leq \vert \varphi_n^{(\lfloor \alpha\rfloor)}(x)\vert+\vert \varphi_n^{(\lfloor \alpha\rfloor)}(y)\vert \leq 2 \Tilde{C}_{L,\alpha,\beta,\lfloor \alpha\rfloor}  n^{-\lambda}.\end{equation}

Applying (\ref{majderivphi}) to $k=\lfloor \alpha\rfloor+1$, we get for $x\in L$ that $$\vert \varphi_n^{(\lfloor \alpha\rfloor+1)}(x)\vert \leq \Tilde{C}_{L,\alpha,\beta,\lfloor \alpha\rfloor+1}  n^{1-\lambda},$$so that for $x,y\in L$, \begin{equation}
    \label{case12}
\vert \varphi_n^{(\lfloor \alpha\rfloor)}(x)-\varphi_n^{(\lfloor \alpha\rfloor)}(y)\vert \leq  \underset{t\in K}{\sup} \vert \varphi_n^{(\lfloor \alpha\rfloor+1)}(t)\vert \vert x-y\vert \leq \Tilde{C}_{L,\alpha,\beta,\lfloor \alpha\rfloor+1}  n^{1-\lambda}\vert x-y\vert.\end{equation}
Set $M=(2\Tilde{C}_{L,\alpha,\beta,\lfloor \alpha\rfloor})^{1-\lambda}\Tilde{C}_{L,\alpha,\beta,\lfloor \alpha\rfloor+1}^\lambda$. Combining (\ref{case11}) and (\ref{case12}) we have $$\vert \varphi_n^{(\lfloor \alpha\rfloor)}(x)-\varphi_n^{(\lfloor \alpha\rfloor)}(y)\vert \leq M(n^{-\lambda})^{1-\lambda} (n^{1-\lambda}\vert x-y\vert)^\lambda=M \vert x-y\vert^\lambda, $$ which is the result we wanted.

\textbf{Case 2}: assume that $\alpha-\lfloor \alpha\rfloor\geq\frac{1}{2}$, which is equivalent to $\lfloor \alpha+\frac{1}{2}\rfloor=\lfloor \alpha\rfloor+1$. Set $$\lambda=\alpha-\lfloor \alpha\rfloor-\frac{1}{2}=\alpha+\frac{1}{2}-\lfloor \alpha+\frac{1}{2}\rfloor.$$We have $\lambda \in [0,\frac{1}{2}]$. Similarly to the first case, we apply (\ref{majderivphi}) to $k=\lfloor \alpha\rfloor+1$ and $k=\lfloor \alpha\rfloor+2$ to get that for all $x,y\in L$, \begin{equation}
    \label{case21}
\vert \varphi_n^{(\lfloor \alpha\rfloor+1)}(x)-\varphi_n^{(\lfloor \alpha\rfloor+1)}(y)\vert \leq  2\Tilde{C}_{L,\alpha,\beta,\lfloor \alpha\rfloor+1} n^{-\lambda}\end{equation}and 
\begin{equation}
    \label{case22}
\vert \varphi_n^{(\lfloor \alpha\rfloor+2)}(x)-\varphi_n^{(\lfloor \alpha\rfloor+2)}(y)\vert \leq \Tilde{C}_{L,\alpha,\beta,\lfloor \alpha\rfloor+2}  n^{1-\lambda}\vert x-y\vert,\end{equation}
which combines as in the first case, giving$$\vert \varphi_n^{(\lfloor \alpha\rfloor)}(x)-\varphi_n^{(\lfloor \alpha\rfloor)}(y)\vert \leq M'(n^{-\lambda})^{1-\lambda}(n^{1-\lambda}\vert x-y\vert)^\lambda\leq M' \vert x-y\vert^\lambda, $$ with $M'=(2\Tilde{C}_{L,\alpha,\beta,\lfloor \alpha\rfloor+1})^{1-\lambda}\Tilde{C}_{L,\alpha,\beta,\lfloor \alpha\rfloor+2}^\lambda$, which is the result we wanted.
\end{proof}

\begin{coro}\label{matrixcoef}Let $(G,K)$ be a compact symmetric Gelfand pair of rank one. Let $\varphi$ be a $K$-bi-invariant matrix coefficient of a unitary representation of $G$. Then $\varphi\circ \exp\in C^{(\lfloor \alpha_\infty \rfloor,\alpha_\infty-\lfloor \alpha_\infty \rfloor)}(Q)$.
\end{coro}
\begin{proof}It follows directly from Lemma \ref{roptid}.
\end{proof}

\subsection{Schatten norm}\label{schattensection}
Given a Hilbert space $H$, $1\leq p \leq +\infty$ and $T$ an operator on $H$, the Schatten $p$-norm of $T$ is $$\Vert T\Vert_{S_p}=\Tr(\vert T\vert^p)^{1/p}$$defined by functional calculus. If $p=+\infty$, this is the operator norm. Then $S_p(H)$ is the space of operators $T$ such that $\Vert T\Vert_{S_p}<+\infty$.

Let $(G,K)$ be a pair as in Section \ref{sectionspheriques}. Define $\Tilde{T}_g=\int_{K\times K} \lambda(kgk')\,dk\,dk'$ where $dk$ is the normalized Haar measure on $K$ and $\lambda$ the regular representation of $G$ on $L^2(G)$. Since the map $\Tilde{T}$ is $K$-bi-invariant, the function $T=\Tilde{T}\circ \exp:\overline{Q}\mapsto B(L^2(G))$ determines $\Tilde{T}$. We want to investigate the regularity of $T$ as map into $S_p(L^2(G))$. Let $p>2+\frac{2}{\dim G/K-1}$ and set $$\alpha_p=\alpha_\infty-\frac{\dim G/K}{p}=\frac{\dim G/K-1}{2}-\frac{\dim G/K}{p}$$ (so $\alpha_p>0$).

\begin{prop}
\label{schatteninfty}The map $\delta\mapsto T_\delta$ belongs to $$C^{(\lceil \alpha_\infty\rceil-1,\alpha_\infty-\lceil \alpha_\infty\rceil+1)}\left(Q,S_\infty(L^2(G))\right).$$    
\end{prop}
\begin{proof}There is an orthonormal basis of $L^2(G)$ such that for any $\delta\in Q$, $T_\delta$ is diagonal with eigenvalues $\psi_n(\delta)$ of multiplicity $m_n$ the dimension of the representation associated to $\psi_n$.

By Theorem \ref{mainthm}, we know that $\varphi_n$ is bounded in $C^{(\lfloor \alpha_\infty\rfloor,\alpha_\infty-\lfloor \alpha_\infty\rfloor)}(Q)$, thus by Lemma \ref{existencederiv}, we get the result.
\end{proof}
\begin{remark}
If $\alpha=0$ in Lemma \ref{existencederiv}, we cannot show that $T$ is $C^r$. The derivative will exist in a weak sense, but we cannot show convergence in norm.

In particular, if $\alpha_\infty\in \Z$, we do not get that $\partial^{\alpha_\infty}T$ exists. In fact, we will see in the next subsection that the result in Proposition \ref{schatteninfty} is optimal.
\end{remark}

The following theorem is a generalisation of a result from \cite{dlSPR} on the group $SO(n)$.

\begin{theorem}
\label{schatten}Let $2+\frac{2}{\dim G/K-1}<p<+\infty$. The map $\delta\mapsto T_\delta$ belongs to $C^{(\lfloor\alpha\rfloor,\alpha-\lfloor\alpha\rfloor)}(Q,S_p(L^2(G)))$ where 
$$\alpha=\left\lbrace \begin{aligned}\alpha_p & \textrm{ if }\alpha_p\not\in \Z\\ \alpha_p-\varepsilon& \textrm{ if }\alpha_p\in \Z\end{aligned}\right.$$ with $\varepsilon>0$ arbitrarily small.    
\end{theorem}
\begin{proof}Let \begin{equation}\label{rappel}\varphi_n(x)=\frac{\Gamma(u+1)\Gamma(n+1)}{\Gamma(n+u+1)}P_n^{(u,v)}(x)\end{equation} be the spherical functions of $(G,K)$ viewed on $\cos(Q)=]-1,1[$ with $u=\frac{\dim G/K}{2}-1$. By Lemma \ref{precomposition}, we can assume that $T$ is defined on $]-1,1[$ and is such that there is an orthogonal basis such that $T_\delta$ is diagonal with eigenvalues $\varphi_n(\delta)$ of multiplicity $m_n$, the dimension of the associated representation (see Section \ref{sectionspheriques}). Notice that there is a constant $C>0$ which depends only on $(G,K)$ such that for all $n\in \N$, \begin{equation}\label{dimension}m_n\leq C(n+1)^{\dim G/K-1}.\end{equation}
Let $r<\alpha_\infty$ and $\delta\in ]-1,1[$. By Lemma \ref{existencederiv}, we have that $\partial^r T$ exists and $$\Vert \partial^r T_\delta\Vert_{S_p}=\left(\sum_{n\geq 0} m_n \vert \varphi_n^{(r)}(\delta)\vert^p\right)^{1/p}.$$Thus, using the inequality (\ref{majderivphi}) from the proof of Theorem \ref{jacobi}, we have that there is a constant $C(G,r,\delta)>0$ such that \begin{multline*}m_n\vert \varphi_n^{(r)}(x)\vert^p \leq C(G,r,\delta)(n+1)^{\dim G/K-1+p(r-\frac{\dim G/K-1}{2})}\\=C(G,r,\delta)(n+1)^{p\left(r-\left(\frac{\dim G/K-1}{2}-\frac{\dim G/K}{p}\right)\right)-1}.\end{multline*}Thus the sum converges as soon as $r<\frac{\dim G/K-1}{2}-\frac{\dim G/K}{p}=\alpha_p$, and we get $\partial^r T_\delta\in S_p(L^2(G)))$ for these $r$. This happens when $r=\lfloor \alpha\rfloor$.

In that case, for $\delta,\delta'$ in some compact $I$ of $]-1,1[$, which we can assume to be an interval up to taking its convex hull, we have $$\Vert \partial^{\lfloor \alpha\rfloor}T_\delta - \partial^{\lfloor \alpha\rfloor}T_{\delta'}\Vert_{S_p}^p=\sum_{n=0}^\infty m_n \vert \varphi_n^{(\lfloor \alpha\rfloor)}(\delta)-\varphi_n^{(\lfloor \alpha\rfloor)}(\delta')\vert^p.$$We bound again $m_n$ by $C(n+1)^{\dim G/K-1}$. For the term $\vert \varphi_n^{(\lfloor \alpha\rfloor)}(\delta)-\varphi_n^{(\lfloor \alpha\rfloor)}(\delta')\vert$, we have two different bounds. By the expression of $\varphi_n$ in (\ref{rappel}) and using (\ref{majderivphi}), on the one hand, there is $C_1>0$ such that\begin{align*}
    \vert \varphi_n^{(\lfloor \alpha\rfloor)}(\delta)-\varphi_n^{(\lfloor \alpha\rfloor)}(\delta')\vert &\leq \vert \varphi_n^{(\lfloor \alpha\rfloor)}(\delta)\vert+\vert \varphi_n^{(\lfloor \alpha\rfloor)}(\delta')\vert \\&\leq C_1(n+1)^{\lfloor \alpha\rfloor-u-\frac{1}{2}} \\& \leq C_1(n+1)^{\lfloor \alpha\rfloor-\frac{\dim G/K-1}{2}},
\end{align*}which we use as soon a $n\vert \delta-\delta'\vert >1$ and on the other hand, there is $C_2>0$ such that\begin{align*}
    \vert \varphi_n^{(\lfloor \alpha\rfloor)}(\delta)-\varphi_n^{(\lfloor \alpha\rfloor)}(\delta')\vert & \leq \underset{t\in I}{\sup }\vert \varphi_n^{(\lfloor \alpha\rfloor+1)}(t)\vert \vert \delta-\delta'\vert \\ &\leq C_2(n+1)^{\lfloor \alpha\rfloor+1-\frac{\dim G/K-1}{2}}\vert \delta-\delta'\vert,
\end{align*} which we use a soon as $n\vert \delta-\delta'\vert \leq 1$. So we have, since $$\dim G/K-1+p\left(\lfloor \alpha \rfloor-\frac{\dim G/K-1}{2}\right)=p(\lfloor \alpha\rfloor-\alpha_p)-1,$$ that \begin{align*}\Vert \partial^{\lfloor \alpha\rfloor}T_\delta - \partial^{\lfloor \alpha\rfloor}T_{\delta'}\Vert_{S_p}^p &\leq CC_1^p\sum_{n> \frac{1}{\vert \delta-\delta'\vert}}(n+1)^{p(\lfloor \alpha\rfloor-\alpha_p)-1}\\& \phantom{\leq}+CC_2^p\sum_{n\leq \frac{1}{\vert \delta-\delta'\vert}} (n+1)^{p(\lfloor \alpha\rfloor+1-\alpha_p)-1}\vert \delta-\delta'\vert^p.\end{align*}

There are now two cases to consider.

\textbf{Case 1:} if $\alpha_p\not\in \Z$, we have $\alpha=\alpha_p$. Let $n_0=\lfloor \frac{1}{\vert \delta-\delta'\vert}\rfloor+1$. We get \begin{align*}
    \Vert \partial^{\lfloor \alpha\rfloor}T_\delta - \partial^{\lfloor \alpha\rfloor}T_{\delta'}\Vert_{S_p}^p & \leq  CC_1^p\displaystyle \sum\limits_{n=n_0}^{+\infty} \frac{1}{(n+1)^{1+p(\alpha-\lfloor \alpha\rfloor)}} \\ &\phantom{\leq } +CC_2^{p}\vert \delta-\delta'\vert^p \sum\limits_{n=0}^{n_0-1} \frac{1}{(n+1)^{1+p(\alpha-\lfloor \alpha\rfloor-1)}}.
\end{align*} There are again two subcases to consider. First, assume that $1+p(\alpha-\lfloor \alpha\rfloor-1)\geq 0$, then 
\begin{multline*}
     \Vert \partial^{\lfloor \alpha\rfloor}T_\delta - \partial^{\lfloor \alpha\rfloor}T_{\delta'}\Vert_{S_p}^p \\
     \begin{aligned}
         & \leq CC_1^p \sum\limits_{n=n_0}^{+\infty} \frac{1}{(n+1)^{1+p(\alpha-\lfloor \alpha\rfloor)}}+CC_2^{p}\vert \delta-\delta'\vert^p \sum\limits_{n=0}^{n_0-1} \frac{1}{(n+1)^{1+p(\alpha-\lfloor \alpha\rfloor-1)}}\\
        & \leq  CC_1^p \int_{n_0}^\infty \frac{1}{x^{1+p(\alpha-\lfloor \alpha\rfloor)}}dx + CC_2^p \vert \delta-\delta'\vert^p \int_0^{n_0}\frac{1}{x^{1+p(\alpha-\lfloor \alpha\rfloor-1)}}dx \\
        & \leq \frac{CC_1^p}{p(\alpha-\lfloor \alpha\rfloor) n_0^{p(\alpha-\lfloor \alpha\rfloor)}}+  \frac{CC_2  \vert \delta-\delta'\vert^p}{-p(\alpha-\lfloor \alpha\rfloor-1)}n_0^{p(\lfloor \alpha\rfloor+1-\alpha)}\\
        & \leq \frac{CC_1^p}{p(\alpha-\lfloor \alpha\rfloor)}\vert \delta-\delta'\vert^{p(\alpha-\lfloor \alpha\rfloor)} + \frac{CC_2^p \vert \delta-\delta'\vert^p}{p(\lfloor \alpha\rfloor+1-\alpha)}\left(\frac{1+\vert \delta-\delta'\vert}{\vert\delta-\delta'\vert}\right)^{p(\lfloor \alpha\rfloor+1-\alpha)} \\
        & \leq \frac{CC_1^p}{p(\alpha-\lfloor \alpha\rfloor)}\vert \delta-\delta'\vert^{p(\alpha-\lfloor \alpha\rfloor)} + \frac{CC_2^p (1+\pi)^{p(\lfloor \alpha\rfloor+1-\alpha)}}{p(\lfloor \alpha\rfloor+1-\alpha)}\left(\vert \delta-\delta'\vert\right)^{p(\alpha-\lfloor \alpha\rfloor)}\\
        & \leq \Tilde{C}^p \vert \delta-\delta'\vert^{p(\alpha-\lfloor \alpha\rfloor)}
     \end{aligned}
\end{multline*}where $\Tilde{C}>0$ is a constant which do not depend on $\delta,\delta'\in I$. Thus, we finally have $$\Vert \partial^{\lfloor \alpha\rfloor}T_\delta - \partial^{\lfloor \alpha\rfloor}T_{\delta'}\Vert_{S_p}\leq \Tilde{C}\vert \delta-\delta'\vert^{\alpha-\lfloor \alpha\rfloor}.$$
In the other subcase, we have $1+p(\alpha-\lfloor \alpha\rfloor-1)\leq 0$, then \begin{multline*}
   \Vert \partial^{\lfloor \alpha\rfloor}T_\delta - \partial^{\lfloor \alpha\rfloor}T_{\delta'}\Vert_{S_p}^p \\
   \begin{aligned}
    & \leq CC_1^p \sum\limits_{n=n_0}^{+\infty} \frac{1}{(n+1)^{1+p(\alpha-\lfloor \alpha\rfloor)}}+CC_2^{p}\vert \delta-\delta'\vert^p \sum\limits_{n=0}^{n_0-1}(n+1)^{-1+p(\lfloor \alpha\rfloor+1-\alpha)}  \\
     & \leq  CC_1^p \int_{n_0}^\infty \frac{1}{x^{1+p(\alpha-\lfloor \alpha\rfloor)}}dx + CC_2^p \vert \delta-\delta'\vert^p \int_1^{n_0+1}x^{-1+p(\lfloor \alpha\rfloor+1-\alpha)}dx\\
     & \leq  \frac{CC_1^p}{p(\alpha-\lfloor \alpha\rfloor) n_0^{p(\alpha-\lfloor \alpha\rfloor)}}+  \frac{CC_2  \vert \delta-\delta'\vert^p}{-p(\alpha-\lfloor \alpha\rfloor-1)}\left((n_0+1)^{p(\lfloor \alpha\rfloor+1-\alpha)}-1\right)\\
     & \leq \frac{CC_1^p}{p(\alpha-\lfloor \alpha\rfloor)}\vert \delta-\delta'\vert^{p(\alpha-\lfloor \alpha\rfloor)} + \frac{CC_2^p \vert \delta-\delta'\vert^p}{p(\lfloor \alpha\rfloor+1-\alpha)}\left(\frac{1+2\vert \delta-\delta'\vert}{\vert\delta-\delta'\vert}\right)^{p(\lfloor \alpha\rfloor+1-\alpha)} \\
     & \leq  \frac{CC_1^p}{p(\alpha-\lfloor \alpha\rfloor)}\vert \delta-\delta'\vert^{p(\alpha-\lfloor \alpha\rfloor)} + \frac{CC_2^p (1+2\pi)^{p(\lfloor \alpha\rfloor+1-\alpha)}}{p(\lfloor \alpha\rfloor+1-\alpha)}\left(\vert \delta-\delta'\vert\right)^{p(\alpha-\lfloor \alpha\rfloor)}\\
     & \leq  \Tilde{C'}^p \vert \delta-\delta'\vert^{p(\alpha-\lfloor \alpha\rfloor)}
\end{aligned}\end{multline*}where $\Tilde{C'}>0$ does not depend on $\delta,\delta'\in I$. Again, we finally get $$\Vert \partial^{\lfloor \alpha\rfloor}T_\delta - \partial^{\lfloor \alpha\rfloor}T_{\delta'}\Vert_{S_p}\leq \Tilde{C'}\vert \delta-\delta'\vert^{\alpha-\lfloor \alpha\rfloor}.$$

\textbf{Case 2:} if $\alpha_p\in \Z$, let $0<\varepsilon<1$. We have $\lfloor \alpha\rfloor=\lfloor \alpha_p-\varepsilon\rfloor=\alpha_p-1$. Let $n_0=\lfloor \frac{1}{\vert \delta-\delta'\vert}\rfloor+1$. We get \begin{multline*}
   \Vert \partial^{\lfloor \alpha\rfloor}T_\delta - \partial^{\lfloor \alpha\rfloor}T_{\delta'}\Vert_{S_p}^p \\ 
   \begin{aligned}
   & \leq  CC_1^p \sum\limits_{n=n_0}^{+\infty} \frac{1}{(n+1)^{1+p}}+CC_2^{p}\vert \delta-\delta'\vert^p \sum\limits_{n=0}^{n_0-1} \frac{1}{n+1}  \\
     & \leq CC_1^p \int_{n_0}^\infty \frac{1}{x^{1+p}}dx + CC_2^p \vert \delta-\delta'\vert^p + CC_2^p \vert \delta-\delta'\vert^p \int_1^{n_0-1} \frac{1}{x}dx\\
     & \leq  \frac{CC_1^p}{p n_0^p}+  CC_2^p \vert \delta-\delta'\vert^p + CC_2^p \vert \delta-\delta'\vert^p \ln (n_0-1)\\
     & \leq  \frac{CC_1^p}{p}\vert \delta-\delta'\vert^p + CC_2^p \vert \delta-\delta'\vert^p + CC_2^p \vert \delta-\delta'\vert^p \ln \left(\frac{1}{\vert \delta-\delta'\vert}\right)\\
     & \leq D^p \vert \delta-\delta'\vert^p \vert \ln \vert \delta-\delta'\vert \vert
\end{aligned}\end{multline*}where $D>$ does not depend on $\delta,\delta'\in I$.
But $x^\varepsilon \vert\ln x\vert^{1/p}\underset{x\to 0}{\to}0$ so there is a constant $C_\varepsilon>0$ such that $\vert \ln \vert \delta-\delta'\vert \vert^{1/p} \leq C_\varepsilon \vert \delta-\delta'\vert^{-\varepsilon}$ and thus, \begin{equation*}\Vert \partial^{\lfloor \alpha\rfloor}T_\delta - \partial^{\lfloor \alpha\rfloor}T_{\delta'}\Vert_{S_p}\leq DC_{\varepsilon}\vert \delta-\delta'\vert^{1-\varepsilon}=DC_{\varepsilon}\vert \delta-\delta'\vert^{\alpha-\lfloor\alpha\rfloor}.\qedhere\end{equation*}
\end{proof}

For $2+\frac{2}{\dim G/K-1}<p\leq +\infty$, we now denote $(r_p,\delta_p)$ the regularity of $\delta\mapsto T_\delta$ as a map with values in $S_p(L^2(G))$ obtained in Corollary \ref{schatteninfty} and Theorem \ref{schatten}. We can use these two results to study the regularity of $S_p$-multipliers and give a generalisation of \cite[Prop. 4.2]{dlSPR}. Let $\varphi:G\to \C$ be a bounded measurable map. We can consider the map $$\fonction{S_\varphi}{S_2(L^2(G))}{S_2(L^2(G))}{(a_{g,h})_{g,h\in G}}{(\varphi(gh^{-1}))_{g,h\in G}}.$$If $S_\varphi:S_p(L^2(G))\cap S_2(L^2(G))\mapsto S_p(L^2(G))$ is bounded, by density we can extend it to $S_p(L^2(G))$ and if it remains bounded, we say that $\varphi$ is an $S_p$-multiplier.

\begin{coro}\label{spmult}Let $1<p\leq +\infty$. Let $\varphi$ be a $K$-bi-invariant $S_p$-multiplier of $G$ and $\psi=\varphi\circ \exp$. Then $\psi\in C^{(r_p,\delta_p)}(Q)$.
\end{coro}

\begin{proof}If $\delta\in Q$, we have $S_\varphi(T_\delta)=\psi(\delta)T_\delta$ by \cite[Remark 4.5]{dlSPR}. Thus, if $\mathbf{1}$ denote the constant function, which belongs to $L^2(G)$, we have $$\psi(\delta)=\langle S_\varphi(T_\delta)\mathbf{1},\mathbf{1}\rangle.$$So $\psi$ is at least as regular as $T$ is.\end{proof}

\subsection{Optimality of the results}\label{optimalitysection}
In this section, we will show that the results obtained before are optimal. We keep the notations of Theorem \ref{schatten}, where we view the functions on $]-1,1[$ instead of $Q$ (using Lemma \ref{precomposition}).

\begin{theorem}
 \label{optimality}Let $1<p<+\infty$. For $\delta\in ]-1,1[$, the operator $\partial^r T_\delta$ does not belong to $S_p(L^2(G))$ as soon as $r\geq \alpha_p$. Furthermore, if $$r=\left\lbrace\begin{aligned}\lfloor \alpha_p\rfloor & \textrm{ if }\alpha_p\not\in \Z\\\lfloor \alpha_p\rfloor-1 & \textrm{ if }\alpha_p\in \Z\end{aligned}\right.$$ then for each compact interval $I$ of $]-1,1[$, there exists a constant $C_{I,p}$ such that for all $x,y\in I$, $$\Vert \partial^r T_x-\partial^r T_y\Vert_{S_p} \geq C_{I,p}\vert x-y\vert^{\alpha_p-r}.$$   
\end{theorem}

\begin{lem}
\label{equirepcos}Let $\alpha,\beta\geq 0$ two fixed reals. For any $\theta\in ]0,\pi[$, there exists $C>0$ and an integer $N_0$ such that for $N\geq N_0$, we have $$\left(\sum_{n=N+1}^{2N} n\vert P_n^{(\alpha,\beta)}(\cos \theta)\vert^2\right)^{1/2}\geq C\sqrt{N}.$$
\end{lem}
\begin{proof}Denote $$S=\left(\sum_{n=N+1}^{2N} n\vert P_n^{(\alpha,\beta)}(\cos \theta)\vert^2\right)^{1/2}.$$By Darboux's formula (Proposition \ref{darboux}), there is a constant $M$ such that for all $n\in \N^*$, $$\vert P_n^{(\alpha,\beta)}(\cos \theta)-n^{-\frac{1}{2}}k(\theta)\cos(N\theta+\gamma)\vert\leq M n^{-\frac{3}{2}}.$$By triangle inequality, we thus have $$\vert k(\theta)\cos(N\theta+\gamma)\vert \leq \sqrt{n}\vert P_n^{(\alpha,\beta)}(\cos \theta)\vert + M n^{-1} $$so $$\frac{1}{2}\vert k(\theta)\cos(N\theta+\gamma)\vert^2 \leq n\vert P_n^{(\alpha,\beta)}(\cos \theta)\vert^2 + M^2 n^{-2}.$$Thus, we get \begin{align*}
    \sum_{n=N+1}^{2N} \frac{1}{2}\vert k(\theta)\cos(N\theta+\gamma)\vert^2 & \leq  S^2 + M^2  \sum_{n=N+1}^{2N} n^{-2} \\
     & \leq  S^2+ M^2  \sum_{n=N+1}^{+\infty} n^{-2}\\
     & \leq  S^2 + M^{2} N^{-1}.
\end{align*}Taking square root, since $\sqrt{a+b}\leq \sqrt{a}+\sqrt{b}$, we get $$S \geq \left(\sum_{n=N+1}^{2N} \frac{1}{2}\vert k(\theta)\cos(N\theta+\gamma)\vert^2\right)^{1/2}-MN^{-1/2}.$$
Denote $u=(N+1+\frac{1}{2}(\alpha+\beta+1))\theta+\gamma$. Then \begin{align*}\sum_{n=N+1}^{2N} \frac{1}{2}\vert k(\theta)\cos(N\theta+\gamma)\vert^2 &= \frac{k(\theta)^2}{2}\sum_{n=0}^{N-1} \vert \cos(u+n\theta)\vert^2\\&=\frac{Nk(\theta)^2}{4}+\frac{k(\theta)^2}{2}\sum_{n=0}^{N-1}\cos(2u+2n\theta).\end{align*}
But we have that $\vert\sum_{n=0}^{N-1}\cos(x+ny)\vert = \vert \mathrm{Re}\sum e^{i(x+ny)}\vert \leq \frac{1}{\vert \sin y/2\vert}$, thus $$\sum_{n=N+1}^{2N} \frac{1}{2}\vert k(\theta)\cos(N\theta+\gamma)\vert^2 \geq \frac{Nk(\theta)^2}{4}-\frac{k(\theta)^2}{2\vert\sin\theta\vert}\geq C'N$$for some $C'>0$ and $N$ large enough.

Finally, $S\geq \sqrt{C'N}-M N^{-1/2}\geq C\sqrt{N}$ for $C=\frac{C'}{2}$ and $N$ large enough.
\end{proof}

\begin{lem}
\label{lemmeutile}Let $\alpha,\beta\geq 0$ two fixed reals. Let $J$ be a compact interval of $]0,\pi[$. There exists a real $C>0$ and an integer $m_0\geq 1$ such that for all $\theta,\phi\in J$, $m\geq m_0$ and $N\geq \frac{m_0}{\vert \theta-\phi\vert}$, we have $$\left(\sum_{n=m+1}^{m+N} n\vert P_n^{(\alpha,\beta)}(\cos \theta)-P_n^{(\alpha,\beta)}(\cos \phi)\vert^2\right)^{1/2}\geq C\sqrt{N}.$$
\end{lem}

\begin{proof}
Denote $$S=\left(\sum_{n=m+1}^{m+N} n\vert P_n^{(\alpha,\beta)}(\cos \theta)-P_n^{(\alpha,\beta)}(\cos \phi)\vert^2\right)^{1/2}.$$By Darboux's formula (Proposition \ref{darboux}), there is a constant $M_J$ such that for all $\theta\in J$ and $n\in \N^*$, $$\vert P_n^{(\alpha,\beta)}(\cos \theta)-n^{-\frac{1}{2}}k(\theta)\cos(N\theta+\gamma)\vert\leq M_Jn^{-\frac{3}{2}}.$$By triangle inequality, we thus have $$\vert k(\theta)\cos(N\theta+\gamma)-k(\phi)\cos(N\phi+\gamma)\vert \leq \sqrt{n}\vert P_n^{(\alpha,\beta)}(\cos \theta)-P_n^{(\alpha,\beta)}(\cos \phi)\vert + 2M_J n^{-1} $$so $$\frac{1}{2}\vert k(\theta)\cos(N\theta+\gamma)-k(\phi)\cos(N\phi+\gamma)\vert^2 \leq n\vert P_n^{(\alpha,\beta)}(\cos \theta)-P_n^{(\alpha,\beta)}(\cos \phi)\vert^2 + 4M_J^2 n^{-2}.$$Thus, we get \begin{align*}
    \sum_{n=m+1}^{m+N} \frac{1}{2}\vert k(\theta)\cos(N\theta+\gamma)-k(\phi)\cos(N\phi+\gamma)\vert^2 & \leq  S^2 + 4M_J^2  \sum_{n=m+1}^{m+N} n^{-2} \\
     & \leq  S^2+ 4M_J^2  \sum_{n=m+1}^{+\infty} n^{-2}\\
     & \leq  S^2 + 4M_J^{2} m^{-1}.
\end{align*}
    
Taking square root, since $\sqrt{a+b}\leq \sqrt{a}+\sqrt{b}$, we get $$S \geq \left(\sum_{n=m+1}^{m+N} \frac{1}{2}\vert k(\theta)\cos(N\theta+\gamma)-k(\phi)\cos(N\phi+\gamma)\vert^2\right)^{1/2}-2M_Jm^{-1/2}.$$
Denote $u=(m+\frac{1}{2}(\alpha+\beta+1)+1)\theta+\gamma$, $v=(m+\frac{1}{2}(\alpha+\beta+1)+1)\phi+\gamma$, $a=\frac{k(\theta)}{\sqrt{2}}$ and $b=\frac{k(\phi)}{\sqrt{2}}$. Then $$\sum_{n=m+1}^{m+N} \frac{1}{2}\vert k(\theta)\cos(N\theta+\gamma)-k(\phi)\cos(N\phi+\gamma)\vert^2 = \sum_{n=0}^{N-1} \vert a\cos(u+n\theta)-b\cos(v+n\phi)\vert^2.$$Making use of the formula $2\cos(x)\cos(y)=\cos(x+y)+\cos(x-y)$, we have \begin{multline*}
    \vert a\cos(u+n\theta)-b\cos(v+n\phi)\vert^2  \\ \begin{aligned} 
    & =a^2\cos^2(u+n\theta)+b^2\cos^2(v+n\phi)-2ab\cos(u+n\theta)\cos(v+n\phi) \\
     & =  \frac{1}{2}[a^2+b^2-2ab\cos(u-v+n(\theta-\phi))+a^2\cos(2u+2n\theta)\\
     &\phantom{=}+b^2\cos(2v+2n\phi)-2ab\cos(u+v+n(\theta+\phi))].
\end{aligned}\end{multline*}

Again using $\vert\sum_{n=0}^{N-1}\cos(x+ny)\vert = \vert \mathrm{Re}\sum e^{i(x+ny)}\vert \leq \frac{1}{\vert \sin y/2\vert}$, we get \begin{multline*}
    \frac{2}{N}\sum_{n=0}^{N-1} \vert a\cos(u+n\theta)-b\cos(v+n\phi)\vert^2 \geq\\ a^2+b^2-\frac{2ab}{N\vert \sin \frac{\theta-\phi}{2}\vert}-\frac{a^2}{N\vert \sin \theta\vert}-\frac{b^2}{N\vert \sin \phi\vert}-\frac{2ab}{N\vert \sin \frac{\theta+\phi}{2}\vert}.
\end{multline*}Now if $N\geq \frac{m_0}{\vert \theta-\phi\vert}$ for some $m_0$ depending on $J$, we get that $$\frac{2}{N}\sum_{n=0}^{N-1} \vert a\cos(u+n\theta)-b\cos(v+n\phi)\vert^2\geq \frac{a^2+b^2}{2}\geq C'$$where $C'>0$ is a constant depending on $\alpha,\beta,J$.

Finally, $S\geq \sqrt{\frac{C'N}{2}}-2M_J m^{-1/2}\geq C\sqrt{N}$ for $m$ large enough.
\end{proof}

\begin{proof}[Proof of Theorem \ref{optimality}]Let $(\varphi_n)_{n\in\N}$ be the family of spherical functions of the pair $(G,K)$ and $m_n$ the dimension of the associated representation of $G$. Then for $r\in \N$, $x\in ]-1,1[$, we have $$\Vert \partial^r T_x\Vert^p_{S_p}=\sum_{n\geq 0}m_n\vert \varphi_n^{(r)}(x)\vert^p.$$
As in (\ref{dimension}), there is $C_0>0$ depending only on the pair $(G,K)$ such that $$m_n \geq C_0 n^{\dim G/K -1}$$for all $n\in \N$. Furthermore, recall that there is $\theta\in ]0,\pi[$ such that $x=\cos \theta$ and $$\vert\varphi_n^{(r)}(x)\vert=\frac{\Gamma(\alpha+1)}{2^r}\frac{\Gamma(n+1)}{\Gamma(n+\alpha+1)}\frac{\Gamma(\alpha+\beta+n+1+r)}{\Gamma(\alpha+\beta+n+1)}\vert P_{n-r}^{(\alpha+r,\beta+r)}(\cos \theta)\vert.$$
But $$\frac{\Gamma(\alpha+\beta+n+1+r)}{\Gamma(\alpha+\beta+n+1)}=\prod_{j=1}^r (\alpha+\beta+n+j)\geq n^r$$
and as in (\ref{gamma}), if we define $f(x)=\sqrt{2\pi}x^{x+\frac{1}{2}}e^{-x}$, there are constants $C_1,C_2>0$ such that for $x\geq 0$, $$C_1f(x)\leq \Gamma(x+1)\leq C_2f(x).$$Thus for all $n\in \N$, \begin{align*}
    \frac{\Gamma(n+1)}{\Gamma(n+\alpha+1)} & \geq  \frac{C_1}{C_2}\frac{n^{n+\frac{1}{2}}}{(n+\alpha)^{n+\alpha+\frac{1}{2}}}\frac{e^{-n}}{e^{-n-\alpha}}\\
     & \geq  \frac{C_1e^\alpha}{C_2}n^{-\alpha}\sqrt{\frac{n}{n+\alpha}}\left(1-\frac{\alpha}{n+\alpha}\right)^{n+\alpha}\\
     & \geq  C_3 n^{-\alpha} 
\end{align*} where $C_3>0$ is a constant depending only on $\alpha$ thus on the pair $(G,K)$.

Combining these estimates, reindexing the sum and using $n-k\geq C_4 n$ for some $C_4>0$ and for all $n>k$, we get that there is a constant $C>0$ which depends only on $(G,K)$ such that \begin{equation}
    \label{minorationnormepbase} \Vert \partial^r T_x\Vert^p_{S_p}\geq C \sum_{n\geq 1} n^{p(r-\alpha)+\dim G/K -1} \vert P_{n}^{(\alpha+r,\beta+r)}(\cos \theta)\vert^p
.\end{equation}Recall (Section \ref{sectionspheriques}) that $\alpha=\alpha(G,K)=\frac{\dim G/K}{2}-1$ so that \begin{align*}p(r-\alpha)+\dim G/K -1 &= p\left(r-\left(\frac{\dim G/K-1}{2}-\frac{\dim G/K}{p}\right)+\frac{1}{2}\right)-1\\&=p(r-\alpha_p+\frac{1}{2})-1.\end{align*}
To simplify notations, denote $\kappa=p(r-\alpha_p+\frac{1}{2})-1$. Note that for now, we did not make assumptions on $r$ so we do not know if the right-hand side converges.

For any $N\in \N$, we have by Hölder's inequality that \begin{equation}
    \label{holderbase}    \sum_{n=N+1}^{2N} n^\kappa \vert P_{n}^{(\alpha+r,\beta+r)}(\cos \theta)\vert^p \geq \frac{\left(\sum_{n=N}^{2N} n \vert P_{n}^{(\alpha+r,\beta+r)}(\cos \theta)\vert^2\right)^{p/2}}{\left(\sum_{n=N+1}^{2N} n^{\frac{p-2\kappa}{p-2}}\right)^{p/2-1}}.
\end{equation}
By Lemma \ref{equirepcos}, for $N$ large enough, the numerator is greater than $C^pN^{p/2}$. For the denominator, we have $$\left(\sum_{n=N+1}^{2N} n^{\frac{p-2\kappa}{p-2}}\right)^{p/2-1}\leq \left(N (2N)^{\frac{p-2\kappa}{p-2}}\right)^{\frac{p-2}{2}}\leq C' N^{p-\kappa-1}$$ thus there is $\Tilde{C}>0$ such that $$\sum_{n=N+1}^{2N} n^\kappa \vert P_{n}^{(\alpha+r,\beta+r)}(\cos \theta)\vert^p \geq \Tilde{C}N^{\kappa+1-\frac{p}{2}}.$$ Then if $r\geq \alpha_p$, we get $\kappa+1-\frac{p}{2}\geq 0$, thus by (\ref{holderbase}), $$\sum_{n=N+1}^{2N} n^\kappa \vert P_{n}^{(\alpha+r,\beta+r)}(\cos \theta)\vert^p \underset{N\to \infty}{\not\to} 0$$so by (\ref{minorationnormepbase}), the $p$-norm of $\partial^r T_x$ is not finite so $\partial^r T_x\not\in S_p(L^2(G))$.

Now set $r<\alpha_p$ as in Theorem \ref{schatten}. Using what is above, for $x,y\in I$, we have $$\Vert \partial^r T_x-\partial^r T_y\Vert^p_{S_p}=\sum_{n\geq 0}m_n\vert \varphi_n^{(r)}(x)-\varphi_n^{(r)}(y)\vert^p.$$ Let $J=\arccos(I)$ which is a compact interval of $]0,\pi[$. There are $\theta,\phi\in J$ such that $x=\cos \theta$ and $y=\cos \phi$. Using the same inequalities on $m_n$ and $\varphi_n$ as before, we get that there is a constant  $C>0$ which depends only on $(G,K)$ such that \begin{equation}
    \label{minorationnormep} \Vert \partial^r T_x-\partial^r T_y\Vert^p_{S_p}\geq C \sum_{n\geq 1} n^{\kappa} \vert P_{n}^{(\alpha+r,\beta+r)}(\cos \theta)- P_{n}^{(\alpha+r,\beta+r)}(\cos\phi)\vert^p
.\end{equation}
Set $m_0$ as in Lemma \ref{lemmeutile} and $m_k=m_0+kN$. Then we get \begin{equation}
    \label{normep2} \Vert \partial^r T_x-\partial^r T_y\Vert^p_{S_p}\geq C \sum_{k\geq 0}\sum_{n=m_k+1}^{m_{k+1}} n^\kappa \vert P_{n}^{(\alpha+r,\beta+r)}(\cos \theta)- P_{n}^{(\alpha+r,\beta+r)}(\cos\phi)\vert^p.
\end{equation}
By Hölder's inequality, we have that \begin{multline}
    \label{holder}    \sum_{n=m_k+1}^{m_{k+1}} n^\kappa \vert P_{n}^{(\alpha+r,\beta+r)}(\cos \theta)- P_{n}^{(\alpha+r,\beta+r)}(\cos\phi)\vert^p \geq\\ \frac{\left(\sum_{n=m_k+1}^{m_{k+1}} n \vert P_{n}^{(\alpha+r,\beta+r)}(\cos \theta)- P_{n}^{(\alpha+r,\beta+r)}(\cos\phi)\vert^2\right)^{p/2}}{\left(\sum_{n=m_k+1}^{m_{k+1}} n^{\frac{p-2\kappa}{p-2}}\right)^{p/2-1}}.
\end{multline}

By Lemma \ref{lemmeutile}, for $N\geq \frac{m_0}{\vert \theta-\phi\vert}$, the numerator is greater than $C'^pN^{p/2}$. For the denominator, first notice that $$m_{k+1}=m_0+(k+1)N=N\left(\frac{m_0}{N}+k+1\right)\leq (\pi+k+1)N$$ so that \begin{align*}
    \left(\sum_{n=m_k+1}^{m_{k+1}} n^{\frac{p-2\kappa}{p-2}}\right)^{p/2-1} & \leq  \left(N m_{k+1}^{\frac{p-2\kappa}{p-2}} \right)^{\frac{p-2}{2}}\\
     & \leq  \left( (\pi+k+1)^{\frac{p-2\kappa}{p-2}} N^{\frac{p-2\kappa}{p-2}+1}\right)^{\frac{p-2}{2}}\\
     & \leq  (\pi+k+1)^{\frac{p-2\kappa}{2}} N^{p-\kappa-1}.
\end{align*}

Using theses two inequalities in (\ref{holder}) and then in (\ref{normep2}) we finally get \begin{equation}
    \label{normep3} \Vert \partial^r T_x-\partial^r T_y\Vert^p_{S_p}\geq CC'^p N^{\kappa+1-\frac{p}{2}} \sum_{k\geq 0} (1+\pi+k)^{\kappa-\frac{p}{2}}.
\end{equation}

Notice that this sum converges because $\kappa - \frac{p}{2}  <  -1$. We get $$\sum_{k\geq 0} (1+\pi+k)^{\kappa-\frac{p}{2}} \geq \int_0^\infty (1+\pi+x)^{\kappa-\frac{p}{2}}= \frac{\left(1+\pi\right)^{1+\kappa-\frac{p}{2}}}{\frac{p}{2}-\kappa-1}.$$Now since we want $N\geq \frac{m_0}{\vert \theta-\phi\vert}$, we can choose $N$ so that $$N\leq \frac{m_0}{\vert \theta-\phi\vert}+1 \leq \frac{m_0+\pi}{\vert \theta-\phi\vert}.$$
Thus taking the $p$-th root in (\ref{normep3}) and using that $\frac{\kappa+1}{p}-\frac{1}{2}=r-\alpha_p$, we have \begin{equation}
    \label{normep4} \begin{split}
        \Vert \partial^r T_x-\partial^r T_y\Vert_{S_p}&\geq  C^{1/p}C' \left(\frac{m_0+\pi}{\vert \theta-\phi\vert}\right)^{\frac{\kappa+1}{p}-\frac{1}{2}} \frac{\left(1+\pi\right)^{\frac{1+\kappa}{p}-\frac{1}{2}}}{\left(\frac{p}{2}-\kappa-1\right)^{1/p}} \\
         & \geq  C^{1/p}C' \frac{\left((m_0+\pi)(1+\pi)\right)^{r-\alpha_p}}{\left(\frac{p}{2}-\kappa-1\right)^{1/p}} \vert \theta-\phi\vert^{\alpha_p-r}.
    \end{split}
\end{equation}and finally, we conclude using that $\vert\theta-\phi\vert \geq \Tilde{C}\vert \cos\theta-\cos\phi\vert$ for some $\Tilde{C}>0$.
\end{proof}

\begin{coro}
\label{coropti}We have the following results:\begin{enumerate}
    \item If $p\leq 2+\frac{2}{\dim G/K-1}$ and $\delta\in Q$, $T_\delta\not\in S_p(L^2G)$.
    \item If $2+\frac{2}{\dim G/K-1}<p<+\infty$, the regularity obtained in Theorem \ref{schatten} is optimal.
    \item If $p=+\infty$, the regularity obtained in Proposition \ref{schatteninfty} is optimal.
    \item Theorem \ref{mainthm} is optimal: for any $(r,\alpha)>(\lfloor \alpha_\infty\rfloor,\alpha_\infty - \lfloor \alpha_\infty \rfloor)$ in lexicographic order, the family of spherical functions is not bounded in $C^{(r,\alpha)}(Q)$.
    \item Corollary \ref{matrixcoef} is optimal: with the same notations as above, there are $K$-bi-invariant matrix coefficients of unitary representations that are not in $C^{(r,\alpha)}(Q)$.
\end{enumerate}
\end{coro}

\begin{remark}
    Even though the regularity of $T$ is optimal, we do not know whether Corollary \ref{spmult} is optimal or not. Indeed, we only say that $S_p$-multipliers are at least as regular as $T$, but we cannot construct specific multipliers with the exact regularity of $T$.
\end{remark}
\begin{proof}
\begin{enumerate}
    \item If $p\leq 2+\frac{2}{\dim G/K-1}$, then $\alpha_p\leq 0$ so this is the first part of Theorem \ref{optimality}.
    \item There are two cases to consider.
In the first case, if $\alpha_p\in \Z$, we know that $T\in C^{(\alpha_p-1,1-\varepsilon)}(]-1,1[,S_p(L^2(G)))$ for any $\varepsilon\in ]0,1[$. By the first part of Theorem \ref{optimality}, $T\not\in C^{(\alpha_p,0)}(]-1,1[,S_p(L^2(G)))$ (because $\partial^{\alpha_p}T_\delta\not\in S_p$). Then, assume $\partial^{\alpha_p-1}T$ is Lipschitz on a compact interval $J$ of $]-1,1[$. Since $1<p<+\infty$, $S_p(L^2(G))$ is reflexive, thus by \cite[Corollary 5.12]{benyamini1998geometric}, $\partial^{\alpha_p-1}T:\mathrm{Int}(J)\to S_p(L^2(G))$ is differentiable almost everywhere (for the Lebesgue measure on $J$) which contradicts the fact that $\partial^{\alpha_p}T_\delta\not\in S_p(L^2(G))$ for any $\delta\in ]-1,1[$. Thus, $T\not\in C^{(\alpha_p-1,1)}(]-1,1[,S_p(L^2(G)))$.

In the second case, we consider $\alpha_p\not\in \Z$. Let $r=\lfloor\alpha_p\rfloor$. Then we know that $T\in C^{(r,\alpha_p-r)}(]-1,1[,S_p(L^2(G)))$. Let $1-(\alpha_p-r)>\varepsilon>0$, and $I$ compact interval of $]-1,1[$. Assume that there is $C_I>0$ such that for all $x,y\in I$, $$\Vert \partial^r T_x-\partial^r T_y\Vert \leq C_I\vert x-y\vert^{\alpha_p-r+\varepsilon}.$$
Then for $x\neq y\in I$, by the second part of Theorem \ref{optimality}, we get $$C_{I,p}\leq C_I\vert x-y\vert^{\varepsilon},$$which is impossible when $x\to y$.

Thus, for any $\varepsilon>0$, $T\not\in C^{(r,\alpha_p-r+\varepsilon)}(]-1,1[,S_p(L^2(G)))$.

\item Once again we distinguish two cases. First, if $\alpha_\infty\not\in \Z$, let $1<p<+\infty$ be large enough so that $r=\lfloor \alpha_p\rfloor = \lfloor \alpha_\infty\rfloor$. Then, by the second part of Theorem \ref{optimality}, we get for $x,y\in I$ compact subset of $]-1,1[$ that $$\Vert \partial^r T_x-\partial^r T_y\Vert_{S_p} \leq C_{I,p}\vert x-y\vert^{\alpha_p-r}.$$From the expression of $C_{I,p}$, we see that $\underset{p\to +\infty}{\lim} C_{I,p}=C_I$ exists and is finite. Then we get with $p\to +\infty$ that $$\Vert \partial^r T_x-\partial^r T_y\Vert_{S_\infty}\leq C_I\vert x-y\vert^{\alpha_\infty-r}.$$Thus, as above, $T\not\in C^{(r,\alpha_\infty-r+\varepsilon)}(]-1,1[,S_\infty(L^2(G)))$ for any $\varepsilon>0$.

Next, if $\alpha_\infty\in \Z$, we have $T\in C^{(\alpha_\infty-1,1)}(]-1,1[,S_\infty(L^2(G)))$. Let $x=\cos \theta$, $y=\cos \phi$, $\theta,\phi\in [\varepsilon,\pi-\varepsilon]$. We have \begin{align*}
\Vert \partial^{\alpha_\infty}T_x-\partial^{\alpha_\infty}T_y\Vert_{S_\infty} &=  \underset{n}{\sup}\vert \varphi_n^{(\alpha_\infty)}(\cos\theta)-\varphi_n^{(\alpha_\infty)}(\cos \phi)\vert\\
& \geq  \left(\frac{1}{N}\sum_{n=m}^{m+N}\vert \varphi_n^{(\alpha_\infty)}(\cos\theta)-\varphi_n^{(\alpha_\infty)}(\cos \phi)\vert^2\right)^{1/2}\\
& \geq  C(G,K,\varepsilon)
\end{align*}by Lemma \ref{lemmeutile}. Thus, we get that $\partial^{\alpha_\infty}T$ is not continuous, so we showed that $T\not\in C^{(\alpha_\infty,0)}(]-1,1[,S_\infty(L^2(G)))$.

\item It follows from $(3)$ and Lemma \ref{existencederiv}.
\item It is a consequence of $(4)$ and Lemma \ref{roptid}.\qedhere
\end{enumerate}
\end{proof}

\section{Higher rank symmetric spaces}\label{higherranksection}

\subsection{The case of a Lie group seen as a symmetric space}\label{caracteres}
We consider $G$ a compact semisimple Lie group and $\Delta(G)$ the diagonal subgroup of $G\times G$. We study the symmetric Gelfand pair $(G\times G,\Delta(G))$. Note that unlike what we did in Section \ref{sectionpreli}, we do not need to assume that the symmetric space is simply connected and we do not even need it to be irreducible (which means $G$ does not need to be simply connected nor simple).

In this case, spherical functions can be described in a more efficient way than what the Cartan-Helgason theorem (Theorem \ref{cartanhelgason}) tells us. An irreducible representation of $G\times G$ is of the form $(\pi\otimes \sigma,V\otimes W)$ where $\pi,\sigma$ are irreducible representations of $G$. We have that $V\otimes W\simeq \Hom(V^*,W)$, and the representation is given by $$\left((\pi\otimes \sigma)(g,g')f\right)(v)=\sigma(g')\left(f(\pi(g)^*v)\right).$$ Assume that there is a $\Delta(G)$-invariant vector $\xi\in \Hom(V^*,W)$. Then for all $g\in G$, we have $$\sigma(g)(\xi(\pi(g)^*v))=\xi(v)$$so $\xi$ is $G$-equivariant from $V^*\to W$, so by Schur's lemma, $V^*\simeq W$. Thus, the classes of irreducible representations of $G\times G$ with a $\Delta(G)$-invariant vector are in bijection with the classes of irreducible representations of $G$. On $\Hom(V,V)$, the scalar product is $\langle f,g\rangle=\Tr(fg^*)$. For $(\pi,V)$ a (unitary) irreducible representation of $G$, the associated spherical function of $(G\times G,\Delta(G))$ is $$\varphi_\pi:(g,g')\mapsto \frac{\langle(\pi\otimes \pi^*)(g,g')(\Id),\Id\rangle}{\langle\Id,\Id\rangle}=\frac{\langle\pi(g),\pi(g')\rangle}{\dim V}=\frac{\chi_\pi(gg'^{-1})}{\dim V}.$$So the spherical functions of the pair $(G\times G,\Delta(G))$ are just the normalized characters of $G$.

Let $\mathfrak{g}$ be the Lie algebra of $G$. The Lie algebra of $G\times G$ is $\mathfrak{g}\oplus\mathfrak{g}$, and the Lie algebra of $\Delta(G)$ is the subspace $$\mathfrak{k}=\{(H,H) \vert H\in \mathfrak{g}\}.$$Then, we have $\mathfrak{g}\oplus\mathfrak{g}=\mathfrak{k}\oplus\mathfrak{m}$ where $$\mathfrak{m}=\{(H,-H) \vert H\in \mathfrak{g}\}.$$

Let $\mathfrak{b}$ be a maximal abelian subalgebra of $\mathfrak{g}$, $\Phi$ the root system associated to $(\mathfrak{g}_\C,\mathfrak{b}_\C)$, $\Phi^+$ a choice of positive roots and $\Delta=\{\alpha_1,\cdots,\alpha_\ell\}$ be the associated basis, where $\ell=\rank G$. We have $$\mathfrak{g}_\C=\mathfrak{b}_\C\oplus\bigoplus_{\alpha\in \Phi}\mathfrak{g}_\alpha$$with $\mathfrak{g}_\alpha=\{X\in \mathfrak{g}_\C \vert \forall H\in \mathfrak{b}_\C, [H,X]=\alpha(H)X\}$ and $\dim g_\alpha=1$ for all $\alpha\in \Phi$.

Let $\Lambda_G$ be the set of dominant analytically integral element. We know that there is a bijection between $\Lambda_G$ and the set of equivalence classes of finite dimensional irreducible representations of $G$ (\cite[Thm. 5.110]{knapp2002lie}). For $\lambda\in \Lambda_G$, let $\pi_\lambda$ be an associated representation, $\chi_\lambda$ its character and $d_\lambda$ its dimension. Then, the previous result amount to the following :
\begin{prop}
For $\lambda\in \Lambda_G$, define $\varphi_\lambda:(g,g')\mapsto \frac{\chi_\lambda(gg'^{-1})}{d_\lambda}$. Then $\lambda\mapsto \varphi_\lambda$ is a bijection from $\Lambda_G$ to the set of spherical functions of $((G\times G),\Delta(G))$.
\end{prop}

Now consider $\mathfrak{a}=\{(H,-H)\vert H\in \mathfrak{b}\}$. We have that $\mathfrak{a}$ is a maximal abelian subspace of $\mathfrak{m}$. Thus, we can consider $\Sigma_{\mathfrak{a}}$ root system of $\mathfrak{a}_\C$ in $(\mathfrak{g}\times \mathfrak{g})_\C$, whose roots are $\Tilde{\alpha}:(H,-H)\mapsto \alpha(H)$ for $\alpha\in \Phi$. We choose as a positive root system $\Sigma_\mathfrak{a}^+=\{\Tilde{\alpha}\vert \alpha\in \Phi^+\}$. Let $\Tilde{\mathfrak{g}}_\alpha=\{X\in (\mathfrak{g}\times \mathfrak{g})_\C \vert \forall H\in \mathfrak{a}_\C, [H,X]=\Tilde{\alpha}(H)X\}$. Then we have $$(\mathfrak{g}\times \mathfrak{g})_\C=\mathfrak{a}_\C\oplus \left(\mathfrak{k}^{\mathfrak{a}}\right)_\C\oplus \bigoplus_{\alpha\in \Phi} \Tilde{\mathfrak{g}}_\alpha$$and we see that $\mathfrak{k}^{\mathfrak{a}}=\{(H,H) \vert H\in \mathfrak{b}\}$ and for all $\alpha\in \Phi$, $\dim \Tilde{\mathfrak{g}}_\alpha=2$ - in fact, $\Tilde{\mathfrak{g}}_\alpha=(\mathfrak{g}_\alpha \times 0)\oplus (0\times \mathfrak{g_{-\alpha}})$.

Now, recall that we defined $Q$ (Proposition \ref{kakcpt}) as the connected component of $\{H\in \mathfrak{a} \vert  \forall \alpha\in \Sigma_\mathfrak{a}, \alpha(H)\vert\in i\pi\Z\}$ contained in the Weyl chamber associated to $\Sigma_\mathfrak{a}^+$ and containing $0$ in its closure. We want to study the regularity of the functions $\psi_\lambda=\varphi_\lambda\circ \exp\vert_Q$.

By definition, any positive root $\alpha\in \Phi^+$ is a linear combination of the roots in $\Delta$, with coefficients in $\N$. We write $\alpha=\sum_{i=1}^\ell n_i(\alpha)\alpha_i$. We define $$\gamma=\underset{1\leq i\leq \ell}{\min} \vert \{\alpha \in \Phi^+ \vert n_i(\alpha)\geq 1\}\vert.$$
\begin{remark}The number $\gamma$ was used by Cowling and Nevo in \cite{cowling2001uniform} in estimates of spherical functions of the non-compact dual of $(G\times G)/\Delta(G)$, namely $G_\C/G$. A table of the values of $\gamma$ can be found in \cite[Appendix]{10.2307/2001308}.
\end{remark}

As in Section \ref{sectionrank1}, we begin by showing that $\Delta(G)$-bi-invariant matrix coefficients of unitary representations of $G\times G$ are at least of regularity $\gamma$.
\begin{theorem}
 \label{regcar}The family $(\psi_\lambda)_{\lambda\in \Lambda_G}$ is bounded in $C^{(\gamma,0)}(Q)$.   
\end{theorem}
\begin{proof}
Let $\Tilde{H}=(H,-H)\in Q$, $H\in \mathfrak{b}$. Then $$\psi_\lambda(\Tilde{H})=\varphi_\lambda(\exp(H),\exp(-H))=\frac{\chi_\lambda(\exp(2H))}{d_\lambda}.$$
Let $q$ be the Weyl denominator, defined on $\mathfrak{b}$ by $q(H)=\prod_{\alpha\in \Phi^+} (e^{\alpha(H)/2}-e^{-\alpha(H)/2})$. Note that the roots $\alpha$ are in $i\mathfrak{b}^*$, so $\alpha(H)\in i\R$. Thus, if $\alpha(H)/2 \not\in i\pi \Z$, we get that $e^{\alpha(H)/2}-e^{-\alpha(H)/2}=2\sinh(\alpha(H)/2)$ is non-zero. So for any $\Tilde{H}=(H,-H)\in Q$, by definition of $Q$, $2H$ is such that $q(2H)\neq 0$.

Now, by the Weyl character formula (\cite[Thm. 10.14]{hall2003lie}), if $q(2H)\neq 0$, we have $$\chi_\lambda(e^{2H})=\frac{\sum_{w\in W}\det(w)e^{(w(\lambda+\rho))(2H)}}{q(2H)}$$where $W$ denote the Weyl group and $\rho=\frac{1}{2}\sum_{\alpha\in \Phi^+}\alpha$.

Since the map $(H,-H)\mapsto q(2H)$ is smooth and non-zero on $Q$, by Lemma \ref{multsmooth}, we can study the functions $\Tilde{\psi}_\lambda:(H,-H)\mapsto q(2H)\psi_\lambda(H,-H)$, which will have the same regularity.

Let $k\in \N$, $(\Tilde{X_1},\cdots,\Tilde{X_k})\in \mathfrak{a}^k$. Then $$D^k\Tilde{\psi}_\lambda(\Tilde{H})(\tilde{X_1},\cdots,\tilde{X_k})=\sum_{w\in W}\frac{\det(w)}{d_\lambda}\left(\prod_{j=1}^k (w(\lambda+\rho))(2X_j)\right)e^{(w(\lambda+\rho))(2H)}.$$Note that $\Vert w(\lambda+\rho)\Vert=\Vert \lambda+\rho\Vert$ for all $w\in W$. Furthermore, $(w(\lambda+\rho))(2H)$ is pure imaginary, so that $\vert e^{(w(\lambda+\rho))(2H)}\vert=1$. So we have \begin{equation}\label{majdiffcar}
    \Vert D^k\Tilde{\psi}_\lambda(\Tilde{H})\Vert \leq \frac{2^k\vert W\vert \Vert \lambda+\rho\Vert^k}{d_\lambda}.
\end{equation}So now, we have to study $d_\lambda$. By the Weyl dimension formula (\cite[Thm. 10.18]{hall2003lie}), we have $$d_\lambda=\prod_{\alpha\in \Phi^+} \frac{\langle\alpha,\lambda+\rho\rangle}{\langle\alpha,\rho\rangle}.$$Given that $\Delta$ is a basis of the finite-dimensional vector space $\mathfrak{a}_\C^*$, there exists $C>0$ such that for any $\lambda\in\mathfrak{a}_\C^*$, we have $C\Vert \lambda\Vert \leq \underset{1\leq j\leq \ell}{\max} \vert \langle\alpha_j,\lambda\rangle\vert$. For $\lambda\in \Lambda_G$ fixed, there exists $j(\lambda)$ such that $\langle\alpha_{j(\lambda)},\lambda+\rho\rangle=\underset{1\leq j\leq \ell}{\max} \langle\alpha_j,\lambda+\rho\rangle$. Note that since $\lambda+\rho$ is a dominant element, all these scalar products are non-negative integers. Now, for any $\alpha\in \Phi^+$, we have $$\langle\alpha,\lambda+\rho\rangle=\sum_{j=1}^\ell n_j(\alpha)\langle\alpha_j,\lambda+\rho\rangle \geq n_{j(\lambda)}(\alpha)\langle\alpha_{j(\lambda)},\lambda+\rho\rangle\geq Cn_{j(\lambda)}\Vert \lambda+\rho\Vert.$$By definition of $\gamma$, there are at least $\gamma$ positive roots $\alpha$ such that $n_{j(\lambda)}(\alpha)\geq 1$. Choose exactly $\gamma$ out of them. Let $D=\underset{\alpha}{\min}\langle\alpha,\rho\rangle>0$. For the $\vert \Phi^+\vert-\gamma$ remaining roots, we have $$\langle\alpha,\lambda+\rho\rangle\geq \langle\alpha,\rho\rangle\geq D.$$
Thus, we get that \begin{equation}\label{mindimension}d_\lambda \geq \frac{C^\gamma D^{\vert \Phi^+\vert-\gamma}\Vert \lambda+\rho\Vert^\gamma}{\prod_{\alpha\in \Phi^+}\langle\alpha,\rho\rangle}.\end{equation}Let $$C_k= \frac{2^k\vert W\vert\prod_{\alpha\in \Phi^+}\langle\alpha,\rho\rangle}{C^\gamma D^{\vert \Phi^+\vert-\gamma}}.$$Combining \eqref{majdiffcar} and \eqref{mindimension}, we have that for any $\lambda\in \Lambda_G$, $\Tilde{H}\in Q$, $$\Vert D^k\Tilde{\psi}_\lambda(\Tilde{H})\Vert \leq C_k\Vert \lambda+\rho\Vert^{k-\gamma}.$$Therefore, as soon as $k\leq \gamma$, the family of differentials of order $k$ are bounded in $\lambda$.
\end{proof}

\begin{coro}
Any $\Delta(G)$-bi-invariant matrix coefficient $\varphi$ of a unitary representation of $G\times G$ is such that $\varphi\circ \exp \in C^{(\gamma,0)}(Q)$.
\end{coro}
\begin{proof}It follows from Lemma \ref{roptid}.\end{proof}

We now prove that Theorem \ref{regcar} is optimal, in the sense that there are $\Delta(G)$-bi-invariant coefficients of unitary representations that are exactly of regularity $\gamma$.
\begin{theorem}
 \label{opticar}For any $0<\delta\leq 1$, the family $(\psi_\lambda)_{\lambda\in \Lambda_G}$ is not bounded in $C^{(\gamma,\delta)}(Q)$.
Hence, there exists a $K$-bi-invariant matrix coefficient of a unitary representation of $G$ that is not in $C^{(\gamma,\delta)}(Q)$.   
\end{theorem}
\begin{proof}As before, we will prove this for the family $(\Tilde{\psi}_\lambda)_{\lambda\in \Lambda_G}$. Assume that there is $\delta>0$ such that for any compact $L\subset Q$, there is $C_L>0$ such that for all $\Tilde{H},\Tilde{H'}\in L$, $\lambda\in \lambda_G$, \begin{equation}\label{hyp}\Vert D^\gamma\Tilde{\psi}_\lambda(\Tilde{H})-D^\gamma\Tilde{\psi}_\lambda(\Tilde{H'})\Vert \leq C_L\Vert \Tilde{H}-\Tilde{H'}\Vert^\delta.\end{equation}
To produce a contradiction, we want to restrict to a subfamily of $\lambda$ such that the previous estimates were sharp. To give a bound on $d_\lambda$, we used that there were at least $\gamma$ roots non-orthogonal to $\lambda$, but in general there might be more. So we will take $\lambda$ such that there are exactly $\gamma$ such roots.

Let $\lambda_1,\cdots,\lambda_\ell$ be the fundamental weights defined by $$\frac{2\langle\lambda_i,\alpha_j\rangle}{\langle\alpha_j,\alpha_j\rangle}=\delta_{i,j}.$$Let $i_0$ be such that $\gamma=\vert\{\alpha\in \Phi ^+\vert n_{i_0}(\alpha)\geq 1\}\vert$. Up to relabeling, we assume that $i_0=1$. Denote also $\beta_1,\cdots,\beta_\gamma$ the roots $\alpha$ such that $n_{1}(\alpha)\geq 1$, $B^+=\{\beta_1,\cdots,\beta_\gamma\}$ and $B=B^+\cup \left(-B^+\right)$. Then $\langle\lambda_1,\alpha\rangle\neq 0$ if and only if $\alpha\in B$. We will study the subfamily of functions associated to $\{n\lambda_1\}_{n\in \N}\subset \Lambda_G$. For any $\Tilde{X_1},\cdots,\Tilde{X_\gamma}\in \mathfrak{a}$ unit vectors, we have 
\begin{align*}
    \MoveEqLeft[8]\Vert D^\gamma\Tilde{\psi}_{n\lambda_1}(\Tilde{H})-D^\gamma\Tilde{\psi}_{n\lambda_1}(\Tilde{H'})\Vert & \\
     &  \geq \vert D^\gamma\Tilde{\psi}_{n\lambda_1}(\Tilde{H})(\Tilde{X_1},\cdots,\Tilde{X_\gamma})-D^\gamma\Tilde{\psi}_{n\lambda_1}(\Tilde{H'})(\Tilde{X_1},\cdots,\Tilde{X_\gamma})\vert \\
      & \geq \left\vert \sum_{w\in W} \frac{\det(w)}{d_{n\lambda_1}} \left(\prod_{j=1}^\gamma (w(n\lambda_1+\rho))(2X_j)\right)\right.\\
      & \phantom{\geq \quad\quad\quad\quad\quad}\left. \left[e^{(w(n\lambda_1+\rho))(2H)}-e^{(w(n\lambda_1+\rho))(2H')}\right]\right\vert.
\end{align*}

A root $\alpha$ is an element of $i\mathfrak{b}^*$. There is $Y_\alpha\in i\mathfrak{b}$ such that for any $\lambda\in i\mathfrak{b}^*$, $\langle\alpha,\lambda\rangle=\lambda(Y_\alpha)$. Choose $X_j=iY_{\beta_j}$, up to a normalisation constant. Thus, \begin{multline*}\Vert D^\gamma\Tilde{\psi}_{n\lambda_1}(\Tilde{H})-D^\gamma\Tilde{\psi}_{n\lambda_1}(\Tilde{H'})\Vert \geq \\C\left\vert \sum_{w\in W} \frac{\det(w)}{d_{n\lambda_1}} \left(\prod_{j=1}^\gamma \langle w(n\lambda_1+\rho),\beta_j\rangle\right) \left[e^{(w(n\lambda_1+\rho))(2H)}-e^{(w(n\lambda_1+\rho))(2H')}\right]\right\vert\end{multline*} for some constant $C>0$ depending only on the root system.

Let $W'=\{w\in W \vert w(B^+)\subset B\}$. Then $W'$ is a subgroup of $W$, and for $w\in W'$, $w(\beta_i)=\varepsilon_i(w)\beta_{i(w)}$, with $i\mapsto i(w)$ bijection of $\{1,\cdots,\gamma\}$ and $\varepsilon_i(w)\in \{\pm 1\}$. Denote also $\varepsilon(w)=\prod_i \varepsilon_i(w)$.

For $w\in W$, we have \begin{align*} \frac{\prod_{j=1}^\gamma \langle w(n\lambda_1+\rho),\beta_j\rangle}{d_{n\lambda_1}}   & =  \prod_{j=1}^\gamma \langle n\lambda_1+\rho,w^{-1}(\beta_j)\rangle \frac{\prod_{\alpha\in \Phi^+}\langle\rho,\alpha\rangle}{\prod_{\alpha\in \Phi^+}\langle n\lambda_1+\rho,\alpha\rangle} \\
     & =  \left(\prod_{j=1}^\gamma \langle\rho,\beta_j\rangle\right) \prod_{j=1}^\gamma\frac{\langle n\lambda_1+\rho,w^{-1}(\beta_j)\rangle}{\langle n\lambda_1+\rho,\beta_j\rangle}.
\end{align*}
Now, if $w\in W'$, $w^{-1}$ preserves $B^+$ up to signs, so we have $$\frac{\prod_{j=1}^\gamma \langle w(n\lambda_1+\rho),\beta_j\rangle}{d_{n\lambda_1}}=\varepsilon(w)\prod_{j=1}^\gamma \langle\rho,\beta_j\rangle.$$
Otherwise, if $w\in W\setminus W'$, there is $j$ such that $w^{-1}(\beta_j)\not\in B$ and so\linebreak $\langle n\lambda_1+\rho,w^{-1}(\beta_j)\rangle=\langle\rho,w^{-1}(\beta_j)\rangle$ is independent of $n$, and so $$\frac{\prod_{j=1}^\gamma \langle w(n\lambda_1+\rho),\beta_j\rangle}{d_{n\lambda_1}}=O\left(\frac{1}{n}\right).$$

We also have $w(n\lambda_1+\rho))(2H')\in i\R$ so $\vert e^{(w(n\lambda_1+\rho))(2H)}-e^{(w(n\lambda_1+\rho))(2H')}\vert\leq 2$. Thus, setting $C'=C\prod_{j=1}^\gamma \langle\rho,\beta_j\rangle$, there is $M>0$ such that we have for all $n\in \N^*$, $\Tilde{H},\Tilde{H'}\in Q$, \begin{align*}
    \MoveEqLeft[8] C'\left\vert \sum_{w\in W'} \det(w)\epsilon(w)\left[e^{(w(n\lambda_1+\rho))(2H)}-e^{(w(n\lambda_1+\rho))(2H')}\right]\right\vert & \\
    & \leq C\Bigg\vert \sum_{w\in W} \frac{\det(w)}{d_{n\lambda_1}} \left(\prod_{j=1}^\gamma \langle w(n\lambda_1+\rho),\beta_j\rangle\right) \\
    & \phantom{\leq \quad\quad\quad\quad\quad}\left[e^{(w(n\lambda_1+\rho))(2H)}-e^{(w(n\lambda_1+\rho))(2H')}\right]\Bigg\vert\\
    & \phantom{\leq} + C\Bigg\vert \sum_{w\in W\setminus W'} \frac{\det(w)}{d_{n\lambda_1}} \left(\prod_{j=1}^\gamma \langle w(n\lambda_1+\rho),\beta_j\rangle \right)\\
    & \phantom{\leq \quad\quad\quad\quad\quad}\left[e^{(w(n\lambda_1+\rho))(2H)}-e^{(w(n\lambda_1+\rho))(2H')}\right]\Bigg\vert\\
    & \leq \Vert D^\gamma\Tilde{\psi}_{n\lambda_1}(\Tilde{H})-D^\gamma\Tilde{\psi}_{n\lambda_1}(\Tilde{H'})\Vert +\frac{M}{n}
\end{align*}

We can further simplify the expression. Indeed, if $j\neq 1$, $w\in W'$, we have $w^{-1}(\alpha_j)\not\in B$ since $w(B)=B$. Thus, $$\langle w(\lambda_1),\alpha_j\rangle=\langle\lambda_1,w^{-1}(\alpha_j)\rangle=0,$$so $w(\lambda_1)\in \{\alpha_2,\cdots,\alpha_\ell\}^\bot=\R\lambda_1$. Furthermore, $w$ is unitary, so $w(\lambda_1)=s(w)\lambda_1$, $s(w)=\pm 1$. Now denote $$f(H)=\sum_{w\in W'}\det(w)\varepsilon(w)e^{w(p)(2H)}$$ and $$g(H)=i\sum_{w\in W'}\det(w)\varepsilon(w)s(w)e^{w(p)(2H)}.$$
Let $ix=\lambda_1(2H)$ and $iy=\lambda_1(2H')$. We have \begin{multline*}
    C'\left\vert \sum_{w\in W'} \det(w)\epsilon(w)\left[e^{(w(n\lambda_1+\rho))(2H)}-e^{(w(n\lambda_1+\rho))(2H')}\right]\right\vert  \\= C'\vert f(H)\cos(nx)+g(H)\sin(nx)-f(H')\cos(ny)-g(H')\sin(ny) \vert
\end{multline*}

Thus, we get that \begin{multline}
    \label{etape1}C'\vert f(H)\cos(nx)+g(H)\sin(nx)-f(H')\cos(ny)-g(H')\sin(ny) \vert \\\leq \Vert D^\gamma\Tilde{\psi}_{n\lambda_1}(\Tilde{H})-D^\gamma\Tilde{\psi}_{n\lambda_1}(\Tilde{H'})\Vert + \frac{M}{n}
\end{multline}so \begin{multline}
    \label{etape1'}\frac{C'^2}{2}\vert f(H)\cos(nx)+g(H)\sin(nx)-f(H')\cos(ny)-g(H')\sin(ny) \vert^2 \\ \leq \Vert D^\gamma\Tilde{\psi}_{n\lambda_1}(\Tilde{H})-D^\gamma\Tilde{\psi}_{n\lambda_1}(\Tilde{H'})\Vert^2 + \frac{M^2}{n^2}
\end{multline}
Let $m_0,N$ that will be chosen later. For any $m\geq m_0$, set $$S=S(m,N,H,H')=\sum_{n=m+1}^{m+N}\Vert D^\gamma\Tilde{\psi}_{n\lambda_1}(\Tilde{H})-D^\gamma\Tilde{\psi}_{n\lambda_1}(\Tilde{H'})\Vert^2.$$Since $\sum_{n=m+1}^{m+N}\frac{1}{n^2}\leq \sum_{n=m+1}^{+\infty}\frac{1}{n^2}=\frac{1}{m}$, from \eqref{etape1'} we get \begin{equation}
    \label{etape2} \frac{C'^2}{2}\sum_{n=m+1}^{m+N} \vert f(H)\cos(nx)+g(H)\sin(nx)-f(H')\cos(ny)-g(H')\sin(ny) \vert^2 \leq S+\frac{M^2}{m}
\end{equation}Since $\vert z\vert^2 \geq (\mathrm{Re}(z))^2$, if we denote $a_1=\frac{C'}{\sqrt{2}}\mathrm{Re}(f(H))$, $a_2=\frac{C'}{\sqrt{2}}\mathrm{Re}(g(H))$, \linebreak$b_1=\frac{C'}{\sqrt{2}}\mathrm{Re}(f(H'))$ and $b_2=\frac{C'}{\sqrt{2}}\mathrm{Re}(g(H'))$, as well as $u=(m+1)x$ and $v=(m+1)y$, we have \begin{equation}\label{etape3}
    \sum_{n=0}^{N-1} \vert a_1\cos(nx+u)+a_2\cos(nx+u+\frac{\pi}{2})-b_1\cos(ny+v)-b_2\cos(ny+v+\frac{\pi}{2})\vert^2 \leq S+\frac{M^2}{m}
\end{equation}

Since $\rho$ is strictly dominant, its images under the Weyl group are disjoint. Thus, the functions $H\mapsto e^{w(\rho)(H)}$ are linearly independent in $C(\mathfrak{b})$, so $f,g$ are non-zero. Furthermore, their extensions to $\mathfrak{b}_\C$ are holomorphic, so the zeros are isolated. Thus, we can find a compact ball $L$ in $Q$ such that $a_1^2+a_2^2+b_1^2+b_2^2\geq D>0$ for any $\Tilde{H},\Tilde{H'}\in L$, and $a_1,a_2,b_1,b_2$ bounded by $D'$. We can find a smaller ball $L'$ in $L$ such that $x,y,x-y,x+y$ are all in the same compact subset $J$ of $]0,\pi[$, up to $k\pi$, for any $\Tilde{H}\neq \Tilde{H'}\in L'$.

In that case, the same computations as in Lemma \ref{lemmeutile} but with more terms give that for $N\geq \frac{m_0}{\vert x-y\vert}$, where $m_0$ depends on $J,L$, we have \begin{multline}\label{etape4}\frac{1}{N} \sum_{n=0}^{N-1} \vert a_1\cos(nx+u)+a_2\cos(nx+u+\frac{\pi}{2})-b_1\cos(ny+v)-b_2\cos(ny+v+\frac{\pi}{2})\vert^2\\ \geq  \frac{a_1^2+a_2^2+b_1^2+b_2^2}{4}\geq \frac{D}{4}.\end{multline}
Thus, combining \eqref{etape3} and \eqref{etape4}, we have $S\geq \frac{DN}{4}-\frac{M^2}{m}$, so for $m$ large enough, we have that for any $\Tilde{H}\neq\Tilde{H'}\in L'$, \begin{equation}
    \label{etape5} S=\sum_{n=m+1}^{m+N}\Vert D^\gamma\Tilde{\psi}_{n\lambda_1}(\Tilde{H})-D^\gamma\Tilde{\psi}_{n\lambda_1}(\Tilde{H'})\Vert^2\geq \frac{DN}{8}.
\end{equation}

But by our first assumption in \eqref{hyp}, we have \begin{equation}
    S=\sum_{n=m+1}^{m+N}\Vert D^\gamma\Tilde{\psi}_{n\lambda_1}(\Tilde{H})-D^\gamma\Tilde{\psi}_{n\lambda_1}(\Tilde{H'})\Vert^2 \leq NC_{L'}^2 \Vert \Tilde{H}-\Tilde{H'}\Vert^{2\delta}.
\end{equation}

Thus, we get for any $\Tilde{H}\neq \Tilde{H'}\in L'$, $$\frac{D}{8}\leq C_{L'}^2\Vert \Tilde{H}-\Tilde{H'}\Vert^{2\delta}$$which gives a contradiction when $H\to H'$.

For the matrix coefficient, it then follows directly from Lemma \ref{roptid}.
\end{proof}

As in rank $1$, we end this subsection with a regularity result on Schur multipliers. Let $\lambda$ be the regular representation on $L^2(G\times G)$. We define a function on $\overline{Q}$ by $$T_{H}=\int\int_{\Delta(G)\times \Delta(G)}\lambda(k\exp(H)k')\,dk\,dk'\in B(L^2(G\times G)).$$
We want to study the regularity of $T$ as a map from $Q$ to $S_p(L^2(G\times G))$.
\begin{prop}
The map $T$ belongs to $C^{(\gamma-1,1)}(Q,S_\infty(L^2(G\times G)))$.
\end{prop}
\begin{proof}There is an orthonormal basis of $L^2(G\times G)$ such that for any $H\in Q$, $T_H$ is diagonal and its eigenvalues are $\psi_{{\lambda}}(H)$ for $\lambda\in \Lambda_G$. Then, by Theorem \ref{regcar} and Lemma \ref{existencederiv}, we get the result.\end{proof}
\begin{prop}
\label{schattencar}Let $p>2+\frac{\ell}{\gamma}$, $\gamma_p=\gamma-\frac{\ell+2\gamma}{p}$and  $$d=\left\lbrace \begin{aligned}\gamma_p & \textrm{ if }\gamma_p\not\in \Z\\ \gamma_p-\varepsilon& \textrm{ if }\gamma_p\in \Z\end{aligned}\right.$$ with $\varepsilon>0$ arbitrarily small. Then the map $T$ belongs to $C^{(r,\delta)}(Q,S_p(L^2(G\times G)))$ where $r=\lfloor d\rfloor$ and $\delta=d-r$.
\end{prop}

\begin{proof}There is an orthonormal basis such that for any $H\in Q$, $T_H$ is diagonal, with eigenvalues $\psi_{{\lambda}}(H)$ of multiplicity the dimension of the representation associated to $\Tilde{\lambda}\in \Lambda$ of $G\times G$. Since this representation is $V\otimes V^*$ where $V$ is the representation of $G$ associated to $\lambda\in \Lambda_G$, this multiplicity is $d_\lambda^2$. Thus, we have $$\Vert T_H\Vert_{S_p}^p=\sum_{\lambda\in \Lambda} d_\lambda^2 \vert \psi_{{\lambda}}(H)\vert^p \in [0,+\infty].$$

Let $k< \gamma$. Then by Theorem \ref{regcar} and Lemma \ref{existencederiv}, $D^kT$ exists and is the map from $Q$ such that for any $H\in Q$, $X_1,\cdots,X_k\in \mathfrak{a}$, the map $$D^kT(H)(X_1,\cdots,X_k)\in B(L^2(G\times G))$$ is diagonal in the previous basis with eigenvalues $D^k\psi_{{\lambda}}(H)(X_1,\cdots,X_k)$ and multiplicities $d_\lambda^2$. Thus, we have $$\Vert D^k T(H)\Vert_{L(\mathfrak{a}^{\otimes k},S_p)}^p\leq \sum_{\lambda\in \Lambda_G} d_\lambda^2 \Vert D^k\psi_{{\lambda}}(H)\Vert^p\in [0,+\infty].$$ 

Let $L$ be a compact subset of $Q$. Up to replacing $L$ by its convex hull, which is still compact by Carathéodory's theorem and a subset of $Q$, since $Q$ itself is convex, we can assume that $L$ is convex.  Let $H\in L$. By \eqref{majdiffcar}, there is $C_L>0$ such that $$d_\lambda^2 \Vert D^k\psi_{{\lambda}}(H)\Vert^p\leq C_L \Vert \lambda+\rho\Vert^{pk}d_\lambda^{2-p}.$$
Since by assumption, $2-p<0$, by \eqref{mindimension}, there exists a constant $C'_L>0$ such that $$d_\lambda^2 \Vert D^k\psi_{{\lambda}}(H)\Vert^p\leq C'_L\Vert\lambda+\rho\Vert^{pk+(2-p)\gamma}.$$

Given $\lambda_1,\cdots,\lambda_\ell$ the fundamental weights, $\lambda=\sum_{i=1}^\ell n_i(\lambda)\lambda_i$ with $n_i(\lambda)\in \N$ and $\rho=\sum_{i=1}^\ell \lambda_i$. Then, by equivalence of norms, there are $c,C>0$ such that $$c\sqrt{\sum_{i=1}^\ell (n_i(\lambda)+1)^2}\leq \Vert \lambda+\rho\Vert \leq C\sqrt{\sum_{i=1}^\ell (n_i(\lambda)+1)^2}$$

Thus, there is $M_L>0$ such that $$\Vert D^k T(H)\Vert^p\leq M_L\sum_{n_1,\cdots,n_\ell\geq 1} \frac{1}{\left(\sum n_i^2\right)^{\frac{(p-2)\gamma-pk}{2}}}.$$

We know that the sum on the right converges whenever $$(p-2)\gamma-pk>\ell,$$which is equivalent to $$\gamma_p=\gamma-\frac{\ell+2\gamma}{p}>k.$$In particular, this is the case for $k=r=\lfloor d\rfloor$, and we get that $\Vert D^rT(H)\Vert$ is bounded by $\Tilde{M}_L$ on $L$.

Now consider $H,H'\in L$, doing the same thing we get $$\Vert D^rT(H)-D^rT(H')\Vert^p \leq \sum_{\lambda\in \Lambda_G} d_\lambda^2 \Vert D^r\psi_{{\lambda}}(H)-D^r\psi_{{\lambda}}(H')\Vert^p.$$

We use two different bounds. For the first one, by triangle inequality, we have $$\Vert D^r\psi_{{\lambda}}(H)-D^r\psi_{{\lambda}}(H')\Vert\leq \Vert D^r\psi_{{\lambda}}(H)\Vert +\Vert D^r\psi_{{\lambda}}(H')\Vert$$so by \eqref{majdiffcar} and \eqref{mindimension} again, there are $C_1',C_1>0$ such that$$d_\lambda^2 \Vert D^r\psi_{{\lambda}}(H)-D^r\psi_{{\lambda}}(H')\Vert^p\leq C'_1\Vert \lambda+\rho\Vert^{pr-(p-2)\gamma}\leq \frac{C_1}{\left(\sum n_i^2\right)^{\frac{(p-2)\gamma-pr}{2}}}.$$We use this bound when $\left(\sum n_i^2\right)^{1/2}>\frac{1}{\Vert H-H'\Vert}$.

For the second one, since $L$ is convex, we get by the mean value theorem, $$\Vert D^r\psi_{{\lambda}}(H)- D^r\psi_{{\lambda}}(H')\Vert \leq \sup_L \Vert D^{r+1}\psi_{\lambda}(x)\Vert \Vert H-H'\Vert$$so there are $C_2',C_2>0$ such that \begin{multline*}d_\lambda^2 \Vert D^r\psi_{{\lambda}}(H)-D^r\psi_{{\lambda}}(H')\Vert^p\leq C'_2\Vert \lambda+\rho\Vert^{p(r+1)-(p-2)\gamma}\Vert H-H'\Vert^p\\ \leq \frac{C_2}{\left(\sum n_i^2\right)^{\frac{(p-2)\gamma-p(r+1)}{2}}}\Vert H-H'\Vert^p.\end{multline*}We use this bound when $\left(\sum n_i^2\right)^{1/2}\leq\frac{1}{\Vert H-H'\Vert}$.

Assume here that $H,H'$ are close enough, so that $\frac{1}{\Vert H-H'\Vert}>\sqrt{\ell}$.

Let $\kappa=(p-2)\gamma-pr>\ell$. Let $\underline{\mathbf{n}}=(n_1,\cdots,n_\ell)\in \N^\ell$, $\underline{\mathbf{n}}\neq 0$. Then for any $x$ in the cube $C(\underline{\mathbf{n}})=\prod_{i=1}^\ell [n_i,n_{i+1}]$, we have $$\frac{1}{\Vert \underline{\mathbf{n}}+\underline{\mathbf{1}}\Vert^\kappa}\leq \frac{1}{\Vert x\Vert^\kappa}\leq \frac{1}{\Vert\underline{\mathbf{n}}\Vert^\kappa}$$so\begin{equation}\label{sommeriemann}\frac{1}{\Vert \underline{\mathbf{n}}+\underline{\mathbf{1}}\Vert^\kappa}\leq \int_{C(\underline{\mathbf{n}})} \frac{1}{\Vert x\Vert^\kappa}\,dx\leq \frac{1}{\Vert\underline{\mathbf{n}}\Vert^\kappa}.\end{equation} 
Since $\bigcup_{\Vert\underline{\mathbf{n}}+\underline{\mathbf{1}}\Vert >\frac{1}{\Vert H-H'\Vert}} C(\underline{\mathbf{n}})\subset \R^\ell\setminus B(0,\frac{1}{\Vert H-H'\Vert}-\sqrt{\ell})=\R^\ell\setminus B$, we get that \begin{align*}
    \sum_{\Vert \underline{\mathbf{n}}\Vert>\frac{1}{\Vert H-H'\Vert}} \frac{1}{\Vert \underline{\mathbf{n}}\Vert^\kappa}& \leq \int_{\R^\ell\setminus B}\frac{1}{\Vert x\Vert^\kappa}\,dx \\
     & \leq  D_1 \int_{\frac{1}{\Vert H-H'\Vert}-\sqrt{\ell}}^\infty \frac{r^{\ell-1}}{r^\kappa}\,dr\\
     & \leq D_2 \left(\frac{1}{\Vert H-H'\Vert}-\sqrt{\ell}\right)^{\ell-\kappa}\\
     & \leq D_3 \Vert H-H'\Vert^{\kappa-\ell}
\end{align*}
where $D_1,D_2,D_3>0$ are constants independent of $H,H'$. Now since $\kappa-\ell=p(\gamma_p-r)$, we get that in the first case, there is $\Tilde{C}_1$ such that
\begin{equation}
    \label{case1schattencar}\sum_{\lambda,\underline{\mathbf{n}}>\frac{1}{\Vert H-H'\Vert}} \Vert D^r\psi_{{\lambda}}(H)-D^r\psi_{{\lambda}}(H')\Vert^p \leq \Tilde{C}_1 \Vert H-H'\Vert^{p(\gamma_p-r)}
\end{equation}

Now we turn to the second case. Let $\kappa=(p-2)\gamma-p(r+1)=p(\gamma_p-r)-p+\ell$. By the choice of $r$, we have $\kappa\leq \ell$.

We divide again into three subcases.

First, if $0<\kappa< \ell$, then we can use again \eqref{sommeriemann}. This time, $\frac{1}{\Vert x\Vert^\kappa}$ is integrable at $0$. We have that $\bigcup_{\Vert\underline{\mathbf{n}}+\underline{\mathbf{1}}\Vert \leq\frac{1}{\Vert H-H'\Vert}} C(\underline{\mathbf{n}})\subset B(0,\frac{1}{\Vert H-H'\Vert})$, so we get that there are $M_1,M_2>0$ such that\begin{align*}
     \sum_{n_i\geq 1,\Vert \underline{\mathbf{n}}\Vert\leq\frac{1}{\Vert H-H'\Vert}} \frac{1}{\Vert \underline{\mathbf{n}}\Vert^\kappa}& \leq \int_{B(0,\frac{1}{\Vert H-H'\Vert})}\frac{1}{\Vert x\Vert^\kappa}\,dx \\
     & \leq  M_1\int_0^\frac{1}{\Vert H-H'\Vert} \frac{r^{l-1}}{r^\kappa}\,dr\\
     & \leq  M_2 \left(\frac{1}{\Vert H-H'\Vert}\right)^{\ell-\kappa}\\
     & \leq  M_2 \Vert H-H'\Vert^{p(\gamma_p-r)-p}.
\end{align*}

If $\kappa\leq 0$, we do the same thing but the inequalities in \eqref{sommeriemann} are reversed, so we must integrate on $B(0,\frac{1}{\Vert H-H'\Vert}+\sqrt{\ell})$. We still end up with $$\sum_{n_i\geq 1,\Vert \underline{\mathbf{n}}\Vert\leq\frac{1}{\Vert H-H'\Vert}} \frac{1}{\Vert \underline{\mathbf{n}}\Vert^\kappa} \leq M_3 \Vert H-H'\Vert^{p(\gamma_p-r)-p}$$for some constant $M_3>0$.

In these two subcases, corresponding to $\gamma_p\not\in \Z$, we get that there is $\Tilde{C}_2>0$ such that 
\begin{equation}
    \label{case2schattencar}\sum_{\lambda,\underline{\mathbf{n}}\leq\frac{1}{\Vert H-H'\Vert}} \Vert D^r\psi_{{\lambda}}(H)-D^r\psi_{{\lambda}}(H')\Vert^p \leq \Tilde{C}_2 \Vert H-H'\Vert^{p(\gamma_p-r)-p}\Vert H-H'\Vert^p
\end{equation}

So combining \eqref{case1schattencar} and \eqref{case2schattencar}, we have for $H,H'\in L$ such that $\Vert H-H'\Vert< \frac{1}{\sqrt{\ell}}$ that there is $D_L>0$ such that $$\Vert D^rT(H)-D^rT(H')\Vert \leq D_L\Vert H-H'\Vert^{\gamma_p-r}.$$
Taking $D'_L=\max(D_L,2\Tilde{M}_L(\sqrt{\ell})^{\gamma_p-r})$, we have for any $H,H'\in L$, $$\Vert D^rT(H)-D^rT(H')\Vert \leq D'_L\Vert H-H'\Vert^{\gamma_p-r}.$$which completes the proof in the case $\gamma_p\not\in \Z$.

Finally, the last of the three subcases is $\kappa=\ell$, which is equivalent to $\gamma_p\in\Z$. Here, we use again \eqref{sommeriemann}, but we cannot integrate in $0$. So there are constants $E_1,E_2>0$ such that \begin{align*}
    \displaystyle\sum_{n_i\geq 1,\Vert \underline{\mathbf{n}}\Vert\leq\frac{1}{\Vert H-H'\Vert}} \frac{1}{\Vert \underline{\mathbf{n}}\Vert^\ell} & \leq \frac{1}{\sqrt{\ell}}+\int_{B(0,\frac{1}{\Vert H-H'\Vert})\setminus B(0,1)}\frac{1}{\Vert x\Vert^\ell}\,dx \\
     & \leq \frac{1}{\sqrt{\ell}}+E_1\int_1^\frac{1}{\Vert H-H'\Vert} \frac{r^{\ell-1}}{r^\ell}\,dr\\
     & \leq \frac{1}{\sqrt{\ell}}+E_1\vert\ln\Vert H-H'\Vert\vert.
\end{align*}
Thus there are constants $\Tilde{C}_3,\Tilde{C}_4>0$ such that\begin{equation}
    \label{case2schattencarbis}\sum_{\lambda,\underline{\mathbf{n}}\leq\frac{1}{\Vert H-H'\Vert}} \Vert D^r\psi_{{\lambda}}(H)-D^r\psi_{{\lambda}}(H')\Vert^p \leq \Tilde{C}_3 \Vert H-H'\Vert^p+\Tilde{C}_4 \Vert H-H'\Vert^p\vert\ln\Vert H-H'\Vert\vert
\end{equation}
So combining \eqref{case1schattencar} (with the exponent which reduces to $p$) and \eqref{case2schattencarbis}, we have for $H,H'\in L$ such that $\Vert H-H'\Vert< \frac{1}{\sqrt{\ell}}$ that there are $D_L,D_{L,\varepsilon}>0$ such that$$\Vert D^rT(H)-D^rT(H')\Vert \leq D_L\Vert H-H'\Vert\vert\ln\Vert H-H'\Vert \vert\leq D_{L,\varepsilon}\Vert H-H'\Vert^{1-\varepsilon}$$for any $1>\varepsilon>0$.
Again, we get the result for any $H,H'\in L$ up to changing the constant $D_{L,\varepsilon}$, and the proof is complete in the case $\gamma_p\in \Z$.
\end{proof}
\begin{remark}
In Theorem \ref{opticar}, we used a specific subfamily to prove optimality. Here, we consider a sum and not a supremum, so we cannot restrict to such a subfamily. The bound we use $d_\lambda^2\vert \psi_\lambda\vert^p$ is only tight for $\lambda=n\lambda_1$, where $\lambda_1$ is as in the proof of Theorem \ref{opticar}. In general, we can find a better bound for $d_\lambda^2\vert \psi_\lambda\vert^p$ if we know the number of roots orthogonal to $\lambda$. Thus, getting an optimal result would require a deeper study of the root system.    
\end{remark}

\subsection{Results for some higher rank symmetric spaces}\label{sectionpartialresult}
\subsubsection{Complex Grassmannians}
\label{sectioncomplexgrass}Let $q\geq p \geq 2$, $G=SU(p+q)$ and a subgroup $K=S(U(p)\times U(q))$ - the case $p=1$ was treated in the previous sections. Let $M=G/K$. Then $M$ is isomorphic to the Grassmann manifold of $p$-dimensional subspaces of $\C^{p+q}$. We have $\dim M=2pq$ and $\rank M=p$. Let $k=q-p$.

As before, let $P_n^{(\alpha,\beta)}$ denote the Jacobi polynomials and let $$\Tilde{P}_n^{(\alpha,\beta)}=\frac{\Gamma(\alpha+1)\Gamma(n+1)}{\Gamma(\alpha+n+1)}P_n^{(\alpha,\beta)}$$be the Jacobi polynomials normalised at $1$.

The spherical functions of $(G,K)$ can be found in \cite[Section 3]{camporesi_spherical_2006}. For a suitable choice of basis and of positive roots, the Weyl chamber is $$C=\{(x_1,\cdots,x_p)\vert x_1>\cdots x_p>0\}.$$Consider $Q\subset C$ as before, open in $\R^p$. Define $$w(X)=\prod_{i<j}\left(\cos(x_i)-\cos(x_j)\right),$$ $$c=2^{p(p-1)/2}\prod_{j=1}^{p-1}j!(j+k)^{p-j}$$and$$c(n)=\left(n+\frac{k+1}{2}\right)^2-\left(\frac{k+1}{2}\right)^2.$$
Then if $\mu=\sum r_i\mu_i\in \Lambda$, let $m_i=r_i+\cdots+r_{p}$, so that $m_1\geq m_2\geq \cdots \geq m_p$, and $n_i=m_i+p-i$, so $n_1>\cdots>n_p$. For $X\in Q$, we have $$\varphi_\mu(\exp(iX))=c\frac{\det \left(\Tilde{P}_{n_i}^{(k,0)}(\cos(x_j))\right)}{w(X)\prod_{i<j} c(n_i)-c(n_j)}.$$

Since $w(X)\neq 0$ for all $x\in Q$, both $\frac{c}{w}$ and $\frac{w}{c}$ are smooth functions, thus the optimal regularity of the family $(\varphi_\mu)$ is the same as the optimal regularity of the family $(\frac{w}{c}\varphi_\mu)$. Using Lemma \ref{precomposition} and since $c(n_i)-c(n_j)=(n_i+n_j+k+1)(n_i-n_j)$, we want to study the optimal regularity of the functions $$\fonction{\psi_\mu}{\Tilde{Q}}{\C}{(t_1,\cdots,t_p)}{\frac{\det \left(\Tilde{P}_{n_i}^{(k,0)}(t_j)\right)}{\prod_{i<j}(n_i+n_j+k+1)(n_i-n_j)}}.$$

We have $\psi_\mu(t)=\frac{\sum_{\sigma \in \mathfrak{S}_p} \varepsilon(\sigma)\prod_{i=1}^{p}\Tilde{P}^{(k,0)}_{n_i}(t_{\sigma(i)})} {\prod_{i<j}(n_i+n_j+k+1)(n_i-n_j)}$. Let $D(l_1,\cdots,l_p)f=\frac{\partial^{\sum l_i}f}{\partial_{x_1}^{l_1}\cdots \partial_{x_p}^{l_p}}$. Then $$D(l_1,\cdots,l_p)\psi_{\mu}(t)=\frac{\sum_{\sigma \in \mathfrak{S}_p} \varepsilon(\sigma)\prod_{i=1}^{p}\partial^{l_{\sigma(i)}}\Tilde{P}_{n_i}^{(k,0)}(t_{\sigma(i)})} {\prod_{i<j}(n_i+n_j+k+1)(n_i-n_j)}.$$

\begin{theorem}
\label{roptGrassC2}Let $M=SU(p+q)/S(U(p)\times U(q))$, with $q\geq p\geq 2$. Set $$\alpha_\infty=\left\lbrace\begin{aligned}2 & \textrm{ if }p=q=2\\\ q+p-\frac{3}{2}& \textrm{ else}\end{aligned}\right..$$Then we have $$r_{opt}(M)\geq (\lfloor \alpha_\infty\rfloor,\alpha_\infty-\lfloor \alpha_\infty\rfloor).$$    
\end{theorem}

\begin{proof}Let $L$ be a compact subset of $\Tilde{Q}$, which we again assume to be convex. Then each $t_i$ is in a compact $L'$ of $]-1,1[$.  From inequality \ref{majderivphi} in Section \ref{mainsection}, we know that there is a constant $C=C(L',k,l)>0$ such that $\forall n\in\N,\forall t\in L'$, we have $$\vert \partial^{l}\Tilde{P}_{n}^{(k,0)}(t) \vert \leq Cn^{l-k-\frac{1}{2}}.$$
Thus there is some constant $C=C(L,k,l_1,\cdots,l_p)$ such that for all $\mu\in \Lambda$ and $t\in L$ \begin{equation}\label{majpartielbase}\vert D(l_1,\cdots,l_p)\psi_{\mu}(t) \vert \leq C \frac{\sum_{\sigma \in \mathfrak{S}_p}\prod_{i=1}^{p} n_i^{l_{\sigma(i)}-k-\frac{1}{2}}} {\prod_{i<j}(n_i+n_j+k+1)(n_i-n_j)}.\end{equation}

Let $I=\left\{i\vert n_{i+1}<\frac{n_i}{2}\right\}$ and $J=\llbracket 1,p-1\rrbracket\backslash I$. Let $i_0=\min I$. For $j>i_0\geq i$, we have $$n_j\leq n_{i_0+1}<\frac{n_{i_0}}{2}\leq \frac{n_i}{2}$$so $(n_i+n_j+k+1)(n_i-n_j)\geq \frac{n_i^2}{2}$. On the other hand, for $i<i_0$, there is $C_a>0$ such that $n_{i+1}^a\leq C_a n_1^a$ for any $a$. Thus, if $l_1+\cdots+l_p=s$, \begin{align*}\vert D(l_1,\cdots,l_p)\psi_\mu(t)\vert & \leq  C \sum_{\sigma\in \mathfrak{S}_p} n_1^s  n_1^{\sum_{i\leq i_0} \left(-k-\frac{1}{2}-(p-i)-(p-i_0)\right)}\prod_{i>i_0} n_i^{-k-\frac{1}{2}-(p-i)}\\ 
& \leq \Tilde{C}_{i_0} n_1^{s-\left(i_0k+\frac{i_0}{2}+i_0p-\frac{i_0(i_0+1)}{2}+i_0(p-i_0)\right)}\\
& \leq \Tilde{C}_{i_0} n_{1}^{s-\left(i_0(q+p)-\frac{3}{2}i_0^2\right)}
\end{align*} for some constant $\Tilde{C}_{i_0}>0$.
Let $\kappa=\underset{i\in\llbracket 1,p\rrbracket}{\min} i(p+q)-\frac{3}{2}i^2$ and $C'=\underset{i\in\llbracket 1,p\rrbracket}{\max} \Tilde{C}_{i}$. Then $\forall \mu\in \Lambda$ and $\forall t\in L$, we have $$\vert D(l_1,\cdots,l_p)\psi_\mu(t)\vert \leq C'n_1^{s-\kappa}.$$Thus, we only need to compute $\kappa$. Let $f:x\mapsto (p+q)x-\frac{3}{2}x^2$. If $q\geq 2p$, $f$ is increasing on $[1,p]$ so its minimum is $f(1)$ and $\kappa=f(1)=q+p-\frac{3}{2}$. If $q\leq 2p$, $f$ is increasing on $[1,\frac{p+q}{3}]$ and decreasing on $[\frac{p+q}{3},p]$. Thus, $\kappa$ is either $f(1)$ or $f(p)$. But $$\begin{array}{rrcl}
   & f(p) & < & f(1) \\
    \Longleftrightarrow & pq-\frac{p^2}{2} & < & p+q-\frac{3}{2}\\
    \Longleftrightarrow & (p-1)q & < &(p-1)\frac{p+3}{2}\\
    \Longleftrightarrow & q & < & \frac{p+3}{2}.
\end{array}$$ But we also have $q\geq p$, thus this last inequality implie $\frac{p+3}{2}>p$ so $p=2$ and then $q=2$. Thus, if $p=q=2$, $\kappa=f(p)=2$. Otherwise, $\kappa=f(1)=q+p-\frac{3}{2}$.

Now, if we take $\Vert .\Vert_1$ on $\R^p$, this tells us that $$\Vert D^s\psi_\mu(t)\Vert \leq \Tilde{C}_s n_1^{s-\kappa}$$so that if $s\leq \kappa$, $D^s\psi_\mu$ is bounded on $L$ independently of $\mu\in \Lambda$.

If $p=q=2$, $\kappa=2$ is an integer so the proof is complete.

Else, $\kappa=r+\frac{1}{2}$, $r=q+p-2$. If $x,y\in L$, on the one hand, $$\Vert D^r\psi_\mu(x)-D^r\psi_\mu(y)\Vert\leq 2\Tilde{C}_rn^{-1/2}.$$On the other hand, we get $$\Vert D^r\psi_\mu(x)-D^r\psi_\mu(y)\Vert\leq \Tilde{C}_{r+1}n^{1/2}\Vert x-y\Vert$$thus, \begin{multline*}
    \Vert D^r\psi_\mu(x)-D^r\psi_\mu(y)\Vert\leq \left(2\Tilde{C}_rn^{-1/2}\right)^{1/2}\left(\Tilde{C}_{r+1}n^{1/2}\Vert x-y\Vert\right)^{1/2} \\ \leq \sqrt{2\Tilde{C}_r\Tilde{C}_{r+1}}\Vert x-y\Vert^{1/2}
\end{multline*}so that $D^r\psi_\mu$ is $\frac{1}{2}$-Hölder on $L$ with a constant independent on $\mu\in \Lambda$, so we get the result.
\end{proof}

\begin{remark}
If in the definition of $\psi_\lambda$ we replace $(k,0)$ with $(\alpha,\beta)\in (\R_+)^2$, we get the same result with $$\kappa=\left\lbrace\begin{aligned}2\alpha+2 & \textrm{ if }p=2,\alpha <\frac{1}{2}\\\ \alpha+2p-\frac{3}{2}& \textrm{ else}\end{aligned}\right..$$
\end{remark}

\subsubsection{Some related spaces}
Consider $\Z^\ell$ with lexicographic order. A polynomial $P$ in $\ell$ variables has degree $n$ if $P=\sum_{m\leq n}c_mx^m$, $c_n\neq 0$. If $P$ is symmetric of degree $n$, then $n_1\geq n_2\geq \cdots n_l$. Let $\Omega=[-1,1]^l$, $\alpha,\beta >-1$ and $\gamma \geq -1/2$. Define on $\Omega$ the function $$w_{\alpha,\beta,\gamma}(x)=\prod_{i=1}^\ell (1-x_i)^\alpha(1+x_i)^{\beta}\prod_{i<j}(x_i-x_j)^{2\gamma+1}.$$If $\ell=1$, this is the Jacobi weight of parameter $(\alpha,\beta)$.
\begin{definition}
The polynomials $\left(P^{(\alpha,\beta,\gamma)}_n\right)_{n\in \N}$ in $\ell$ variables are the unique polynomials defined by \begin{itemize}
    \item $P^{(\alpha,\beta,\gamma)}_0=1$,
    \item $P^{(\alpha,\beta,\gamma)}_n$ is symmetric of degree $n$ and dominant coefficient $1$,
    \item $\int_{\Omega} P^{(\alpha,\beta,\gamma)}_n(x)Q(x)w_{\alpha,\beta,\gamma}(x)\,dx=0$ for all $Q$ symmetric of degree $q<n$.
\end{itemize}
\end{definition}

These polynomials are studied in \cite{vretare}. For certain values, these polynomials can be related to Jacobi polynomials (\cite[Thm. 4.5 and 4.6]{vretare}).
\begin{prop}
For $\gamma=-\frac{1}{2}$, $$P^{(\alpha,\beta,\gamma)}_n(x)=\sum_{\sigma\in\mathfrak{S}_l} P^{(\alpha,\beta)}_{n_1}(x_{\sigma(1)})\cdots P^{(\alpha,\beta)}_{n_l}(x_{\sigma(\ell)}).$$
\end{prop}
\begin{prop}
Let $A^{\alpha,\beta}_x(m)=\det\left(P^{(\alpha,\beta)}_{m_i}(x_j)\right)$. Then $$P^{(\alpha,\beta,1/2)}_n(x)=\frac{A^{(\alpha,\beta)}_x(n_1+\ell-1,n_2+\ell-2,\cdots,n_\ell)}{A^{(\alpha,\beta)}_x(\ell-1,\ell-2,\cdots,0)}.$$
\end{prop}
Thus, we can see that, once normalised by $1$ at $1$, the family of $P_n^{(k,0,1/2)}$ are the spherical functions of $(SU(k+2l),S(U(k+l)\times U(l)))$. It turns out that more families of spherical functions appears as polynomials of this kind (\cite[Thm. 4.2]{vretare}). 
\begin{theorem}
\label{spheriquevretare}Let $\psi_n^{(\alpha,\beta,\gamma)}=\frac{P^{(\alpha,\beta,\gamma)}_n}{P^{(\alpha,\beta,\gamma)}_n(1)}$. Then the functions $\psi_n^{(\alpha,\beta,\gamma)}$ are the spherical functions of $(G,K)$ for the following values of $\alpha,\beta,\gamma$:
$$\begin{array}{|c|c|c|c|c|}
    \hline
    (G,K) & \ell & \alpha & \beta & \gamma \\
    \hline
    (SO(p+q),S(O(p)\times O(q))) & p & (q-p-1)/2 & -1/2 & 0\\
    \hline
    (SU(p+q),S(U(p)\times U(q))) & p & q-p & 0 & 1/2\\
    \hline
    (Sp(p+q),Sp(p)\times Sp(q)) & p &  2(q-p)+1 & 1 & 3/2\\
    \hline
    (Sp(k),U(k)) & k & 0 & 0 & 0\\
    \hline
    (SO(4k),U(2k)) & k & 0 & 0 & 3/2\\
    \hline 
    (SO(4k+2),U(2k+1) & k & 2 & 0 & 3/2\\
    \hline
    (Sp(k)\times Sp(k),\Delta(Sp(k))) & k & 1/2 & 1/2 & 1/2\\
    \hline
    (SO(2k+1)\times SO(2k+1),\Delta(SO(2k+1))) & k & 1/2 & -1/2 & 1/2\\
    \hline
\end{array}$$
\end{theorem}
\begin{remark}
Adapting the previous work for general values of $\alpha\geq 0,\beta> -1$, we know that the optimal uniform regularity of the family $(\psi_n^{(\alpha,\beta,1/2)})$ is at least $\alpha+\frac{1}{2}+2(\ell-1)$. From the previous table, this recovers the regularity found in Section \ref{caracteres} for $(Sp(k)\times Sp(k),\Delta(Sp(k)))$ and $(SO(2k+1)\times SO(2k+1),\Delta(SO(2k+1))$.
\end{remark}

In \cite[Thm. 5.1]{vretare}, Vretare shows differents formulas expressing $\psi_n^{(\alpha,\beta,\gamma)}$ as a linear combination of different $\psi_n^{(\alpha',\beta',\gamma')}$. If we were able to control the coefficient, we could investigate more families. It turns out that this is the case in $2$ variables. Set $\ell=2$. Vretare showed the following (\cite[Thm. 6.2]{vretare}):
\begin{theorem}Let $\alpha,\beta>-1$ and $\gamma\geq -1/2$. Let $n\geq m$. Then \begin{multline*}(x-y)^2\psi_{n,m}^{(\alpha,\beta,\gamma+1)}\\=b_{20}\psi_{n+2,m}^{(\alpha,\beta,\gamma)}+b_{10}\psi_{n+1,m}^{(\alpha,\beta,\gamma)}+b_{00}\psi_{n,m}^{(\alpha,\beta,\gamma)}+b_{11}\psi_{n+1,m+1}^{(\alpha,\beta,\gamma)}+b_{1-1}\psi_{n+1,m-1}^{(\alpha,\beta,\gamma)}\end{multline*}and there is a constant $C=C(\alpha,\beta,\gamma)>0$ such that $b_{ij}\leq Cn^{-1}(n-m)^{-1}$ for any $i,j$.
\end{theorem}
\begin{remark}
From this, in the case $\ell=2$, we can recover the regularity obtained in Theorem \ref{roptGrassC2} from the regularity at least $\alpha+\frac{1}{2}$ for the family $(\psi_n^{(\alpha,\beta,-1/2)})$ that is easy to compute.
\end{remark}

With Lemma \ref{multsmooth}, the regularity of the family $\left(\psi_{n,m}^{(\alpha,\beta,3/2)}\right)_{n\geq m}$ is the same as the regularity of the family $\left((x-y)^2\psi_{n,m}^{(\alpha,\beta,3/2)}\right)_{n\geq m}=\left(\phi_{n,m}^{(\alpha,\beta)}\right)$.

Set $L$ a compact subset of $]-1,1[^2$. Let $l_1,l_2\in \N$, $s=l_1+l_2$. Assume first that $m<\frac{n}{2}$. As in Theorem \ref{roptGrassC2}, there is a constant $C>0$ such that $$\vert D(l_1,l_2)\psi_{n,m}^{(\alpha,\beta,1/2)}(x)\vert \leq Cn^{s-\alpha-\frac{5}{2}}$$thus there is $C'>0$ such that $$\vert D(l_1,l_2)\phi_{n,m}^{(\alpha,\beta)}\vert \leq C'n^{s-\alpha-\frac{9}{2}}.$$On the contrary, assume that $n>m\geq \frac{n}{2}$, then there is a constant $D>0$, $$\vert D(l_1,l_2)\psi_{n,m}^{(\alpha,\beta,1/2}(x)\vert \leq Dn^{s-2\alpha-2}$$thus there is $D'>0$ such that $$\vert D(l_1,l_2)\phi_{n,m}^{(\alpha,\beta)}\vert \leq D'n^{s-2\alpha-3}.$$
Now depending on $(n,m)$, we have two possible upper bounds. We want to know which one is the worst, and thus works independently of $(n,m)$.\\We have $2\alpha+3\geq \alpha+\frac{3}{2}$ if and only if $\alpha \geq \frac{3}{2}$. Now, adapting again the argument from the proof of Theorem \ref{roptGrassC2}, we get that the regularity of the family $\left(\psi_{n,m}^{(\alpha,\beta,3/2)}\right)_{n\geq m}$ is at least $$r=\left\lbrace\begin{aligned}2\alpha+3 & \textrm{ if }\alpha < \frac{3}{2}\\ \alpha+\frac{9}{2}& \textrm{ if} \alpha\geq \frac{3}{2}\end{aligned}\right.$$

From Theorem \ref{spheriquevretare}, this gives lower bound for new pairs.
\begin{theorem}\label{highrankbis}Let $M=SO(8)/U(4)$, $\rank M=2$ and $\dim M=12$. We have $$r_{opt}(M)\geq (3,0).$$
Let $M=SO(10)/U(5)$, $\rank M=2$ and $\dim M=20$. We have $$r_{opt}(M) \geq (6,\frac{1}{2}).$$
Let $n\geq 4$, $M=Sp(n)/(Sp(2)\times Sp(n-2))$, $\rank M=2$ and $\dim M=8(n-2)$. We have $$r_{opt}(M)\geq \left\lbrace\begin{aligned}(2n-3,\frac{1}{2})&\textrm{ if }n>4\\(5,0) & \textrm{ if }n=4 \end{aligned}\right.$$
\end{theorem}

\subsection{A conjecture on the optimal regularity}\label{sectionconj}
In this section, we will see on an example that there is a difference in the estimates when $\mu$ is close to the walls of the Weyl chamber and away from the walls. We will use this to give a conjecture on the optimal regularity in the general case.

Let $G=SU(q+2)$ and $K=S(U(2)\times U(q))$. Let $M=G/K$ of rank $2$ and dimension $4q$. In Theorem \ref{roptGrassC2}, we showed that $r_{opt}(M)\geq (q,\frac{1}{2})$ if $q>2$, and $(2,0)$ otherwise. Furthermore, spherical functions are determined by the heighest weight $\mu=r_1\mu_1+r_2\mu_2$ of the associated representation of $G$. There is a basis $\{\alpha,\beta\}$ of the root system $\Sigma_\mathfrak{a}$ such that any positive root is a linear combination of $\alpha,\beta$ with positive integral coefficients. The walls of the Weyl chamber are defined by the hyperplanes $H_\alpha=\alpha^\bot$ and $H_\beta=\beta^\bot$. We have $\mu_1\in H_\alpha$ and $\mu_2\in H_\beta$. Finally, the positive roots are $\alpha,\beta,2\alpha,\alpha+\beta,\alpha+2\beta,2(\alpha+\beta)$ of multiplicities $2q-4,2,1,2q-4,2,1$ respectively (see \cite{camporesi_spherical_2006}).

Let $\mu_0=\mu_1+\mu_2$. Then, for $n\in \N$, $n\mu_0$ is in a cone with compact base in the open Weyl chamber, so away from the walls. Then, $$\psi_{n\mu_0}(x,y)=\frac{\Tilde{P}_{2n+1}^{(q-2,0)}(x)\Tilde{P}_{n}^{(q-2,0)}(y)-\Tilde{P}_{2n+1}^{(q-2,0)}(y)\Tilde{P}_{n}^{(q-2,0)}(x)}{(3n+q)(n+1)}.$$The family $\left(\psi_{n\mu_0}\right)_{n\in \N}$ is a subfamily of the family of all spherical functions, and we can show with the same estimates as in Theorem \ref{roptGrassC2} that this family is bounded in the Hölder space of regularity $2q-1$.

So, for this subfamily corresponding to highest weight away from the walls, the regularity is strictly better than the regularity of the whole family. Note that being in a cone with compact base means that $\langle n\mu_0,\alpha\rangle,\langle n\mu_0,\beta\rangle$ have the same growth rate as $\Vert n\mu_0\Vert$.

In fact, for any of the pairs we considered above, with $\mu_0=\sum_{i=1}^l \mu_i$, the subfamily corresponding to $\{n\mu_0\}$ is of regularity at least $\frac{\dim M-\rank M}{2}$.

This shows that it should be highest weights close to the walls that gives a bad regularity. Let us see what happens when the highest weight is exactly on a wall, which corresponds to $\mu$ being orthogonal to some roots. For this, set either $r_1=0$ or $r_2=0$.

If $r_1=0$, $r_2=n$, then $$\psi_{n}(x,y)=\frac{\Tilde{P}_{n+1}^{(q-2,0)}(x)\Tilde{P}_{n\mu_2}^{(q-2,0)}(y)-\Tilde{P}_{n+1}^{(q-2,0)}(y)\Tilde{P}_{n}^{(q-2,0)}(x)}{(2n+q)}$$and the uniform regularity of this subfamily is at least $2q-2$. Note that in that case, $\mu$ is orthogonal to $\beta$, which is the root of the basis with small multiplicity.

If $r_1=n-1,r_2=0$, since $\Tilde{P}_0^{(q-2,0)}=1$, the expression is easier. We have $$\psi_{n\mu_1}(x,y)=\frac{\Tilde{P}_n^{(q-2,0)}(x)-\Tilde{P}_n^{(q-2,0)}(y)}{(n+q+1)n}.$$Here, we get uniform regularity at least $q+\frac{1}{2}$, so this is where the worst happens (except for $q=2$). It is not clear whether this is the optimal regularity or not, but we will check that $(q+1)$-th differential is unbounded in $n$ - so the uniform regularity is at most $q+1$. Let $x=\cos \theta$. From Proposition \ref{derivee} and Proposition \ref{darboux}, we have \begin{align*}\frac{d^r}{dx^r}\Tilde{P}_n^{(q-2,0)}(x) & = \frac{\Gamma(q-1)\Gamma(n+1)}{\Gamma(q+n-1)} \frac{d^r}{dx^r} P_n^{(q-2,0)}(x)\\
& = \frac{\Gamma(q-1)\Gamma(n+1)}{\Gamma(q+n-1)} \frac{\Gamma(q+n+r-1)}{2^r\Gamma(q+n-1)} P_{n-r}^{(q-2+r,r)}(x)\\
&=n^{r-q+3/2}u(\theta)\cos(n\theta+\gamma)+O(n^{r-q+1/2})\end{align*} where $u$ is a smooth function of $]-1,1[$, $\gamma$ depends on $q,r,\theta$ but not on $n$. We used the same asymptotics on $\Gamma$ as in \eqref{gamma}. Also, $\Tilde{P}_n^{(q-2,0)}(y)=O(1)$. Now let $r=q+1$ and consider the partial derivative $$D(q+1,0)\psi_{n}(x,y)=n^{1/2}u(\theta)\cos(n\theta+\gamma)+O(n^{-1/2})$$which is unbounded in $n$. So we get that when $\mu$ is orthogonal to $\alpha$, the subfamily is bounded in $C^{(q,\frac{1}{2})}$ but not bounded in $C^{q+1}$. Thus, we can conclude that it has lower uniform regularity than the previous subfamilies. Note that $\alpha$ is the root with high multiplicity.

So this example shows that when $\mu$ is not orthogonal to any roots, spherical functions are well-behaved and we get good estimates as the spectral parameter grows. But when $\mu$ becomes orthogonal to some roots, we see that the spherical functions tend to be unbounded in more Hölder spaces. We also see that the optimal uniform regularity decreases with the sum of multiplicities of non-orthogonal roots.\smallskip

We now try to give a more quantitative interpretation of this. For a root $\alpha\in \Sigma_\mathfrak{a}$, let $m(\alpha)$ be its multiplicity. For $\mu\in \Lambda$, define $$S_\mu=\sum_{\alpha\in \Sigma_{\mathfrak{a}}^+,<\alpha,\mu>\neq 0}m(\alpha).$$
Our conjecture is that if we set$$r=\frac{1}{2}\underset{\mu\in \Lambda\setminus\{0\}}{\inf} S_\mu,$$ then $r_{opt}(M)=(\lfloor r\rfloor,r-\lfloor r\rfloor).$ 
  
More generally, the optimal uniform regularity of a subfamily of spherical functions indexed by $\Lambda'\subset \Lambda$ should be be given by a similar formula involving only weights of $\Lambda'$. However, removing a finite number of elements of $\Lambda'$ will not change the regularity of the family. Furthermore, we can notice with the example above that the behaviour of the family $\{(0,n)\}_n$ will be the same as the behaviour of $\{(k,n)\}_n$ for $k$ fixed, so what is important is not the orthogonality of the family with roots, but the boundedness of the scalar products, thus making the formula more complicated.
  
By \cite{vretare_orth}, there is a basis $\{\alpha_1,\cdots,\alpha_\ell\}$ associated to the positive root system $\Sigma_\mathfrak{a}^+$ such that the fundamental weights $\mu_i$ satisfy $$\frac{\langle\mu_i,\alpha_j\rangle}{\langle\alpha_j,\alpha_j\rangle}=\left\lbrace \begin{aligned}0 & \textrm{ if }i\neq j\\ 1& \textrm{ if }i=j,2\alpha_j\not\in \Sigma_\mathfrak{a}^+\\ 2& \textrm{ if }i=j,2\alpha_j\in \Sigma_\mathfrak{a}^+\end{aligned}\right.$$

Then we have $r=\underset{1\leq i \leq \ell}{\min} \sum_{<\alpha,\mu_i>\neq 0} \frac{m(\alpha)}{2}$. In Table \ref{tab:symtypegroup} and \ref{tab:valuesconj}, we compute the values of $r$ for all simply connected irreducible symmetric spaces of compact type, using the classification and multiplicities given in \cite[Ch. VII]{loos1969symmetric2}. We can verify that these values agree with the results found for rank $1$ (Theorem \ref{mainthm}) and with the lower bound found for some higher rank spaces in Theorem \ref{roptGrassC2} and Theorem \ref{highrankbis}.

\begin{table}[ht]
\caption{Values of $r$ for the irreducible symmetric spaces of the form $G\times G/\Delta(G)\simeq G$.}
\label{tab:symtypegroup}
\renewcommand{\arraystretch}{1.2}
\centering

\begin{tabular}{|c|c|c|c|}
\hline
$M$                 & $\dim M$  & $\rank M$ & $r$    \\ \hline
$SU(n), n\geq 2$    & $n^2-1$   & $n-1$     & $n-1$  \\ \hline
$SO(2n+1), n\geq 1$ & $n(2n+1)$ & $n$       & $2n-1$ \\ \hline
$Sp(n), n\geq 1$    & $n(2n+1)$ & $n$       & $2n-1$ \\ \hline
$SO(2n), n\geq 2$   & $n(2n-1)$ & $n$       & \begin{tabular}{@{}c@{}} $1$ if $n=2$\\$3$ if $n=3$ \\ $2n-2$ else\end{tabular} \\ \hline
$G_2$               & $14$      & $2$       & $5$    \\ \hline
$F_4$               & $52$      & $4$       & $15$   \\ \hline
$E_6$               & $78$      & $6$       & $16$   \\ \hline
$E_7$               & $133$     & $7$       & $27$   \\ \hline
$E_8$               & $248$     & $8$       & $57$   \\ \hline
\end{tabular}

\end{table}

\begin{table}[ht]
\caption{Values of $r$ for irreducible symmetric spaces.}
\label{tab:valuesconj}
\renewcommand{\arraystretch}{1.3}
\centering
\centerline{\begin{tabular}{|cc|c|c|c|}
\hline
\multicolumn{2}{|c|}{$M$}                                                                  & $\dim M$        & $\rank M$ & $r$ \\ \hline
\multicolumn{1}{|c|}{$AI$}                    & $SU(n)/SO(n),n\geq 2$                              & $\frac{(n-1)(n+2)}{2}$ & $n-1$              &    $\frac{n-1}{2}$ \\ \hline
\multicolumn{1}{|c|}{$AII$}                   & $SU(2n)/Sp(n),n\geq 2$                             & $(n-1)(2n+1)$          & $n-1$              &  $2(n-1)$   \\ \hline
\multicolumn{1}{|c|}{$AIII$}                  & $SU(p+q)/S(U(p)\times U(q)), p+q\geq 3$               & $2pq$                  & $\min(p,q)$        &  \begin{tabular}{@{}c@{}}$2$ if $p=q=2$ \\ $p+q-\frac{3}{2}$ else\end{tabular}   \\ \hline
\multicolumn{1}{|c|}{$BDI$}                   & $SO(p+q)/SO(p)\times SO(q),p+q\geq 3$ & $pq$                   & $\min(p,q)$        &   \begin{tabular}{@{}c@{}}$\frac{1}{2}$ if $p=q=2$ \\ $\frac{3}{2}$ if $p=q=3$\\$\frac{p+q}{2}-1$ else\end{tabular}  \\ \hline
\multicolumn{1}{|c|}{$CI$}                    & $Sp(n)/U(n),n\geq 1$                               & $n(n+1)$               & $n$                &   $n-\frac{1}{2}$  \\ \hline
\multicolumn{1}{|c|}{$CII$}                   & $Sp(p+q)/Sp(p)\times Sp(q),p+q\geq 2$                & $4pq$                  & $\min(p,q)$        &   \begin{tabular}{@{}c@{}}$5$ if $p=q=2$ \\ $2(p+q)-\frac{5}{2}$ else\end{tabular}   \\ \hline
\multicolumn{1}{|c|}{\multirow{2}{*}{$DIII$}} & $SO(4n)/U(2n), n\geq 1$                             & $2n(2n-1)$             & $n$                &  \begin{tabular}{@{}c@{}}$n\left(n-\frac{1}{2}\right)$ if $n\leq 3$ \\ $4n-\frac{7}{2}$ if $n>3$\end{tabular}   \\ \cline{2-5} 
\multicolumn{1}{|c|}{}                        & $SO(4n+2)/U(2n+1),n\geq 1$                         & $2n(2n+1)$             & $n$                &   $4n-\frac{3}{2}$  \\ \hline
\multicolumn{2}{|c|}{$EI$}                                                                 & $42$                   & $6$                &  $8$   \\ \hline
\multicolumn{2}{|c|}{$EII$}                                                                & $40$                   & $4$                &   $\frac{21}{2}$  \\ \hline
\multicolumn{2}{|c|}{$EIII$}                                                               & $32$                   & $2$                & $\frac{21}{2}$    \\ \hline
\multicolumn{2}{|c|}{$EIV$}                                                                & $26$                   & $2$                &  $8$   \\ \hline
\multicolumn{2}{|c|}{$EV$}                                                                 & $70$                   & $7$                &   $\frac{27}{2}$  \\ \hline
\multicolumn{2}{|c|}{$EVI$}                                                                & $64$                   & $4$                &   $\frac{33}{2}$  \\ \hline
\multicolumn{2}{|c|}{$EVII$}                                                               & $54$                   & $3$                &  $\frac{27}{2}$   \\ \hline
\multicolumn{2}{|c|}{$EVIII$}                                                              & $128$                  & $8$                &  $\frac{57}{2}$   \\ \hline
\multicolumn{2}{|c|}{$EIX$}                                                                & $112$                  & $4$                &   $\frac{57}{2}$  \\ \hline
\multicolumn{2}{|c|}{$FI$}                                                                 & $28$                   & $4$                &   $\frac{15}{2}$ \\ \hline
\multicolumn{2}{|c|}{$FII$}                                                                & $16$                   & $1$                &   $\frac{15}{2}$  \\ \hline
\multicolumn{2}{|c|}{$G$}                                                                  & $8$                    & $2$                &   $\frac{5}{2}$  \\ \hline
\end{tabular}}

\end{table}

Furthermore, if $M=M_1\times M_2$, the root system of $M$ is the direct sum of the root systems of $M_i$, and so $r(M)=\min(r(M_1),r(M_2))$. Thus, Table \ref{tab:symtypegroup} and \ref{tab:valuesconj} are sufficient to compute the value of $r$ for any symmetric space of compact type. Also, this show that the conjecture agrees with the fact that $r_{opt}(M)=\min(r_{opt}(M_i))$.

Finally, assume that $M=(G\times G)/\Delta(G)$. We saw in Subsection \ref{caracteres} that if $\Phi^+$ is a choice of positive roots with basis $\alpha_1,\cdots,\alpha_\ell$, we have $$r_{opt}(M)=\gamma = \underset{1\leq i\leq \ell}{\min} \vert \{ \alpha \in \Phi^+ \vert n_i(\alpha)\geq 1\}\vert.$$We also saw that a choice of positive roots for $M$ was given by $\Tilde{\alpha}:(H,-H)\mapsto \alpha(H)$ with $\alpha\in \Phi^+$, and that $m({\Tilde{\alpha}})=2$ for any $\alpha\in \Phi^+$. The bijection $\alpha\mapsto \Tilde{\alpha}$ extends to a map $\mathfrak{b}_\C^*\mapsto \mathfrak{a}_\C^*$, which sends $\Lambda_G$ to $\Lambda$ and such that $\langle\Tilde{\lambda},\Tilde{\mu}\rangle=\langle\lambda,\mu\rangle$. Thus, from this we get that for any $\lambda\in \Lambda_G$, $$\frac{1}{2}S_{\Tilde{\lambda}}=\vert \{\alpha\in \Phi^+ \vert \langle\alpha,\lambda\rangle\neq 0\}\vert.$$
Let $\pi_i,1\leq i\leq \ell$ be the fundamental weights of the root system of $G$, defined by $$\frac{2\langle\pi_i,\alpha_j\rangle}{\langle\alpha_j,\alpha_j\rangle}=\delta_{i,j}.$$We have that $\lambda\in \Lambda_G$ if and only if $\lambda=\sum_{i=1}^\ell m_i(\lambda)\pi_i$, with $m_i(\lambda)\in \N$. Then, for any $\lambda\in \Lambda_G,\alpha\in \Phi^+$, $$\langle\lambda,\alpha\rangle=\sum_{i=1}^\ell\sum_{j=1}^\ell m_i(\lambda)n_i(\alpha)\langle\pi_i,\alpha_j\rangle=\frac{1}{2}\sum_{i=1}^{\ell} m_i(\lambda)n_i(\alpha)\langle\alpha_i,\alpha_i\rangle.$$If $\lambda\neq 0$, there is $i$ such that $m_i(\lambda)\neq 0$. Then for any $\alpha$ such that $n_i(\alpha)\neq 0$, we get $\langle\lambda,\alpha\rangle>0$, so $$\frac{1}{2}S_{\Tilde{\lambda}}\geq \vert \{ \alpha \in \Phi^+ \vert n_i(\alpha)\geq 1\}\vert\geq \gamma.$$But on the other hand, $\frac{1}{2}S_{n\pi_{i_0}}=\gamma$ for all $n\in \N^*$ where $i_0$ is such that\linebreak $\gamma=\vert \{ \alpha \in \Phi^+ \vert n_{i_0}(\alpha)\geq 1\}\vert$. Thus, we get $$\gamma=\underset{\Tilde{\lambda}\in \Lambda\setminus \{0\}}{\inf} \frac{1}{2}S_{\Tilde{\lambda}},$$so the regularity found in this case fits again the conjecture.

\section{Regularity of K-finite matrix coefficients}\label{kfinitesection}
In this section, we keep the notations introduced in the previous section. Thus $(G,K)$ is a symmetric compact Gelfand pair with $G$ connected such that $G/K$ is simply connected. Recall that $Q$ is an open subset of $\mathfrak{a}$ such that, by Proposition \ref{kakcpt}, $G=K\exp(\overline{Q})K$. Set $G_1=K\exp(Q)K$. Then $G_1$ is a dense open subset of $G$.

\begin{definition}
Let $\pi$ be a unitary representation of $G$ on $\mathcal{H}$ and $(\rho,V)$ a representation of $K$. We say that $\xi\in\mathcal{H}$ is \begin{itemize}
    \item $K$-finite if $\vspan(\pi(K)\xi)$ is finite dimensional,
    \item of $K$-type $V$ if $\vspan(\pi(K)\xi)\simeq V$ as a representation of $K$.
\end{itemize}
\end{definition}
Note that this definition of $K$-type $V$ is not standard.

The goal of this section is to prove the following result.
\begin{theorem}
 \label{kfinitecoef}Assume that for any $K$-bi-invariant matrix coefficient $\varphi$ of a unitary representation of $G$, the function $\psi=\varphi\circ\exp$ is of class $C^{(r,\delta)}$ on $Q$. Let $\pi$ be a unitary representation of $G$ on $\mathcal{H}$, $\xi,\eta\in \mathcal{H}$ $K$-finite. Consider the associated matrix coefficient $c:g\mapsto \langle\pi(g)\xi,\eta\rangle$. Then $c$ is of class $C^{(r,\delta)}$ on $G_1$.   
\end{theorem}

\begin{remark}
This shows that $r_{opt}(M)$, which is the optimal uniform regularity of spherical functions, and also (by Lemma \ref{roptid}) the optimal regularity of $K$-bi-invariant matrix coefficients of unitary representations of $G$, is even the optimal regularity of $K$-finite matrix coefficients.
\end{remark}

Since $K$-invariant vectors are $K$-finite, the first step is to extend the regularity of $\psi$ on $Q$ to a regularity of $\varphi$ on $G_1$. For this, we need to study properties of the decomposition of Proposition \ref{kakcpt}.

\begin{lem}
\label{kaksubmersion}The map $$\fonction{q}{K\times K\times Q}{G_1}{(k_1,k_2,H)}{k_1\exp(H)k_2^{-1}}$$is a submersion.
\end{lem}
\begin{proof}If $g\in G$, denote $L_g$ and $R_g$ the translations by $g$ on the left and right respectively. Let $m:G\times G\to G$ be the multiplication map, its tangent map at $(a,b)$ is $$\fonction{T_{(a,b)}m}{T_aG\times T_bG}{T_{ab}G}{(X_a,X_b)}{T_aR_b(X_a)+T_bL_a(X_b)}.$$We can identify $T_gG$ with $\mathfrak{g}$ by the isomorphism $T_eL_g$. Under this identification, we have $\forall g,h\in G$, $T_hL_g=\Id$ and $T_hR_g=\Ad(g^{-1})$, so that the tangent map of the multiplication becomes $T_{(a,b)}m(X_a,X_b)=\Ad(b^{-1})(X_a)+X_b$. Furthermore, if $k\in K$, since $L_k(K)=K$, $T_kK\subset T_kG$ is identified with $\mathfrak{k}\subset \mathfrak{g}$. Thus by the chain rule we have $$\fonction{T_{(k_1,k_2,H)}q}{\mathfrak{k}\times \mathfrak{k}\times \mathfrak{a}}{\mathfrak{g}}{(X_1,X_2,Y)}{\Ad(k_2)(\Ad(\exp(-H))(X_1)+T_H\exp(Y))-X_2}.$$We know that $\Ad(k)$ is an isomorphism of $\mathfrak{g}$ and an isomorphism of $\mathfrak{k}$ in restriction. Furthermore, $T_H\exp:\mathfrak{a}\mapsto \mathfrak{a}$ is also an isomorphism. Thus, the map $T_{(k_1,k_2,H)}q$ is surjective if and only if $u=\Ad(k_2^{-1})\circ T_{(k_1,k_2,H)}q\circ (\Id,\Ad(k_2),(T_H\exp)^{-1})$ is surjective. We have $$u(X_1,X_2,Y)=\Ad(\exp(-H))(X_1)-X_2+Y.$$

Consider the decomposition $\mathfrak{g}=\mathfrak{k}\oplus\mathfrak{p}$ in eigenspaces of $\sigma$. In this decomposition, $\mathfrak{a}$ is a maximal abelian subspace of $\mathfrak{p}$. Let $\mathfrak{m}=\mathfrak{k}^\mathfrak{a}$, $\Sigma_{\mathfrak{a}}$ the root system of $(\mathfrak{g}_\C,\mathfrak{a}_\C)$.

For $\lambda\in \Sigma_\mathfrak{a}$, let $\mathfrak{g}_\lambda=\{X\in \mathfrak{g}_\C \vert \forall H\in \mathfrak{a}_\C, [H,X]=\lambda(H)X\}$. We have $$\mathfrak{g}_\C=\mathfrak{a}_\C\oplus \mathfrak{m}_{\C}\oplus \bigoplus_{\lambda\in \sigma_{\mathfrak{a}}} \mathfrak{g}_\lambda.$$Let also $\mathfrak{k}_\lambda=\mathfrak{k}\cap (\mathfrak{g}_\lambda \oplus \mathfrak{g}_{-\lambda})$ and $\mathfrak{p}_\lambda=\mathfrak{p}\cap (\mathfrak{g}_\lambda \oplus \mathfrak{g}_{-\lambda})$. From \cite[Ch. VI, Prop. 1.4]{loos1969symmetric2}, we get $$\mathfrak{k}=\mathfrak{m}\oplus \bigoplus_{\lambda\in \Sigma_{\mathfrak{a}}^+} \mathfrak{k}_\lambda=\mathfrak{m}\oplus \mathfrak{l},$$ $$\mathfrak{p}=\mathfrak{a}\oplus \bigoplus_{\lambda\in \Sigma_{\mathfrak{a}}^+} \mathfrak{p}_\lambda=\mathfrak{a}\oplus \mathfrak{b}.$$We also get that for $\lambda\in \Sigma^+_\mathfrak{a}$, there is $Z_{\lambda,1},\cdots,Z_{\lambda,r_\lambda}$ a $\C$-basis of $\mathfrak{g_\lambda}$, such that setting $Z_{\lambda,i}^{+}=Z_{\lambda,i}+\sigma(Z_{\lambda,i})$ and $Z_{\lambda,i}^-=i(Z_{\lambda,i}-\sigma(Z_{\lambda,i}))$, $\{Z_{\lambda,i}^{+}\}$ gives an $\R$-basis of $\mathfrak{k}_\lambda$ and $\{Z_{\lambda,i}^{-}\}$ gives an $\R$-basis of $\mathfrak{p}_\lambda$.

Let also $H_1,\cdots,H_\ell$ be a basis of $\mathfrak{a}$ and $Y_1,\cdots,Y_r$ a basis of $\mathfrak{m}$. Then for $H\in \mathfrak{a}$, we have $[H,Y_i]=0$, $[H,Z_{\lambda,i}^{+}]=-i\lambda(H)Z_{\lambda,i}^{-}$ and $[H,Z_{\lambda,i}^{-}]=i\lambda(H)Z_{\lambda,i}^{+}$. Thus, we see that \begin{itemize}
    \item $u(0,0,H_i)=H_i$,
    \item $u(Y_i,0,0)=e^{-\ad(H)}(Y_i)=Y_i$,
    \item $u(0,Y_i,0)=-Y_i$,
    \item $u(Z_{\lambda,i}^{+},0,0)=e^{-\ad(H)}(Z_{\lambda,i}^{+})=\sin(i\lambda(H))Z_{\lambda,i}^{-}+\cos(i\lambda(H))Z_{\lambda,i}^{+}$,
    \item $u(0,Z_{\lambda,i}^{+},0)=-Z_{\lambda,i}^{+}$.
\end{itemize} 
Since $H\in Q$, $\lambda(H)\not\in i\pi\Z$ so $\sin(i\lambda(H))\neq 0$ and $u$ is indeed surjective.
\end{proof}

The following proposition is found in \cite[Ch.V, Thm 3.3]{borel1998semisimple}.
\begin{prop}
\label{centralisateurconnexe}Let $G$ be a compact, connected, simply connected Lie group and $f$ an automorphism of $G$. Then the set of fixed points of $f$ is connected.
\end{prop}
\begin{remark}
This result implies that the subgroup $K$ is automatically connected if $G$ is simply connected.
\end{remark}

\begin{lem}
\label{egalitestab}Let $M=Z_K(\exp \overline{Q})=\{k\in K\vert \forall a\in \exp\overline{Q}, ka=ak\}$. Consider the action of $K\times K$ on $G$ by $(k_1,k_2).g=k_1gk_2^{-1}$. Then for any $H\in Q$, we have $$\Stab(\exp H)=\{(k,k)\vert k\in M\}=\Delta(M).$$
\end{lem}
\begin{proof}First, we assume that $G$ is simply connected. Let $H\in 2Q$. We have that $$Z_G(\exp H)=\{g\in G \vert g\exp (H)=\exp (H)g\}=\{g\in G \vert c_{\exp(-H)}(g)=g\}$$is the set of fixed points of $c_{\exp(-H)}$. Thus, $Z_G(\exp H)$ is connected by Proposition \ref{centralisateurconnexe}.

Furthermore, \begin{align*}
    \mathrm{Lie}(Z_G(\exp H)) & =  \{X\in \mathfrak{g} \vert \forall t\in \R, \exp(tX)\in Z_G(\exp H) \}  \\
     & =  \{X\in \mathfrak{g} \vert \forall t\in \R, c_{\exp(-H)}(\exp(tX))=\exp(tX)\}\\
     & =  \{X\in \mathfrak{g}\vert \Ad({\exp(-H)})(X)=X\}\\
     & =  \ker(f_H)
\end{align*}where $f_H(X)=\Ad(\exp(-H))(X)-X$.

Consider the basis $\{Y_i\}\cup \{H_i\}\cup \{Z_{\lambda,i}^+\}\cup \{Z_{\lambda,_i}^-\}$ of $\mathfrak{g}$ introduced in Lemma \ref{kaksubmersion}.

Then, we have \begin{itemize}
    \item $f_H(Y_i)=0$,
    \item $f_H(H_i)=0$,
    \item $f_H(Z_{\lambda,i}^+)=\sin(i\lambda(H))Z_{\lambda,i}^-+(\cos(i\lambda(H))-1)Z_{\lambda,i}^+$,
    \item $f_H(Z_{\lambda,i}^-)=-\sin(i\lambda(H))Z_{\lambda,i}^++(\cos(i\lambda(H))-1)Z_{\lambda,i}^-$.
\end{itemize}
Since $H\in 2Q$, for any root $\lambda$, $\lambda(H)\not\in 2i\pi\Z$, and so $\ker(f_H)=\mathfrak{m}\oplus \mathfrak{a}$.

Thus, $\mathrm{Lie}(Z_G(\exp H))$ does not depend on $H\in 2Q$, and since $Z_G(\exp H)$ is connected, $Z_G(\exp H)$ and $Z_K(\exp H)=Z_G(\exp H)\cap K$ do not depend on $H\in 2Q$, and so $Z_K(\exp H)=M$.

Now take $a=\exp(H)\in Q$, and $(k,k')\in \Stab(a)$, then $ka=ak'$. The automorphism $\sigma$ of $G$ is such that $K=G^\sigma$, and $\sigma(a)=a^{-1}$, so we get $ka^{-1}=a^{-1}k'$, thus $ka^2=ak'a=a^2k$. So $k\in Z_K(a^2)$, but $a^2=\exp(2H)$, $2H\in 2Q$ so $Z_K(a^2)=Z_K(a)=M$. Thus $ka=ak=ak'$, so $k=k'$ and $(k,k')\in \Delta(M)$.

The other inclusion is clear, thus $\Stab(a)=\Delta(M)$.

For the general case, since $G/K$ is assumed to be simply connected, by the Remark \ref{cartanhelgremark}, we have $p:\Tilde{G}\twoheadrightarrow G$ the universal cover such that $\ker p\subset Z(\Tilde{G})^\sigma$ and $\Tilde{K}=\Tilde{G}^\sigma$. Then the previous case gives that $\Stab_{\Tilde{K}\times \Tilde{K}}(\exp_{\Tilde{G}}(H))=\Delta(\Tilde{M})$ for any $H\in Q$. Clearly, if $(k_1,k_2)\in \Stab_{\Tilde{K}\times \Tilde{K}}(\exp_{\Tilde{G}}(H))$, then the projection $(p(k_1),p(k_2))\in \Stab_{K\times K}(\exp_G(H))$.

Conversely, let $(k_1,k_2)\subset \Stab_{K\times K}(\exp_G(H))$. There exists $\Tilde{k}_i\in \Tilde{K}$ such that $p(\Tilde{k}_i)=k_i$. Then $k_1\exp_G(H)k_2^{-1}=\exp_G(H)$ implies that there exists $x\in \ker p$ such that $$\Tilde{k}_1\exp_{\Tilde{G}}(H)\Tilde{k}_2^{-1}x=\exp_{\Tilde{G}}(H).$$Thus, $(\Tilde{k}_2,x^{-1}\Tilde{k}_2)\in \Delta(\Tilde{M})$ so $\Tilde{k}_2=x^{-1}\Tilde{k}_2$, thus $k_1=k_2\in Z_K(\exp_G(H))$.

So this tells us that $\Stab_{K\times K}(\exp_G(H))=p(\Stab_{\Tilde{K}\times \Tilde{K}}(\exp_{\Tilde{G}}(H)))$ does not depend on $H\in Q$ and is equal to $\Delta(Z_K(\exp_G(H)))$ so $Z_K(\exp_G(H))=M$ and $\Stab_{K\times K}(\exp_G(H))=\Delta(M)$.
\end{proof}

The following proposition is a refinement of Proposition \ref{kakcpt}.
\begin{prop}
\label{regkak}For any $g\in G$, there exists a decomposition $$g=k_1(g)\exp (H(g))k_2(g)^{-1}$$where $k_1(g),k_2(g)\in K$ and $H(g)\in \mathfrak{a}$. The map $g\mapsto H(g)$ is smooth on $G_1$. Furthermore, for each $g\in G_1$, there exists a neighborhood $U_g$ of $g$ in $G_1$ and a choice of $g\mapsto k_i(g)$ such that $k_i$ is smooth on $U_g$, $i=1,2$.
\end{prop}
\begin{proof}By Lemma \ref{egalitestab}, the map $$\fonction{\Tilde{q}}{(K\times K)/\Delta(M)\times Q}{G_1}{((k_1,k_2) \mathrm{ mod }M,H)}{k_1\exp(H)k_2^{-1}}$$is a well-defined smooth bijection between manifolds of the same dimension.

Let $p:K\times K\to (K\times K)/\Delta(M)$ be the projection. It is a surjective submersion. Let $q$ be the submersion defined in Lemma \ref{kaksubmersion}. We have $q=\Tilde{q}\circ (p,\Id)$. Thus, for any $(x,H)\in (K\times K)/\Delta(M)\times Q$, we have $T_{(x,H)}\Tilde{q}$ surjective. But it is a linear map between vector spaces of the same dimension, so it is invertible. Thus, by the local inversion theorem and since $\Tilde{q}$ is bijective, $\Tilde{q}$ is a smooth diffeomorphism.

Let $(x,H):G_1\to (K\times K)/\Delta(M)\times Q$ be a smooth inverse. We get that $H$ is a smooth map. From \cite[Proposition 4.26]{lee2003introduction}, since $p$ is a submersion, any $(k_1,k_2)\in K\times K$ is in the image of a smooth local section of $p$. Let $g\in G_1$. Since $p$ is surjective, $x(g)=p(k_1,k_2)$. There exists a neighborhood $V$ of $x(g)$ and a smooth section $s=(s_1,s_2):V\mapsto K\times K$ such that $s(x(g))=(k_1,k_2)$.

Let $U=x^{-1}(V)$ neighborhood of $g$. Then $k_i=s_i\circ x$ is smooth on $U$ and $g=k_1(g)\exp(H(g))k_2(g)^{-1}$.
\end{proof}

\begin{coro}\label{regmatrixcoefG1}Let $\varphi$ be a $K$-bi-invariant function on $G$. Then $\varphi\in C^{(r,\delta)}(G_1)$ if and only if $\varphi\circ \exp\in C^{(r,\delta)}(Q)$.
\end{coro}
\begin{proof}Since $\exp$ is smooth, the first implication is clear. For the converse, assume $\psi=\varphi\circ \exp\vert_Q \in C^{(r,\delta)}(Q)$. By the previous proposition, the map $H$ is smooth on $G_1$ and $\varphi=\psi\circ H$ by $K$-bi-invariance, thus $\varphi\in C^{(r,\delta)}(G_1)$ by Lemma \ref{precomposition}.
\end{proof}

Let $\pi$ be a unitary representation of $G$ on $\mathcal{H}$ and $\xi,\eta\in\mathcal{H}$ of $K$-type $V,W$ respectively, for $V,W$ irreducible representations of $K$. Denote $V_\xi=\vspan(\pi(K)\xi)$. Then there is an isomorphism $i_\xi:V\to V_\xi\subset \mathcal{H}$, denote $\xi_0=i_\xi^{-1}(\xi)$. Similarly, define $V_\eta$ and $i_\eta$. Then the map $$\fonction{f}{B(\mathcal{H})}{L(V,W^*)\simeq V^*\otimes W^*}{A}{i_\eta^*Ai_\xi}$$is $K\times K$ equivariant.

For the associated matrix coefficient, we have $c(g)=\langle\pi(g)\xi,\eta\rangle=\langle f(\pi(g))\xi_0,\eta_0\rangle$.

Now denote $U=K\times K$, $(\rho,V_\rho)$ the irreducible representation of $U$ on $V^*\otimes W^*$ and $\lambda$ the regular representation of $G$ on $L^2(G)$.

If $g\in G$, let also $U_g$ be the stabiliser of $g$ for the left-right action of $U$, $V_g\subset V_\rho$ the space of $\rho(U_g)$-invariant vectors and $P_g$ the orthogonal projection on $V_g$.

The $U$-equivariance of $f$ means that for any $(k,k')\in U$ and $A\in B(\mathcal{H})$, we have $$f(\pi(k)A\pi(k')^{-1})=\rho(k,k')(f(A)).$$Furthermore, there are $v_1,\cdots,v_n\in V_\rho$ and $\xi_1,\cdots,\xi_n,\eta_1,\cdots,\eta_n\in\mathcal{H}$ such that$$f(A)=\sum_{i=1}^n \langle A\xi_i,\eta_i\rangle v_i.$$

Then in this setting, we can apply \cite[Lemma 2.2]{AIF_2016__66_5_1859_0} to get the following result.
\begin{lem}
\label{dlmdls}For any $g_0\in G$, there exists a smooth map $\psi:G\to B(V_\rho)$ such that:\begin{enumerate}
    \item $\forall u\in U,g\in G$, $\psi(u.g)=\psi(g)\circ \rho(u^{-1})$,
    \item $\forall v_1,v_2\in V_\rho$, $g\mapsto \langle \psi(g)v_1,v_2\rangle$ is a coefficient of $\lambda$,
    \item $\psi(g_0)=P_{g_0}$.
\end{enumerate}
\end{lem}
\begin{remark}The lemma only states that $\psi$ is Lipschitz but the proof shows it is smooth.
\end{remark}

\begin{prop}
\label{regktypeV}If for any $K$-bi-invariant matrix coefficient $\varphi$ of a unitary representation of $G$, the function $\varphi\circ\exp$ is of class $C^{(r,\delta)}$ on $Q$, then the map $f\circ \pi$ is in $C^{(r,\delta)}(G_1)$.
\end{prop}
\begin{proof}
Let $g_0\in G_1$ and $\psi$ given by Lemma \ref{dlmdls}. Let $\Tilde{f}:g\mapsto \psi(g)(f(\pi(g)))$. By the equivariance of $f$ and $(1)$ of Lemma \ref{dlmdls}, we have $$\Tilde{f}(u.g)=\psi(u.g)(f(\pi(u.g)))=\psi(g)\rho(u)^{-1}\rho(u)(f(\pi(g))=\Tilde{f}(g)$$
so $\Tilde{f}$ is a $K$-bi-invariant map.

Let $(e_1,\cdots,e_d)$ be an orthonormal basis of $V_\rho$. By $(2)$ of Lemma \ref{dlmdls} there are $a_{ij},b_{ij}\in L^2(G)$ such that $$\langle\psi(g)v_i,e_j\rangle=\langle\lambda(g)a_{ij},b_{ij}\rangle$$so $\psi(g)v_i=\sum_{j=1}^d \langle\lambda(g)a_{ij},b_{ij}\rangle e_j$ and finally $$\Tilde{f}(g)=\sum_{i=1}^n\sum_{j=1}^d \langle(\lambda\otimes \pi)(g)(a_{ij}\otimes\xi_i),b_{ij}\otimes \xi_i\rangle e_j.$$

Hence, $\Tilde{f}$ is a sum of $K$-bi-invariant matrix coefficients of unitary representations of $G$, so by the hypothesis and Corollary \ref{regmatrixcoefG1}, $\Tilde{f}\in C^{(r,\delta)}(G_1)$.

By Lemma \ref{egalitestab}, if $a\in \exp Q$, we have $U_a=\Delta(M)$. Thus, $V_a=V_0$ is independent of $a\in \exp Q$. If $g=(k_1,k_2).a=k_1ak_2^{-1}$, we have $(k,k')\in U_g$ if an only if $(k_1^{-1}kk_1,k_2^{-1}k'k_2)\in \Delta(M)$ and so $V_g=\rho(k_1,k_2)V_0$.

Let $g_0=k_0a_0k_0^{'-1}$ and $V_1=V_{g_0}$. Since $\psi(g_0)=P_{g_0}$, there is an orthonormal basis adapted to $V_1$ such that $$\psi(g_0)=\begin{pmatrix}\Id&0\\0&0\end{pmatrix}.$$ Furthermore, since $\psi$ is smooth, there is $A_{g_0}$ neighborhood of $g_0$ such that $$\psi(g)=\begin{pmatrix}A(g)&*\\ *&*\end{pmatrix}$$
with $g\mapsto A(g)$ smooth, $A(g)$ invertible for any $g\in A_{g_0}$. Up to restricting $A_{g_0}$, by Proposition \ref{regkak}, we have $g=k_1(g)\exp(H(g))k_2(g)^{-1}=k_1(g)a(g)k_2(g)^{-1}$ with $k_1,k_2$ smooth on $A_{g_0}$.

By the $K$-bi-invariance of $\Tilde{f}$, for any $g\in A_{g_0}$, we have $$\Tilde{f}(g)=\Tilde{f}(a(g))=\Tilde{f}(k_0a(g)k_0^{'-1}).$$
But then $f(\pi(k_0a(g)k_0^{'-1}))\in V_{k_0a(g)k_0^{'-1}}=\rho(k_0,k_0')V_0=V_1$. Set $$\Phi(g)=\rho(k_1(g)k_0^{-1},k_2(g)k_0^{'-1})\begin{pmatrix}
    A(k_0a(g)k_0^{'-1})^{-1} & 0\\0&0
\end{pmatrix},$$ it is a smooth map on $A_{g_0}$ because $A$ is smooth invertible, $k_1,k_2$ are smooth and $\rho$ is a finite dimensional representation of $U$. Since $f(\pi(k_0a(g)k_0^{'-1}))\in V_1$, we have \begin{align*}
    \Phi(g)(\Tilde{f}(g)) & =  \Phi(g)(\Tilde{f}(\pi(k_0a(g)k_0^{'-1})))  \\
     & =  \Phi(g)\psi(g)(f(\pi(k_0a(g)k_0^{'-1})))\\
     & =  \rho(k_1(g)k_0^{-1},k_2(g)k_0^{'-1})(f(\pi(k_0a(g)k_0^{'-1})))\\
     & =  f(\pi(k_1(g)a(g)k_2(g)^{-1}))\\
     & =  f(\pi(g)).
\end{align*}

Now let $B:B(V)\times V\to V$ be the bilinear map sending $(f,v)$ to $f(v)$. We showed that on $A_{g_0}$, $f\circ \pi=B\circ (\Phi,\Tilde{f})$. Since $\Phi$ is smooth on $A_{g_0}$ and $\Tilde{f}\in C^{(r,\delta)}(G_1)$, we get by Leibniz formula that $f\circ \pi\in C^{(r,\delta)}(A_{g_0})$.

So for any $g_0\in G_1$, there exists a neighborhood $A_{g_0}$ such that $f\circ \pi \in C^{(r,\delta)}(A_{g_0})$. Thus, $f\circ \pi \in C^{(r,\delta)}(G_1)$.
\end{proof}

\begin{proof}[Proof of Theorem \ref{kfinitecoef}]If $\xi,\eta$ are of $K$-type $V,W$ respectively, with $V,W$ irreducible representations of $K$, we showed that $c(g)=\langle\pi(g)\xi,\eta\rangle=\langle f(\pi(g))\xi_0,\eta_0\rangle$ and in Proposition \ref{regktypeV} that $f\circ \pi\in C^{(r,\delta)}(G_1)$, thus $c\in C^{(r,\delta)}(G_1)$.

For the general case, if $\xi,\eta$ are $K$-finite, $V_\xi,V_\eta$ are finite dimensional representations of $K$, so they decompose into a finite number of irreducible representations. Thus, $c$ is a finite sum of matrix coefficients of the previous case, so $c\in C^{(r,\delta)}(G_1)$.\end{proof}

\providecommand{\bysame}{\leavevmode\hbox to3em{\hrulefill}\thinspace}
\providecommand{\MR}{\relax\ifhmode\unskip\space\fi MR }
\providecommand{\MRhref}[2]{%
  \href{http://www.ams.org/mathscinet-getitem?mr=#1}{#2}
}
\providecommand{\href}[2]{#2}


\begin{thebibliography}{dLMdlS16}

\bibitem[BFH20]{zimmer2}
A.~Brown, D.~Fisher, and S.~Hurtado, \emph{Zimmer’s conjecture for actions of {SL(m,Z)}}, Inventiones mathematicae \textbf{221} (2020), no.~3, 1001--1060.

\bibitem[BFH21]{zimmer3}
\bysame, \emph{Zimmer's conjecture for non-uniform lattices and escape of mass}, arXiv preprint arXiv:2105.14541 (2021).

\bibitem[BFH22]{zimmer}
\bysame, \emph{{Zimmer's conjecture: Subexponential growth, measure rigidity, and strong property (T)}}, Annals of Mathematics \textbf{196} (2022), no.~3, 891 -- 940.

\bibitem[BL98]{benyamini1998geometric}
Y.~Benyamini and J.~Lindenstrauss, \emph{Geometric {N}onlinear {F}unctional {A}nalysis}, Amer. Math. Soc. Colloq. Publ., no. ptie.~1, Amer. Math. Soc., Providence, 1998.

\bibitem[Bor98]{borel1998semisimple}
A.~Borel, \emph{Semisimple {G}roups and {R}iemannian {S}ymmetric {S}paces}, Texts and Readings in Mathematics, Hindustan Book Agency, 1998.

\bibitem[Bro22]{brown2022lattice}
A.~Brown, \emph{Lattice subgroups acting on manifolds}, Proc. Int. Cong. Math, vol.~5, 2022, pp.~3388--3411.

\bibitem[Cam06]{camporesi_spherical_2006}
R.~Camporesi, \emph{The spherical {Paley}-{Wiener} theorem on the complex {Grassmann} manifolds {SU}(p+q)/{S}({U}(p)×{U}(q))}, Proc. Amer. Math. Soc. \textbf{134} (2006), no.~9, 2649--2659 (en).

\bibitem[Cle88]{clerc}
J.-L. Clerc, \emph{Fonctions sphériques des espaces symétriques compacts}, Trans. Amer. Math. Soc. \textbf{306} (1988), no.~1, 421--431.

\bibitem[CM89]{10.2307/2001308}
M.~Cowling and C.~Meaney, \emph{On a {M}aximal {F}unction on {C}ompact {L}ie {G}roups}, Trans. Amer. Math. Soc. \textbf{315} (1989), no.~2, 811--822.

\bibitem[CN01]{cowling2001uniform}
M.~Cowling and A~Nevo, \emph{Uniform estimates for spherical functions on complex semisimple {L}ie groups}, Geom. Funct. Anal. \textbf{11} (2001), no.~5, 900--932.

\bibitem[CW75]{wolfcahn}
R.S. Cahn and J.A. Wolf, \emph{{Zeta functions and their asymptotic expansions for compact locally symmetric spaces of negative curvature}}, Bulletin of the American Mathematical Society \textbf{81} (1975), no.~6, 1086.

\bibitem[dLdlS15]{delaatdls}
T.~de~Laat and M.~de~la Salle, \emph{{Strong property {(T)} for higher-rank simple {L}ie groups}}, Proceedings of the London Mathematical Society \textbf{111} (2015), no.~4, 936--966.

\bibitem[dLdlS18]{deLaatdelaSalle+2018+49+69}
\bysame, \emph{Approximation properties for noncommutative {L}p-spaces of high rank lattices and nonembeddability of expanders}, Journal für die reine und angewandte Mathematik (Crelles Journal) \textbf{2018} (2018), no.~737, 49--69.

\bibitem[dLMdlS16]{AIF_2016__66_5_1859_0}
T.~de~Laat, M.~Mimura, and M.~de~la Salle, \emph{On strong property {(T)} and fixed point properties for {Lie} groups}, Ann. Inst. Fourier \textbf{66} (2016), no.~5, 1859--1893 (en).

\bibitem[dlS22]{de2022analysis}
M.~de~la Salle, \emph{Analysis on simple lie groups and lattices}, Proc. Int. Cong. Math, vol.~4, 2022, pp.~3166--3188.

\bibitem[Fis22]{fisher2022rigidity}
D.~Fisher, \emph{Rigidity, lattices, and invariant measures beyond homogeneous dynamics}, Proc. Int. Cong. Math, vol.~5, 2022, pp.~3484--3507.

\bibitem[Gor94]{gordon}
L.~Gordon, \emph{A {S}tochastic {A}pproach to the {G}amma {F}unction}, The American Mathematical Monthly \textbf{101} (1994), no.~9, 858--865.

\bibitem[Hal03]{hall2003lie}
B.C. Hall, \emph{Lie {G}roups, {L}ie {A}lgebras, and {R}epresentations: {A}n {E}lementary {I}ntroduction}, Grad. Texts in Math., Springer, 2003.

\bibitem[HC53]{harishchandra}
Harish-Chandra, \emph{Representations of a {S}emisimple {L}ie {G}roup on a {B}anach {S}pace. {I}}, Trans. Amer. Math. Soc. \textbf{75} (1953), no.~2, 185--243.

\bibitem[HdL13]{Haagerup_2013}
U.~Haagerup and T.~de~Laat, \emph{Simple {L}ie groups without the approximation property}, Duke Math. J. \textbf{162} (2013), no.~5, 925--964.

\bibitem[HdL16]{haagerupdelaat2}
U.~Haagerup and T.~de~Laat, \emph{Simple {L}ie groups without the {A}pproximation {P}roperty {II}}, Transactions of the American Mathematical Society \textbf{368} (2016), no.~6, 3777--3809.

\bibitem[Hel79]{helgason1979differential}
S.~Helgason, \emph{Differential {G}eometry, {L}ie {G}roups, and {S}ymmetric {S}paces}, Grad. Stud. Math., Amer. Math. Soc., 1979.

\bibitem[Hel00]{helgason2000groups}
\bysame, \emph{Groups and {G}eometric {A}nalysis: Integral {G}eometry, {I}nvariant {D}ifferential {O}perators, and {S}pherical {F}unctions}, Math. Surveys Monogr., Amer. Math. Soc., 2000.

\bibitem[HS14]{haagschli}
U.~Haagerup and H.~Schlichtkrull, \emph{Inequalities for jacobi polynomials}, The Ramanujan Journal \textbf{33} (2014), no.~2, 227--246.

\bibitem[Kir76]{kirillov1976elements}
A.A. Kirillov, \emph{Elements of the {T}heory of {R}epresentations}, Grundlehren Math. Wiss., Springer, 1976.

\bibitem[Kna02]{knapp2002lie}
A.W. Knapp, \emph{Lie {G}roups {B}eyond an {I}ntroduction}, Progr. Math., Birkh{\"a}user Boston, 2002.

\bibitem[Laf08]{lafforgue}
V.~Lafforgue, \emph{{Un renforcement de la propriété (T)}}, Duke Math. J. \textbf{143} (2008), no.~3, 559 -- 602.

\bibitem[LdlS11]{laffdls}
V.~Lafforgue and M.~de~la Salle, \emph{{Noncommutative $L^{p}$-spaces without the completely bounded approximation property}}, Duke Math. J. \textbf{160} (2011), no.~1, 71 -- 116.

\bibitem[Lee03]{lee2003introduction}
J.M. Lee, \emph{Introduction to {S}mooth {M}anifolds}, Grad. Texts in Math., Springer, 2003.

\bibitem[Lee19]{lee2019introduction}
\bysame, \emph{Introduction to riemannian manifolds}, Graduate Texts in Mathematics, Springer International Publishing, 2019.

\bibitem[Loo69a]{loos1969symmetric}
O.~Loos, \emph{Symmetric {S}paces vol.1 : General {T}heory}, Mathematics Lecture Note Series, W.A. Benjamin, 1969.

\bibitem[Loo69b]{loos1969symmetric2}
\bysame, \emph{Symmetric {S}paces vol.2 : Compact spaces and classification}, Mathematics Lecture Note Series, W. A. Benjamin, 1969.

\bibitem[PRdlS22]{dlSPR}
J.~Parcet, {\'E}.~Ricard, and M.~de~la Salle, \emph{Fourier multipliers in {{\({\text{SL}_n}(\mathbf{R})\)}}}, Duke Math. J. \textbf{171} (2022), no.~6, 1235--1297 (English).

\bibitem[Sze39]{szego1939orthogonal}
G.~Szeg{\H{o}}, \emph{Orthogonal {P}olynomials}, Amer. Math. Soc. Colloq. Publ., American Mathematical Society, 1939.

\bibitem[vD09]{Dijk+2009}
G.~van Dijk, \emph{Introduction to harmonic analysis and generalized {G}elfand pairs}, De Gruyter Stud. Math., De Gruyter, 2009.

\bibitem[Vre76]{vretare_orth}
L.~Vretare, \emph{Elementary spherical functions on symmetric spaces}, Mathematica Scandinavica \textbf{39} (1976), no.~2, 343--358.

\bibitem[Vre84]{vretare}
\bysame, \emph{Formulas for {E}lementary {S}pherical {F}unctions and {G}eneralized {J}acobi {P}olynomials}, SIAM J. Math. Anal. \textbf{15} (1984), no.~4, 805--833.

\end{thebibliography}
\end{document}